\documentclass[11pt]{article}

\usepackage{amsmath,amsmath,amsthm,amscd,amssymb,verbatim,epsf}

\numberwithin{equation}{section}

\usepackage{amsthm, amssymb, amscd}

\newcommand{\arrow}{\rightarrow}

\newtheorem{pro}{Proposition}[section]
\newtheorem{thm}[pro]{Theorem}
\newtheorem{lem}[pro]{Lemma}

\newtheorem{cor}[pro]{Corollary}
\newtheorem{rem}[pro]{Remark}

\input diagrams

\def\G{{\Gamma}}
 \def\d{{\delta}}
 
 \def\e{{\epsilon}}
 
 \def\L{{\Lambda}}
 \def\m{{\mu}}
 \def\n{{\nu}}
 
  \def\O{{\Omega}}
   \def\s{{\sigma}}
 
 \def\a{{\alpha}}
 \def\b{{\beta}}
 \def\p{{\partial}}
 \def\r{{\rho}}
 \def\ra{{\rightarrow}}
 \def\lra{{\longrightarrow}}
 \def\g{{\gamma}}
 \def\D{{\Delta}}

 \def\x{{\xi}}
 \def\c{{\mathbb C}}
 
 \def\z{{\mathbb Z}}
 \def\2{{\mathbb Z_2}}
 
 \def\t{{\tau}}

 \def\da{{\downarrow}}

 \def\sl2{{SL(2,\mathbb C)}}
 
 \def\qed{{\hspace{2mm}{\small $\diamondsuit$}}}

 \def\pf{{\noindent{\bf Proof.\hspace{2mm}}}}
 \def\ni{{\noindent}}

 \def\sl{{{\mbox{\tiny $\L$}}}}

\usepackage{pictexwd,dcpic}
\usepackage{rcs}

\def\ni{\noindent}

\def\lra{\longrightarrow}

\def\d{\delta}
\def\D{\Delta}
\def\a{\alpha}
\def\b{\beta}
\def\t{\tau}
\def\L{\Lambda}
\def\e{\epsilon}

\def\g{\gamma}

\def\s{\sigma}
\begin{document}

\begin{center}
{\Large {\bf Closed  Quasi-Fuchsian Surfaces In  Hyperbolic Knot
Complements} }
\end{center}

\begin{center}

{\large Joseph D. Masters and Xingru Zhang\footnote{{ Partially
supported by NSF grant DMS    0204428.}}}
\end{center}

\begin{abstract}
We show that every hyperbolic knot complement contains a closed
quasi-Fuchsian surface.
\end{abstract}

\ni {\it Keywords:} Hyperbolic manifold, quasi-fuchsian surface,
$\pi_1$-injective surface.

\ni {\it AMS classification:} Primary 57N35, Secondary 57M25

\section{Introduction}

By a {\it knot complement} we mean, in this paper, the complement
of a knot in a connected closed orientable $3$-manifold (which is
not necessarily $S^3$). A knot complement is said to be {\it hyperbolic}
if it admits a complete  hyperbolic metric of finite volume.
 By a \textit{surface} we mean, in this paper,
 the complement of a finite (possibly empty) set of points in
 the interior of  a compact,  connected, orientable $2$-manifold.
By a \textit{surface in a $3$-manifold} $M$, we mean a continuous,
proper map $f:S\ra M$ from a surface $S$ into $M$.
  A surface $f:S\ra M$ in a 3-manifold $M$ is said to be
 \textit{incompressible}
 if $S$ is not a $2$-sphere and the induced homomorphism
 $f^*: \pi_1(S, s)\ra \pi_1(M, f(s))$ is injective for one (and thus for any)
 choice of base  point $s$ in $S$.
 A surface $f:S\ra M$ in a 3-manifold $M$
 is said to be {\it essential} if it is incompressible and
 the map $f: S\ra M$ cannot be homotoped  into a boundary component or an end component
 of $M$.

 Essential surfaces in  hyperbolic knot complements can be divided
into three mutually exclusive geometric types: quasi-Fuchsian
surfaces, geometrically infinite surfaces, and essential surfaces
with accidental parabolics. Now we recall the relevant
terminology. Let $\mathbb H^3$ denote the hyperbolic $3$-space
(always in the upper half space model) and let $S_\infty^2=\c \cup
\{\infty\}$ denote the boundary at infinity,  where $\c$ is
the plane of complex numbers.
 Let $\overline{\mathbb H}^3=\mathbb H^3\cup S_\infty^2$ be the
 compactification of
$\mathbb H^3$, which is topologically a compact $3$-ball.
 The action of every element of the orientation preserving
isometry group $Isom^+(\mathbb{H}^3)$ extends to an action on
$\overline{\mathbb H}^3$. For a discrete subgroup $\G$ of
$Isom^+(\mathbb{H}^3)$, let $\L(\G)$ denote the limit set of $\G$
in $S^2_\infty$ and let $\O(\G)=S^2_\infty-\L(\G)$ denote the
regular set of $\G$ in $S^2_\infty$.
 A discrete, torsion-free
subgroup $\G$ of $Isom^+(\mathbb{H}^3)$
 is called \textit{quasi-Fuchsian} if its limit set $\L(\G)$ in
 $S^2_\infty$ is a Jordan circle and each of the two components of
 $\O(\G)$ is invariant under the action of $\G$.
 In the special case that the
 Jordan circle is a geometric circle, the subgroup is said to be
 \textit{Fuchsian}.

 If $M$ is a hyperbolic knot complement, then its
 fundamental group can be identified as a discrete torsion free subgroup
$\G$ of  $Isom^+(\mathbb H^3)$.
 A surface $f:S\ra M$ in a hyperbolic knot manifold $M$ is said to be
 \newline
 (a) {\it quasi-Fuchsian}
  if it is essential and $f^*(\pi_1(S))$ is a quasi-Fuchsian subgroup
  of $\Gamma \subset Isom^+(\mathbb{H}^3)$; or
  \newline
  (b) {\it geometrically infinite} if
 it is essential and the limit set of  $f^*(\pi_1(S))$
 is the entire $S^2_\infty$; or
 \newline
  (c) {\it essential with accidental parabolics}
 if it is essential and some non-peripheral element of $\pi_1(S)$
 has a  parabolic image
 in $f^*(\pi_1(S))\subset \pi_1(M)\subset Isom^+(\mathbb{H}^3)$.\newline
 A quasi-Fuchsian surface $f:S\ra M$  is further called a \textit{Fuchsian}
 or \textit{totally geodesic} surface if the map
 lifts to a totally geodesic plane in $\mathbb H^3$ with respect to
 the universal covering  $\mathbb H^3\ra M$.
In such case the image group $f^*(\pi_1(S))$ is a Fuchsian
subgroup of $Isom^+(\mathbb{H}^3)$.

Work of Marden (\cite{Mar}), Thurston (\cite{T}) and Bonahon
(\cite{B})
 implies that every essential surface
 falls into one of these categories.  Another consequence
 of their work is that every geometrically infinite
 surface is homotopic to a virtual fiber.
 (It is still an open question whether every
 hyperbolic knot complement is virtually fibered.)
 In particular, if a \textit{closed} essential surface in a hyperbolic
 knot complement has no accidental parabolics, then it is quasi-Fuchsian.

 Examples of quasi-Fuchsian surfaces in hyperbolic knot complements have
 been scarce.
 It was shown in  \cite{CLR} that every hyperbolic
 knot complement contains closed essential surfaces,
 but the surfaces constructed there (via Freedman tubing)
 always contain accidental parabolic elements.  Similarly, the closed
essential surfaces constructed in \cite{O},
\cite{CL1}, \cite{CL2}, \cite{Li}, all contain accidental parabolics.
 It was shown in \cite{Men} that the complement of an alternating
 knot in $S^3$ contains no closed, \textit{embedded} quasi-Fuchsian surface,
 a result which was extended in \cite{A}.
 On the positive side,
 there are well-known examples, such as the figure-eight knot complement,
 which contain closed, totally geodesic surfaces.
 Also, hyperbolic knot complements in $S^3$ which
 contain closed, embedded, quasi-Fuchsian surfaces
 are constructed in \cite{AR}.
 In this paper we prove the following  general existence theorem.

\begin{thm}\label{qfs}
Every hyperbolic knot complement  contains a closed quasi-Fuchsian
surface.
\end{thm}

   A closed quasi-Fuchsian  surface in a hyperbolic knot complement
  $M$ has the nice property
 that it remains essential in all but finitely many Dehn fillings
 of $M$ (see, for example, Theorem 5.3 of \cite{W}).
Theorem \ref{qfs} thus has  the following topological consequence:

\begin{cor}\label{cor}
Suppose that  $M$ is a hyperbolic knot complement. Then $M$
 contains a closed  essential surface which remains essential in all but
finitely many Dehn fillings of $M$.
\end{cor}

 It was first shown in \cite{CL2} and later, by a different method, in \cite{Li}
 that for any hyperbolic knot complement $M$, all but finitely many
 Dehn fillings of $M$ contain a closed essential surface.
 What's new in Corollary \ref{cor} is that for every hyperbolic knot complement
 $M$, there is a {\it single} closed essential surface in $M$
 which survives all but finitely many Dehn fillings of $M$.

 We wish to thank the referee for many helpful comments.

\section{Outline of proof and plan of paper}

Let $M$ be a hyperbolic knot complement, and let $C$ be a geometric
cusp of $M$.  The complement of the interior of $C$ in $M$, which
we denote by $M^-$,  is a compact (connected and orientable)
$3$-manifold whose boundary is a torus. We call $M^-$ the
\textit{truncated knot complement}. The idea is to construct a metrically
complete convex hyperbolic $3$-manifold $Y$ with the following
properties:
\newline
(1)  $Y$ has non-empty boundary;
\newline (2)   there is a local
isometry $f$ from $Y$ into the knot complement $M$, and thus
 an injective homomorphism $f^*$  of $\pi_1(Y)$ into
$\pi_1(M)$ (by Lemma \ref{inj});
\newline (3) $Y$ has a single cusp, $C_0$, such that\\
(i) the fundamental group of $C_0$ is a free  abelian
group of rank two, which injects into the fundamental group of $Y$ under the
inclusion map,\\
 (ii) the image of $\pi_1(C_0)$ under the
map $f^*$ is a finite index subgroup of $\pi_1(C)$, and\\
 (iii)
every Dehn filling of $Y$ along the cusp $C_0$ results a compact
$3$-manifold which is $\p$-irreducible.

Restricting $f$ to any boundary component of $Y$ gives a closed surface in $M$,
 and the above properties imply
 that the surface is quasi-Fuchsian.  The proof of this implication is
 given at the end of Section \ref{pf}.

To construct such a manifold $Y$, we start with a pair
of (non-compact) embedded, quasi-Fuchsian surfaces $S_i$, $i=1,2$,
 in $M$ such that $S_i^-=S_i\cap M^-$, $i=1,2$, are properly embedded
 essential surfaces with  different boundary slopes on $\partial M^-$.
 The existence of
 such a pair of  surfaces follows from work of Culler and Shalen \cite{CS}.
 Let $n_i$ be the number of  components of $\p S_i^-$
 and let $\D$ be the geometric intersection number
 between a component of $\p S_1^-$ and a component of $\p S_2^-$.
 The fundamental group of $S_i$ can be naturally identified with a
fixed quasi-Fuchsian subgroup $\G_i$ of $\G=\pi_1(M)$. The limit
set $\L_i$ of $\G_i$ is a Jordan circle in $S^2_\infty$. Let $H_i$
be the convex hull of $\L_i$ in $\mathbb H^3$, and let $X_i$ be
the $\e$-collared neighborhood of $H_i$ in $\mathbb H^3$ for some
fixed number $\e>0$. Then each of $H_i$ and $X_i$ is a convex
$3$-submanifold of $\mathbb H^3$ invariant under the action of
$\G_i$. Let $Y_i=X_i/\G_i$. Then $Y_i$ is a metrically complete
convex hyperbolic $3$-manifold with a local isometry $f_i$ into
$M$. Topologically $Y_i$ is a product $I$-bundle over $S_i$, i.e.
$Y_i=S_i\times I$. We have the corresponding truncated $I$-bundle
$Y_i^-=S_i^-\times I$. The ``cusp region'' of $Y_i$ has a
standard shape if the geometric cusp $C$ of $M$ is chosen small
enough. In particular, the parabolic boundary
 $\p_p Y_i^- \equiv \p S_i^-\times I$, is a set of $n_i$
standard Euclidean annuli.

\begin{figure}[!ht]
{\epsfxsize=3in \centerline{\epsfbox{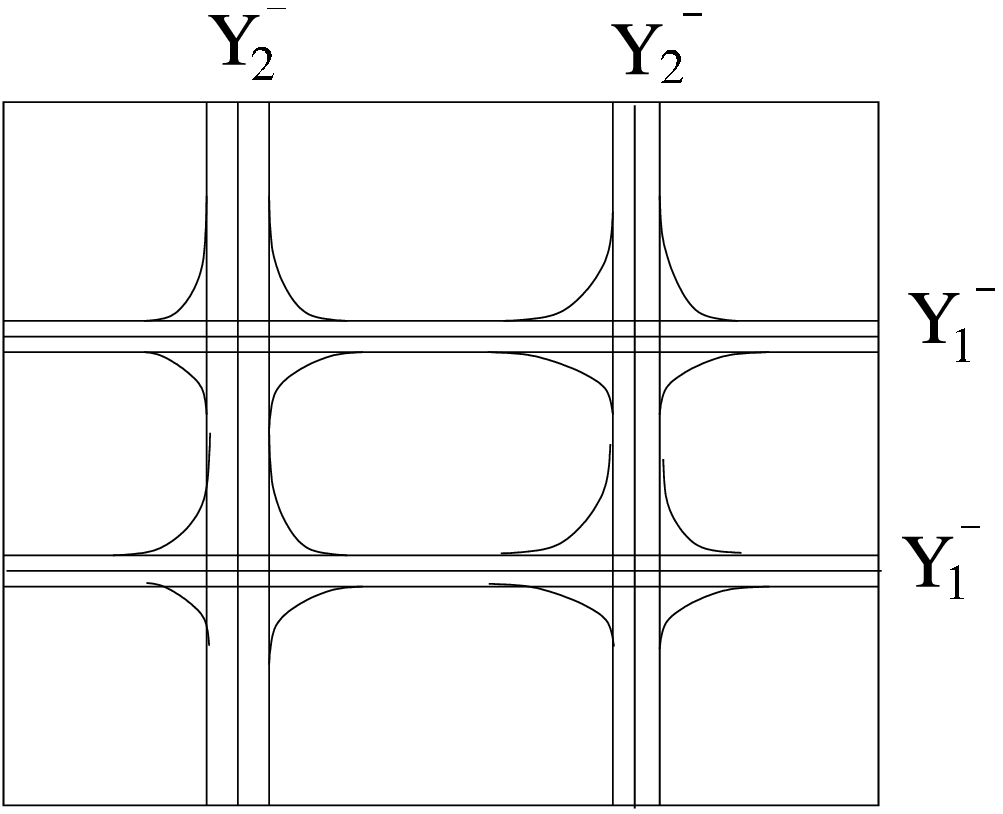}}\hspace{10mm}}
\caption{}\label{intro}
\end{figure}

To illustrate how $Y$ is constructed, let us make some simplifying
 assumptions.  Suppose that each $S_i$ is totally geodesic, that
 $Y_i^-$ is an $\epsilon$-neighborhood of $S_i^-$, embedded in $M$, and
 that $S_1^- \cap S_2^-$ has a large collar neighborhood in both $S_1^-$ and $S_2^-$.
 In this case, we construct $Y$ as an embedded sub-manifold of $M$.
Consider $Y_1^- \cup Y_2^-$, which is a sub-manifold of $M^-$.
 The boundary of this submanifold is convex, except along
 the ``corners'' ($\partial Y_1^- \cap \partial Y_2^-$), and along the truncated
 cusp.  Since we have assumed that the components of $S_1^- \cap S_2^-$
 are well spread out, there is enough room to smooth out the corners,
 as illustrated in Figure \ref{intro} (which shows the part of the smoothing
 near $\p_p Y_1^-\cup\p_p Y_2^-\subset \p M^-$).
 We thus obtain a truncated sub-manifold,
 $Y^- \subset M^-$, whose frontier is convex.
 The complement of $int Y^-$ in $\partial M^-$ consists
 of a finite number of disks, and  the convex hull of each disk
 is a compact subset
 of the cusp $C$. We scoop out each of these convex sets from $C$
 to form a new cusp $C_0$.  The manifold $Y$ is
 the union of $Y^-$ and $C_0$.

 In general, we cannot hope for the manifolds $Y_i^-$ to be embedded
 in $M^-$, and so we must construct $Y$ in a  more abstract way.
 We wish to identify $Y_1^-$ and $Y_2^-$ along
 certain isometric embedded submanifolds
 $K_i^-\subset   Y_i^-$, which correspond to ``intersection"
 components of $Y_1^-$ and $Y_2^-$.
 We then wish to smooth out the corners to form a  hyperbolic 3-manifold $Y^-$
 which is local convex everywhere except on its ``parabolic boundary''
 $\p_p Y^-$. Then we wish to attach a  cusp $C_0$
 along $\p_p Y^-$ to form the required manifold $Y$.

 The gluing and smoothing operations are well-known in the totally
 geodesic case, but to make them work for quasi-Fuchsian surfaces
 is more difficult.
 Furthermore, the gluing can only be performed on manifolds
 with sufficient ``room".  Thus it may be necessary to replace
 the given manifolds $Y_i^-$ with suitable finite covers
 $\breve Y_i^-$.   We construct such covers by proving that
 free groups satisfy a strengthened form
 of the LERF property.

 In Section \ref{property},
 we collect some general facts about hyperbolic geometry.
 Of particular importance is a general fact about convex hulls in
 hyperbolic space (Lemma \ref{convh}), which is essential for
 our gluing constructions.
 In Section \ref{cqf}, we give some general facts about
 cusped, quasi-Fuchsian surfaces, and their convex cores.

 In Section \ref{kiyi}, we construct
 the ``gluing manifolds'' $K_1^-$ and $K_2^-$.
 In the case where each $Y_i^-$ is embedded in $M^-$,
 then $K_i^- \subset Y_i^-$ is just  the intersection $Y_1^- \cap Y_2^-$.
 In general, the fundamental group of each component of $K_i^-$
 is identified with
 the intersection of some conjugate of $\G_1=f^* \pi_1 S_1^-$ and
 some conjugate of $\G_2= f^* \pi_1 S_2^-$,
 and there is an immersion $g_i: K_i^- \rightarrow Y_i^-$.

 The gluing must occur along {\it embedded} sub-manifolds
 of $Y_i^-$, and so we must lift $g_i$ to an embedding.
 For this purpose, it will be useful to isometrically embed $K_i^-$
 into a {\it connected} hyperbolic manifold $J_i^-$,
 whose boundary is convex, outside of a compact set of parabolic
 regions, and which has a local isometry map (still denoted $g_i$) into $Y_i^-$.
 The construction of $J_i^-$ is contained in
 Section \ref{ji}.

 We also wish to control how the parabolic boundary of $J_i^-$ is located in
 $\p Y^-_i$ under the local isometry
 $g_i: J_i^- \ra Y^-_i$,
 and for this purpose we embed $J_i^-$
 isometrically into a certain compact, convex
 hyperbolic manifold $C_n(J_i^-)$.
 We also extend $g_i$ to a local isometry $g_i:C_n(J_i^-) \ra Y_i$.
 The construction of $C_n(J_i^-)$ is contained in Section \ref{cni}.

 Free groups are LERF, and so, using standard arguments,
 it is possible to find a finite cover $\breve Y_i$  of $Y_i$ such that the map
 $g_i:C_n(J_i^-)\ra Y_i$ lifts to an embedding.
 However, for our construction, we require the corresponding truncated cover
 $\breve Y_i^-=\breve S_i^-\times I$ to have
 the same number of parabolic boundary
 components as that of $Y_i^-$. Thus we must show that free groups
 satisfy a strengthened version of the LERF Property.
 This is done in Section \ref{csfg}.
 The proof of this stronger LERF property
  requires much more work than the classical LERF property, and
 may be of independent group theoretic interest.
 The proof applies Stallings' graph-folding techniques.

Finally we need to impose  one more technical condition
 on the covers $\breve Y_i$.
  We require that, after the gluing,
 $\p_p\breve Y_1^-\cup \p_p \breve Y_2^-$
 is isometric to an embedded grid in a
 certain finite  cover $\breve T$ of the Euclidean torus $\p C$.
 The exterior of the grid
 should be a set of Euclidean parallelograms with long sides.
 This requires a further strengthening of the LERF
 property for free groups, which is carried out in Section \ref{wraps}.

 With this property achieved, we can cap off $Y^-$ with a hollowed
 solid cusp $C_0$ along $\p_pY^-$ to get a metrically complete
 convex hyperbolic $3$-manifold $Y$, with non-empty boundary, with a
 single cusp, and with a local isometry into $M$. Thus $Y$
 already has the required properties (1) and (2) given above.
 To show $Y$ has the property (3), we  show  that any Dehn filling
 $Y(\a)$ of  $Y$ with slope $\a$ can be decomposed, in a specified
 way, into handlebody and $I$-bundle pieces.  We call such
 manifold a $HS$-manifold. In Section \ref{hsmf} we show that if an
 $HS$-manifold satisfies certain conditions then its boundary is
 incompressible.

 Our last step is to show that the $HS$-manifold structure of
 $Y(\a)$ satisfies these conditions for incompressibility. The
 final assembly of $Y$, and the proof that $Y$ has all the required
 properties, are given in Section \ref{pf}.

 We remark that Baker and Cooper have recently obtained results on
  gluing convex hyperbolic manifolds (\cite{BC}), which overlap with some of our
  gluing results, for example in Section 4.

\section{Conventions}\label{conv}

In this paper, all  manifolds shall be assumed orientable by
default.  Any $0$-codimension submanifold of an oriented manifold is
given the induced orientation in the obvious way. If $\tilde W$ is a
covering space of an oriented manifold $W$, then the induced
orientation for $\tilde W$ is the one which makes the covering map
orientation preserving. If $W$ is an oriented $n$-manifold ($n\geq
1$) with boundary, then its boundary $\p W$ is given the induced
orientation according the following rule: at each point of $\p W$,
the induced orientation of $\p W$ followed by an inward pointing
tangent vector of $W$ gives the orientation of $W$ at that point.

Let $U_i$ be a submanifold of a manifold $V_i$, $i=1,2$. A map of
pairs $f:(V_1,U_1)\ra (V_2,U_2)$ is called \textit{proper} if the pre-image of
 any compact set is compact, and if
$f(U_1) \subset U_2$.

If $V$ is a hyperbolic $3$-manifold, then for any submanifold $U$ of
$V$ (in particular $\p V$), each component of $U$ is considered as a
metric space with the induced  path metric. If $\tilde V$ is a
connected covering space of $V$, then $\tilde V$ is given the
induced metric so that the covering map from $\tilde V$ to $V$ is a
local isometry.

If $V$ is a metric space and $U$ is a subset of $V$, then $V-U$
denotes the complement of $U$ in $V$,  and $V\setminus U$ denotes
the set obtained by first taking the topological closure of
individual components of
 $V-U$ in $V$ and then taking the disjoint union of these closures.

We say  a connected subspace $U$ of a space $V$ \textit{carries} the
fundamental group of $V$ if the inclusion $U\subset V$ induces a
surjective homomorphism on the fundamental groups.

\section{Some properties of  convex hyperbolic
$3$-manifolds}\label{property}

For standard definitions and  facts about hyperbolic manifolds
(possibly with boundary), the limit set, the convex hull, the
developing map, the holonomy representation, etc.,  we take
\cite{CEG}, \cite{EM} and \cite{Ra} as references.

For any subset $W$ of $\mathbb H^3$, the \textit{limit set} of $W$ in
$S^2_\infty$, denoted $\L(W)$, is the set of intersection points
(possibly empty) between the closure of $W$ in $\overline{\mathbb
H}^3$ and $S^2_\infty$.

Let $V$ be an orientable,  metrically complete, convex (thus
connected) hyperbolic $3$-manifold (possibly with boundary), with
base point $v_0\in V$. Then its universal cover $\tilde V$ is also
a metrically complete, convex, hyperbolic $3$-manifold, and  the
developing map $D:\tilde V\ra \mathbb H^3$ is an isometry of
$\tilde V$ onto its image (Proposition 1.4.2 of \cite{CEG}).
 It follows that the holonomy representation $\r$ of
 $\pi_1(V, v_0)$ into $PSL_2(\c)$ is a discrete and faithful
 representation with no nontrivial elliptic  elements in the image.
 The image group $\G=\r(\pi_1 (V,v_0))$ acts on $D(\tilde V)$
 as a covering transformation group.
So we may consider $\tilde V$ as a submanifold of $\mathbb H^3$
and consider $V$ as the quotient space of $\tilde V$ under the
action of $\G$. Let $p:\tilde V\ra V$ be the quotient map, which
is a universal covering map, and let $\tilde v_0\in\tilde V$ be a
fixed point in $p^{-1}(v_0)$. Then the fundamental group
$\pi_1(V,v_0)$ can be identified with $\G$ in the following way.
Let $\a: ([0,1], \p [0,1])\ra (V,v_0) $ be a loop in $V$, based at
$v_0$, representing a nontrivial  element $\a_*$ of $\pi_1(V,v_0)$,
and let $\tilde \a: ([0, 1], 0)\ra (\tilde V, \tilde v_0)$ be the
unique lift of $\a$ based at the point $\tilde v_0$ with respect
to the covering map $p: (\tilde V, \tilde v_0)\ra (V,v_0)$. Then
the element of $\G$ corresponding to $\a_*$ is the one which maps
$\tilde v_0$ to $\tilde a(1)$.

A nontrivial element $\g$ of $\pi_1(V,v_0)$ is said to be
hyperbolic or parabolic if $\r(\g) \in \G$ is
hyperbolic or parabolic, respectively, in the usual sense;
 i.e. $\g$ has exactly two fixed points or one fixed point,
respectively, in $\overline{\mathbb H}^3$. This definition is
independent of the choices for base points.

Let $V$ be a hyperbolic $3$-manifold and $v_0$ a point in $V$. We
define a \textit{ geodesic loop in $V$ based at $v_0$}
 to be
 a loop $\a: ([0,1], \p [0,1])\ra (V,v_0)$, which is
 geodesic when restricted to $(0,1)$.
 Throughout this paper, a geodesic is always assumed to be
 non-constant.

\begin{lem}\label{georep}
Let $V$ be an orientable, metrically complete,  convex, hyperbolic
$3$-manifold (possibly with boundary), and $v_0\in V$ a base
point. Then every nontrivial element in $\pi_1(V,v_0)$
 is represented uniquely by a   geodesic loop in $V$  based at $v_0$.
\end{lem}

\pf We identify $\tilde{V}$, the universal cover  of $V$,  as a
metrically complete, convex submanifold of $\mathbb H^3$,
 and let $p:\tilde V\ra V$ be the covering map. Fix a point
$\tilde v_0$ in $p^{-1}(v_0)$ as the base point of $\tilde V$. For
a given nontrivial element $\g\in \pi_1(V,v_0)$, let $\a: [0,1]\ra
V$ be a loop in $V$ based at $v_0$ (i.e. $\a(0)=\a(1)=v_0$)
representing $\g$. Let $\tilde \a: [0,1]\ra \tilde V$ be the
unique lift of $\a$ with $\tilde \a(0)=\tilde v_0$. Since $\a$
represents a nontrivial element of $\pi_1(V,v_0)$, $\tilde
\a(0)\ne \tilde \a(1)$. Let $\tilde \s: [0,1]\ra \mathbb H^3$ be
the unique geodesic segment with $\tilde \s(0)=\tilde \a(0)=\tilde
v_0$ and $\tilde \s(1)=\tilde \a(1)$. Since $\tilde V$ is convex,
the geodesic path $\tilde \s$ is contained in $\tilde V$. Thus the
map $\s=p\circ \tilde \s:[0,1]\ra V$ gives a  geodesic loop in $V$
based at $v_0$.
 By convexity, the convex hull of the set $\tilde{\alpha}([0,1]) \cup
 \tilde{\sigma}([0,1])$
 is contained in $\tilde{V}$, and this hull contains a homotopy between
 $\tilde \a$ and $\tilde \s$  with their endpoints fixed.
 Under the covering map $p$, the homotopy
descends  to a homotopy in $V$ between the loop $\a$ and the
geodesic loop $\s$ fixing the base point $v_0$. Hence $\s$ is also
a representative loop of the element $\g$. The uniqueness of such
a based geodesic loop is clear from the argument. \qed

\begin{lem}\label{inj}
Suppose that
 $f:U\ra V$ is a
local isometry between two orientable, metrically complete,
convex, hyperbolic $3$-manifolds $U$ and $V$. Then $f^*:
\pi_1(U,u_0)\ra \pi_1(U, f(u_0))$ is injective for any choice of
the base point $u_0$ in $U$. If in addition $U$ is compact,
then
 $f^*(\pi_1(U,u_0))$ contains no parabolic elements of $\pi_1(V, f(u_0))$.
\end{lem}

\pf Let $v_0=f(u_0)$, let $p:\tilde V\ra V$ be the universal
covering map, where $\tilde V$ is identified as a submanifold of
$\mathbb H^3$, and let $\tilde v_0$ be a fixed point in
$p^{-1}(v_0)$. To prove the first assertion, let $\g$ be a
nontrivial element of $\pi_1(U,u_0)$. By Lemma \ref{georep},
 $\g$ is represented by a   geodesic loop $\s$ in $U$ based at $u_0$.
 Since $f$ is a local
isometry, $f\circ \s$ is a geodesic loop in $V$ based at $v_0$.
 If $f^*(\g)$
is the trivial element of $\pi_1(V, v_0)$, then $f\circ \s$ lifts
to a geodesic loop in $\tilde V$
 based at $\tilde v_0$. But obviously $\mathbb H^3$ contains no
based  geodesic loops.
 Hence $f^*(\g)$ is nontrivial in $\pi_1(V, v_0)$, and
 thus $f^*$ is injective.

Now suppose in addition that $U$ is compact. Let
$H=f^*(\pi_1(U,u_0))$. Let $\bar V=\tilde V/H$, and let $\bar
q:(\tilde V, \tilde v_0) \ra (\bar V, \bar q(\tilde v_0))$, $\bar
p: (\bar V, \bar q(\tilde v_0))\ra (V, v_0)$ be the covering maps.
Since $\bar p^*(\pi_1(\bar V, \bar q(\tilde v_0)))=f^*(\pi_1(U,
u_0))=H$, the map $f:(U, u_0)\ra (V, v_0)$ lifts to a map $\bar f:
(U, u_0)\ra (\bar V, \bar q(\tilde v_0))$. Since $\bar p\circ \bar
f=f$ and since $\bar p$ and $f$ are local isometries, $\bar f$ is
a local isometry.

Let  $p': (\tilde U, \tilde u_0)\ra (U, u_0)$ be the universal
covering map. Then the map $\bar f\circ p'$ lifts to a map $\tilde
f: (\tilde U,\tilde u_0)\ra (\tilde V, \tilde v_0)$. Since $\bar
f\circ p'=\bar q\circ\tilde f$ and since $\bar q, p', \bar f$ are
all local isometries, $\tilde f$ is also a local isometry. Hence
$\tilde f$ sends geodesic arcs to geodesic arcs. Since $\tilde U$
is convex and since $\tilde V$ is a simply connected  submanifold
of $\mathbb H^3$, $\tilde f$ must be an embedding. Since the map
$\tilde f$ is equivariant and  the map $f^*$ is an isomorphism,
from the commutative diagram $$\begin{array}{ccc} (\tilde U,
\tilde u_0)&\stackrel{\tilde f}{\ra} & (\tilde V, \tilde v_0)\\
\da p'& & \da \bar q\\(U,u_0)&\stackrel{\bar f}{\ra} & (\bar V,
\bar q(\tilde v_0)).
 \end{array}$$
we see that $\bar f$ is an embedding. Hence  $\tilde f(\tilde U)$
is a convex submanifold of $\tilde V$ covering the compact
submanifold $\bar f(U)$ of $\bar V$. In fact $\tilde f(\tilde
U)/H=\bar f(U)$.

If $H=f^*(\pi_1(U,u_0))$ contains parabolic elements, then  a
standard hyperbolic geometry argument shows that $\bar f(U)$
contains a non-compact cusp end. In fact if $H_0$ is  a nontrivial
maximal parabolic subgroup of $H$ and if $a\in S_\infty^2$ is the
point fixed by $H_0$, then  there is a horoball $B_a$ of $\mathbb
H^3$, based at $a$, such that $(B_a\cap \tilde f(\tilde U))/H_0$
properly embeds into $\bar f(U)$ as a non-compact end. This is a
contradiction, since $\bar f(U)$ is compact.
 \qed

Every metrically complete, convex subset of $\mathbb H^3$ is a
manifold (Theorem 1.4.3 of \cite{EM}). Obviously  the intersection
of two metrically complete, convex submanifolds
 of a metrically complete, convex $3$-manifold is
 a metrically complete, convex submanifold (when non-empty).
 Every metrically complete, convex $3$-submanifold $U$ of a simply connected,
 metrically complete, convex, hyperbolic $3$-manifold $V$ is
 simply connected (which follows from Lemma \ref{inj}).
A metrically complete, hyperbolic $3$-manifold (possibly with
boundary) is convex or strictly convex if and only if it is everywhere
locally convex or locally strictly convex, respectively
(Corollory 1.3.7 of \cite{CEG}).  These facts will be often used
 in this paper.

Let $V$ be a connected  metric  space and $U$ a subspace of $V$
(possibly disconnected). By an {\it $r$-neighborhood of $U$ in $V$},
denoted $N_{(r,V)}(U)$, we mean the set of points in $V$ whose
distance from $U$ is less than or equal to $r$. Note that the
topology of $N_{(r,V)}(U)$ may be different from that of $U$. An
$r$-neighborhood $N_{(r,V)}(U)$ is further called  an {\it
$r$-collared neighborhood of $U$ in $V$} if, under a universal
covering map $p:\tilde V\ra V$, the components of $p^{-1}(U)$ are
more than  distance $2r$ apart from each other.
 When the ambient space $V$ is clear, we
simply write $N_r(U)$ for $N_{(r,V)}(U)$. The following lemma
follows directly from the definition.

\begin{lem}\label{r-coll}
If $V$ is a simply connected hyperbolic manifold and $U$ a
connected submanifold of $V$, then for any $r>0$, $N_{(r,V)}(U)$
is an $r$-collared neighborhood of $U$ in $V$.\qed
\end{lem}

We also need to define ``$r$-collared neighborhood'' in relative
version, as follows. Let $V$ be a connected, hyperbolic manifold
with boundary and $F$  a submanifold of $\p V$ (possibly  with
infinitely many components). Suppose that $U$ is a submanifold
 of $V$ and let $E=\p U\cap F$
(which possibly has infinitely many components). If there is an $r$-collared
neighborhood $N_{(r,V)}(U)$ of $U$ in $V$ such that for each
component $F_i$ of $F$, $N_{(r,V)}(U)\cap F_i$ is an $r$-collared
neighborhood of $E\cap F_i$ in $F_i$ (where $F_i$ is given the
induced metric as a submanifold of $V$), then we say that the pair
$(U,E)$ has an $r$-collared neighborhood in the pair $(V,F)$. Again
directly from the definition we have the following lemma.

\begin{lem}\label{rr-coll}
Suppose that  $V$ is a simply connected hyperbolic manifold and
$F$ a  submanifold of $\p V$ such that each component of $F$ is
simply connected. Suppose that $U$ is a connected submanifold  of
$V$ and suppose that for each component $F_i$ of $F$,  $F_i\cap \p
U$ is a connected submanifold of $F_i$. Then for any $r>0$, the
pair $(U, \p U\cap F)$ has an $r$-collared neighborhood  in the
pair $(V, F)$. \qed
\end{lem}

For a metrically complete, convex submanifold $V \subset
\mathbb{H}^3$ and a point $v$ in the frontier of $V$ in $\mathbb
H^3$, we use $P_{(v, V)}$ to denote a support plane for $V$ at the
point $v$, i.e. $P_{(v, V)}$ is a hyperbolic plane in $\mathbb H^3$
such that $V$ lies on one side of the plane and such that $V\cap
P_{(v, V)}$ contains the point $v$. A supporting plane always exists
(Lemma 1.4.5 \cite{EM}). Let $\e$ be a fixed  positive number. For a
metrically complete, convex submanifold $V$ in $\mathbb H^3$, the
$\epsilon$-collared neighborhood of $V$ in $\mathbb H^3$, $N_\e(V)$,
is a metrically complete and strictly convex (Lemma 1.4.7 of
\cite{EM}) $3$-dimensional submanifold of $\mathbb H^3$,  with $C^1$
boundary (Lemma 1.3.6
 of \cite{EM}).
 Note that the supporting plane of
$N_\e(V)$ at a point $x$ in the frontier of $N_\e(V)$ (which is
$\p N_\e(V)$ in this case) is unique, and intersects
 $N_\e(V)$ only at the point
$x$, due to the strict convexity of $N_\e(V)$.

 The following
proposition will play a key role.

\begin{pro} \label{convh}For any given
$\epsilon> 0$, there is a number $R=R(\e)>0$ such that the
following holds. If  $V$ and $V'$ are  metrically complete,  convex
submanifolds of $\mathbb{H}^3$ such that $N_{\e}(V)$ and $V'$ have
non-empty intersection, and if $x$ is a point in $\partial
N_{\epsilon}(V)$ such that $d(x, N_\e(V) \cap V') > R$,
 then $P_{(x,N_\e(V))} \cap V' = \emptyset$. In particular
 if we take the
convex hull of the union of $N_\e(V)$ and $N_\e(V')$ then all the
added points are contained in an $R$-collared neighborhood of
$N_\e(V)\cap N_\e(V')$.
\end{pro}

\pf Suppose otherwise that  such $R$ does not exist.
  Let $x\in \p (N_\e(V))$ be a point very far from  $N_{\epsilon}(V) \cap V'$,
 let $A$ be a geodesic segment, tangent to $N_\e(V)$ at $x$,
 contained in the unique supporting plane $P_{(x,N_\e(V))}$, and
 suppose that $A \cap V'$ contains a point $x'$.

\begin{figure}[!ht]
{\epsfxsize=2.5in \centerline{\epsfbox{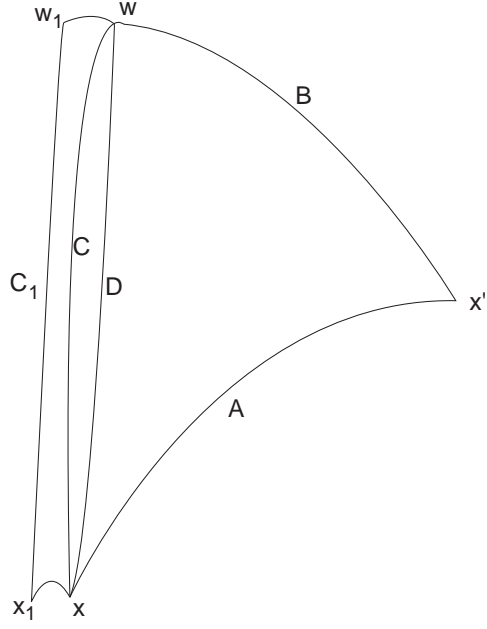}}\hspace{10mm}}
\caption{If $C$ is long, the area between $C$ and $D$ becomes
large}\label{cvn}
\end{figure}

 If $\partial (N_\e(V)) \cap V' = \emptyset$, then
 $V' \subset int N_\e(V)$, and so
 $P_{(x,N_\e(V))} \cap V' = \emptyset$.
 Thus, we may assume that $\partial (N_\e(V)) \cap V'$ contains
 a point $w$. Since every component of $\p N_\e(V)$ separates
 $\mathbb{H}^3$, we may assume that $w$ and $x$ are in the same component
 of $\p N_\e(V)$.
 Let $B$ be a geodesic segment from
$x'$ to $w$, let $C$ be a geodesic segment from $w$ to $x$ , and
let $P_0$ be the unique geodesic plane containing the
 (distinct) points $x, x'$ and $w$. See Figure \ref{cvn}.

Let $x_1$ and $w_1$ be the nearest points in $V$ to $x$ and $w$
respectively, let $C_1$ be a geodesic segment from $x_1$ to $w_1$,
and let $E$ be the geodesic rectangle in $\mathbb{H}^3$ with
vertices $x, x_1, w_1$ and $w$. Since $E$ bounds a surface of area
less than $2\pi$, then if $C$ and $C_1$ are long enough, most of the
arc $C$ is very close to $C_1$; for example, we may assume that
\begin{equation}\label{eq1}
Length(C \cap N_{.01 \e}(C_1)) \geq .99 Length(C).
\end{equation}

Now let $D$ be the segment of the curve $P_0 \cap \partial
(N_\e(V))$ which runs from $x$ to $w$. Since
 $D \subset \p N_\e V$, and $C_1 \subset V$, then
\begin{equation}\label{eq2}
 N_\e C_1 \cap D = \emptyset.
\end{equation}
 By  \ref{eq1} and \ref{eq2}, we have
\begin{eqnarray}\label{eq3}
 Length(C - N_\e D) \geq .99 Length(C).
\end{eqnarray}
 By \ref{eq3} and a simple integration, the area in $P_0$
bounded by $C$ and $D$ is at least $.99 Length(C) * .99 \epsilon$.
But since this region is contained in the triangle region ABC, its
area must be less than $\pi$, which is a contradiction if $C$ is
long enough. \qed

In a similar vein, we have:

\begin{pro} \label{convh2}
Suppose that $X$ is a convex submanifold of $\mathbb{H}^3$,
 that $V_1, ... , V_n$ are convex subsets of $X$, and
 that $N_{(\e,X)}(V_1), ..., N_{(\e,X)} (V_n)$ are all disjoint,
 for some $\e > 0$. Then
 $Hull(N_{(\e,X)}(V_1) \cup ... \cup N_{(\e,X)}( V_n))\setminus
 (N_{(\e,X)}(V_1) \cup ... \cup N_{(\e,X)}( V_n))$
 is compact.
\end{pro}

\pf We first note that for any convex subset $V$ of $X$,
$N_{(\e,X)}(V)=N_\e(V)\cap X$ and that $\p N_{(\e, X)}(V) \cap int
X= \p N_{\e}(V) \cap int X$.

 For every $x$ in $\p N_{(\e, X)}(V_i) \cap int X$,
 there is a geodesic plane $P_x$, such that
 $P_x \cap N_{(\e,X)}(V_i) = \{x \}$.
 Let $\e^{\prime}>0$
 be a number such that $N_{(\e, X)}(V_i) \cap N_{(\e^{\prime}, X)}(V_j) \neq \emptyset$
 for all $1 \leq i,j \leq n$.
 Since $N_{(\e, X)}(V_i)$ and $ N_{(\e, X)}(V_j)$ are disjoint,
the limit set of $N_{(\e, X)}(V_i)$ is disjoint from the limit set
of  $ N_{(\e, X)}(V_j)$. Thus
 $N_{(\e, X)}(V_i)\cap N_{(\e^{\prime},X)}(V_j)$ is
   compact in $\mathbb H^3$. The proof of
 Proposition \ref{convh} shows that there is a compact subset
 $B_i \subset N_{(\e,X)}(V_i)$
  such that   for all
 $x \in (\partial N_{(\e,X)}(V_i)\setminus B_i)\cap int (X)$,
 we have
 $P_x \cap  N_{\e^{\prime}}(V_j)= \emptyset$ for each $j\ne i$.
 It  follows that
 $Hull(N_{(\e, X)}(V_1) \cup ... \cup N_{(\e, X)} (V_n))\setminus
( N_{(\e, X)}(V_1) \cup ... \cup N_{(\e, X)} (V_n))$ has no limit
points in $S^2_\infty$
 and thus is compact.
\qed

Let $\G$ be a  group, $H\subset \G$ a subgroup, and $\g$ an
element in $\G-H$. We say that $H$ is \textit{separable} from $\g$ in $\G$
if there exists a finite index subgroup $G$ of $\G$ such that $G$
contains $H$ but does not contain $\g$. If $H$ is separable from
every element in $\G-H$, then $H$ is said to be separable in
$\G$. It is easy to see that if $H$ is separable in $\G$, then
given any finite set of elements $y_1,..,y_a$
in $\G-H$,  there is a finite index subgroup $G$ of $\G$ such
that $G$ contains $H$ but does not contain any of $y_1,...,y_a$.

\begin{pro}\label{sep1}
Suppose that $U$ is a compact, convex, hyperbolic
$3$-manifold,
 that $V$ is  a metrically complete, convex, hyperbolic
$3$-manifold, and that
 $f: (U,u_0) \ra (V, v_0)$ is a local isometry. Then
 there is a finite (possibly empty) set of
elements  $y_1,...,y_a$ in $\pi_1(V,
v_0)-f^*(\pi_1(U,u_0))$ with the following property:
 if
$G\subset \G=\pi_1(V, v_0)$ is a finite index subgroup which
separates $H=f^*(\pi_1(U, u_0))$ from $y_1,...,y_a$,
 and if $\bar{V}$ is the finite cover of $V$ corresponding
 to $G$, then the map
 $f:(U, u_0)\ra (V, v_0)$ lifts to an embedding
 $\breve f:U \ra \bar{V}$.
\end{pro}

\pf
  Let    $p:(\tilde V,\tilde v_0)\ra  (V, v_0)$ and
   $p': (\tilde U, \tilde u_0)\ra (U,u_0)$ be the universal
covering maps, let $\bar V=\tilde V/H$,  and let $\bar q:(\tilde
V, \tilde v_0) \ra (\bar V, \bar q(\tilde v_0))$
 and  $\bar p: (\bar V, \bar q(\tilde v_0))\ra (V, v_0)$ be the covering
maps. As in the proof of Lemma \ref{inj}, we have the commutative
diagrams
\begin{diagram}[height=2em,w=3em]
 (\tilde U, \tilde u_0)&&\rTo^{\tilde f} &&(\tilde V,
\tilde v_0)\\ \dTo^{p'}& &&& \dTo^{\bar q}\\(U,u_0)&&\rTo^{\bar f}
&&(\bar V, \bar q(\tilde v_0))\\&\rdTo_{f}&&\ldTo_{\bar p}&\\& &(V,
v_0)&&
 \end{diagram}
where both $\tilde f$ and $\bar f$ are embeddings, such that $\tilde
f(\tilde U)$ is a simply connected  convex submanifold of $\tilde V$
covering $\bar f(U)$ with covering group  $H$. Since $\bar f(U)$ is
compact,
 there is  a connected compact submanifold $D$ in
 $\tilde V$ such that
$\bar q(D)$ contains $\bar f(U)$. Since the action of $\G=\pi_1(V,v_0)$ on
$\tilde V$ is properly discontinuous, there are only finitely many elements
$\g$ of $\G$ with $D\cap \g(D)\ne \emptyset$. Let $y_1, ..., y_a$
be all such elements  which are not contained in $H$. Suppose that
$G$ is a finite index subgroup of $\G$ such that $G$ contains $H$
but does not contain any of $y_1,...,y_a$. Let $\breve V=\tilde
V/G$. Then the covering map $\breve q:(\bar V, \bar q(\tilde
v_0))\ra (\breve V, \breve q(\bar q(\tilde v_0)))$ embeds $\bar
f(U)$ into $\breve V$. Let $\breve  p:(\breve V,\breve q(\bar
q(\tilde v_0)))\ra (V,v_0)$ be the finite covering map . One can
easily check that $f=\breve p\circ \breve q\circ \bar f$ (since
$f=\bar p\circ \bar f$ and $\bar p=\breve p\circ \breve q$). Hence
$\breve f=\breve q\circ \bar f:(U,u_0)\ra (\breve V, \breve q(\bar
q(\tilde v_0)))$ is a lift of the map $f:(U, u_0)\ra (V,v_0)$ such
that $\breve f$ is an embedding. \qed

\section{Cusped Quasi-Fuchsian surfaces and their convex
cores}\label{cqf}

Recall that  $M$ denotes an arbitrary fixed, connected, orientable,
complete, finite-volume, hyperbolic 3-manifold with a single cusp.
We consider $M$ as the quotient space of $\mathbb H^3$ under the
action of a fixed, discrete, torsion-free subgroup $\G$ of
$PLS_2(\c)$. A point $a\in S_\infty$  is called a
\textit{parabolic fixed point}
 of a subgroup of $\G$ if $a$ is the fixed point of a
 parabolic element of the subgroup (note that the trivial element
 is not considered as a parabolic element). We may assume that
the point $\infty$ is a parabolic fixed point of $\G$ (up to
replacing $\G$ by
 a conjugate of $\G$ in $PSL_2(\c)$, which we may assume has been done).
The quotient map  $p:\mathbb H^3\ra M=\mathbb H^3/\G$ is a fixed
universal covering map of $M$. Note that $\G$ acts on $\mathbb
H^3$ isometrically as the covering transformation group, and  $p$
is a local isometry. Also $\G$ is isomorphic to the fundamental
group of $M$.

 Let $C$ be an embedded geometric cusp in $M$, i.e.
 ${\cal B}=p^{-1}(C)$ is a set of mutually disjoint
horoballs in $\mathbb{H}^3$ invariant under the action of $\G$.
 Later, we may need to  shrink $C$ if necessary to satisfy some
extra conditions. Note that each component of ${\cal B}$ is based
at a parabolic fixed point of $\G$, and in this fashion the set of
parabolic fixed points of $\G$ is in one-to-one correspondence
 with the set
of components of ${\cal B}$.  Also the set of parabolic fixed points
of $\G$ is invariant under the action of $\G$, and the action is
transitive (since $M$ has a single cusp). Hence all components of
${\cal B}$ are mutually isometric to each other by an element of
$\G$.

   Let $M^-$ be the complement of  the interior of
 $C$ in $M$. By \cite{CS} \cite{CL2}  there are two
connected, embedded, orientable, cusped,  quasi-Fuchsian surfaces
$S_i$ in $M$, such that $S_i^-=S_i\cap M^-$, $i=1,2$,  have
different boundary slopes (we may assume that $S_i\cap \p M^-$  is a
set of embedded simple closed curves each being essential in the
torus $\p M^-$). Let $n_i$ be the number of cusps in $S_i$, i.e.
$n_i$ is the number of components of $\p S_i^-$. By a well-known
duality argument, at least one of the surfaces $S_i$
 must have even number of boundary components, i.e. at least one
 of the integers $n_i$ must be even.

Let $\tilde S_i$ be a fixed component of $p^{-1}(S_i)\subset
\mathbb H^3$ whose closure in $\overline{\mathbb H}^3$ contains
the point $\infty$. Let $Stab_{\G}(\tilde S_i)$ denote the maximal
subgroup of $\G$ which leaves $\tilde S_i$ invariant. Then there
is a finite-index subgroup $\G_i$ of $Stab_{\G}(\tilde S_i)$ such
that $\tilde S_i/\G_i=S_i$ and $\G_i$ is isomorphic to the
fundamental group of $S_i$. As
$\G_i$ is a quasi-Fuchsian subgroup, the limit set $\L_i$ of
$\G_i$ is a
 Jordan  circle in $S^2_\infty$, containing the point   $\infty$ (by our choice).
 Let $H_i$ be the convex hull of $\L_i$ in $\mathbb H^3$.
 Note that $H_i$ is invariant under the action of $\G_i$.

\begin{lem}\label{strip}{\rm \cite{CL2}}
The convex hull $H_i$ lies between two parallel vertical planes in
$\mathbb H^3$. \qed
\end{lem}

The two vertical planes given by Lemma \ref{strip} are
based on two parallel Euclidean lines in $\c$.  Among all pairs of
planes satisfying Lemma \ref{strip}, let $P_{i,j}$, $j=1,2$, be
the pair which are  closest to each other; thus $H_i$
lies between $P_{i,1}$ and $P_{i,2}$, and $P_{i,j}\cap H_i$ is
non-empty for each $j=1,2$. Let $W_i$ be the closed
$3$-dimensional region between the two planes $P_{i,1}$ and
$P_{i,2}$. Let $B_\infty$ be the component of ${\cal B}$ based at
the point $\infty$. So $\p B_\infty$ is a horizontal, Euclidean plane in
$\mathbb H^3$, and $W_i \cap \p B_\infty$ is a \textit{strip}--
 i.e. a region bounded by parallel lines in a Euclidean plane.
 Furthermore, $W_i \cap B_{\infty}$ is the product of the strip
 $W_i \cap \p B_\infty$ with $[0, \infty)$; we call
 this a \textit{3-dimensional strip region, based on }
 $W_i \cap \p B_\infty$.

\begin{lem}\label{strip2}
If the cusp $C$ of $M$ is small enough, or equivalently if the
horizontal  plane $\p B_\infty$ is high enough (i.e. its Euclidean
distance from the complex plane $\c$ is big enough), then $H_i\cap
B_\infty=W_i\cap B_\infty$.
\end{lem}

\pf  Since $H_i$ is convex, we just need to show that if the
horizontal plane $\p B_\infty$ is high enough, then $P_{i,j}\cap
B_\infty$ is contained in $H_i$ for both $j=1,2$. We prove this
for $j=1$; the $j=2$ case being entirely similar. Each of $H_i$,
$W_i$, $P_{i,1}$ is invariant under the action of some parabolic
element $\b_i$ of $\G_i$, which is a horizontal Euclidean
translation. Let $x$ be a point in $H_i\cap P_{i,1}$. Then
$\b_i(x)$ is also contained in $H_i\cap P_{i,1}$, and so is the
hyperbolic geodesic segment $\a$ in $P_{i,1}$ with endpoints $x$
and $\b_i(x)$. Since $\infty$ is a limit point of $H_i$, every
vertical ray in $\mathbb H^3$ based at a point in $H_i$ is
entirely contained in $H_i$. So the part of $P_{i,1}$ lying
directly above $\a$ is contained in $H_i\cap P_{i,1}$. So all the
translations of this set under powers of $\b_i$ are contained in
$H_i\cap P_{i,1}$. So it is clear that if $\p B_\infty$ is higher
than the highest point of the geodesic segment $\a$, then
$P_{i,1}\cap B_\infty$ is contained in $H_i$.
 \qed

Note that the center line of the strip $H_i\cap \p B_\infty$ has
the same slope as that of $\p S_i^-$; that is, its image under the
covering map $p:\mathbb H^3\ra M$ is a simple closed curve in $\p
M^-$ isotopic to a boundary component of $S_i^-$.

 Now let $B_a$ be
any fixed component of ${\cal B}$ based at a parabolic fixed point
$a$ of $\G_i$, and let $\g\in \G$ be any fixed  element which maps
$a$ to $\infty$. Then $\g(B_a)=B_\infty$. Consider the convex set
$\g(H_i)$.  As in Lemma \ref{strip2}, one can show that, after shrinking
 $C$ if necessary,
$\g(H_i)\cap B_\infty$ is a $3$-dimensional strip region, based on a strip
in  $\p B_\infty$. Note that the center line of the strip
$\g(H_i)\cap \p B_\infty$ is parallel to the center line of the
strip $H_i\cap \p B_\infty$, since the boundary curves of $S_i^-$
 are all isotopic in $\p M_i^-$.

\begin{lem}\label{strip3}
Up to replacing the cusp $C$ by a smaller geometric cusp,
$\g(H_i)\cap B_\infty$ is a  $3$-dimensional strip region, for
every $\g\in \G$ which sends a parabolic fixed point of $\G_i$ to
$\infty$. Moreover the center line of  the strip $\g(H_i)\cap \p
B_\infty$ is parallel to the center line of the strip $H_i\cap \p
B_\infty$.
\end{lem}

\pf  The lemma follows from the notes given in the preceding paragraph,
 together with the facts that the set of
parabolic fixed points of $\G_i$ is invariant under the action of
$\G_i$ and that the action has only finitely many orbits (exactly
$n_i$ orbits in fact). \qed

From now on we assume that the cusp $C$ of $M$ has been chosen
small enough so that the conclusion of Lemma \ref{strip3} holds.

 Fix a small positive number $\e$ (e.g. $\e=1$) and let $X_i$
 be the $\e$-collared neighborhood of
 $H_i$ in $\mathbb H^3$ (cf. Lemma \ref{r-coll}).

\begin{cor}\label{strip4}
$\g(X_i)\cap \p B_\infty$ is a strip between two parallel
Euclidean lines in $\p B_\infty$ for every $\g\in \G$ which sends
a parabolic fixed point of $\G_i$ to $\infty$. Moreover the center
line of the strip $\g(X_i)\cap \p B_\infty$ is parallel to the
center line of the strip $X_i\cap \p B_\infty$. \qed
\end{cor}

In fact $\g(X_i)\cap \p B_\infty$ is an $\e$-collared neighborhood
of $\g(H_i)\cap \p B_\infty$ in $\p B_\infty$ for every $\g$ given
in Corollary \ref{strip4}.

 Note that $X_i$  is a   metrically complete and strictly convex
3-submanifold of $\mathbb{H}^3$ with $C^1$ boundary,
 invariant under the action of $\G_i$. Let $${\cal
B}_i=\{X_i\cap B; \mbox{$B$ a component of $ {\cal B}$ based at a
parabolic fixed point of $\G_i$}\}.$$ We call ${\cal B}_i$ the
\textit{horoball region} of $X_i$. Let $X_i^-=X_i\setminus {\cal B}_i$, and
call  $X_i^-\cap \p {\cal B}_i$ the \textit{parabolic boundary} of $X_i^-$,
denoted by $\p_p X_i^-$. Note that $ X_i^-$ is locally
convex everywhere except on its parabolic boundary.

Each of $X_i$,
 ${\cal B}_i$, $X_i^-$ and $\p_p X_i^-$ is invariant under the
action of $\G_i$. Let $Y_i=X_i/\G_i$, which is a metrically
complete and strictly convex hyperbolic $3$-manifold with
boundary. Topologically $Y_i=S_i\times I$,
 where $I = [-1, 1]$. There is a
local isometry $f_i$ of $Y_i$ into $M$, which is induced from the
covering map $\mathbb H^3/\G_i\;\lra \;M$ by restriction on $Y_i$,
 since $Y_i=X_i/\G_i$ is a submanifold of $\mathbb H^3/\G_i$. Also
$p|_{X_i}=f_i\circ p_i$, where $p_i$ is the universal covering map
$X_i\ra Y_i=X_i/\G_i$.
 Let $Y_i^-=X_i^-/\G_i$,  let ${\cal C}_i={\cal B}_i/\G_i$,
and let  $\p_p Y_i^-=\p_p X_i^-/\G_i$. We call
${\cal C}_i$ the cusp part of $Y_i$, and call $\p_pY_i^-$ the
parabolic boundary of $Y_i^-$, which is the frontier of $Y_i^-$ in
$Y_i$ and is also the frontier of ${\cal C}_i$ in $Y_i$. The manifold $Y_i^-$
is locally convex everywhere except on its parabolic boundary.
Topologically $Y_i^-=S_i^-\times I$, where each
component of $\p_p Y_i^-$ is an annulus.

From now on we fix an $I$-bundle structure for $Y_i=S\times I$ as
follows. We first fix an $I$-bundle structure on $Y_i^-=S_i^-\times
I$ such that $\p_p Y_i^-=\p S_i^-\times I$. We may actually assume
that $\p S_i^-\times \{0\}$ are the center horo-circles of $\p_p
Y_i^-$ and that all the $I$-fibers in $\p_pY_i^-$ are perpendicular
to $\p S_i^-\times \{0\}$ with respect to the hyperbolic metric.
Next we extend the $I$-bundle structure to the cusp part ${\cal
C}_i$ of $Y_i$ in the most natural way, i.e. if $C_{i,j}$ is a
component of ${\cal C}_i$ and if we write $C_{i,j}$ as
$A_{i,j}\times [0,\infty)$, where each $A_{i,j}\times \{*\}$ is a
horo-annulus, then we require each $A_{i,j}\times \{*\}$ consists of
$I$-fibers, and  all the $I$-fibers in $A_{i,j}\times \{*\}$ to be
Euclidean geodesics perpendicular to the center horo-circle of
$A_{i,j}\times \{*\}$.

 We let any (free) cover of $Y_i$ have the induced $I$-bundle structure.
 In particular $X_i$ has the induced $I$-bundle structure from that of $Y_i$,
 and this structure is preserved by the action of $\G_i$; i.e. every
 element of $\G_i$ sends an $I$-fiber of $X_i$ to an $I$-fiber of $X_i$.

\begin{lem}\label{length}
For each $i=1,2$, there is a upper bound for the lengths of the
$I$-fibers of $Y_i$.
\end{lem}

\pf Certainly the lengths of the $I$-fibers of
 $Y_i^-=S_i^-\times I$  are bounded,
 since $S_i^-$ is compact. So we only need to show that the
 lengths of the $I$-fibers are bounded in the cusp part ${\cal C}_i$
 of $Y_i$. In turn we just need to show that this is true for
every component of ${\cal C}_i$. Let $C_{i,j}$ be a component of
${\cal C}_i$, and let $\tilde C_{i,j}$ be a component of
$p_i^{-1}(C_{i,j})$. There is an element $\s_{i,j}$ of $\G$ such
that $\s_{i,j}(\tilde C_{i,j})=\s_{i,j}(X_i)\cap B_\infty$. So we
only need to show that the lengths of the $I$-fibers are bounded in
 $\s_{i,j}(X_i)\cap B_{\infty}$. But $\s_{i,j}(X_i)\cap B_\infty$ is the
$\e$-collared neighborhood of $\s_{i,j}(H_i)\cap B_\infty$ in
$B_\infty$ by Lemma \ref{strip3}. Also from Lemma \ref{strip3}, we
see that $\s_{i,j}(H_i)\cap B_\infty$ has the natural $I$-bundle
structure, which is the restriction of the $I$-bundle structure of
$\s_{i,j}(X_i)\cap B_\infty$. Clearly all
 $I$-fibers of $\s_{i,j}(H_i)\cap \p B_\infty$ have the same
 length and every other $I$-fiber of $\s_{i,j}(H_i)\cap B_\infty$
has shorter  length. Similar conclusions hold for $I$-fibers of
$\s_{i,j}(X_i)\cap B_\infty$. \qed

\begin{cor}\label{length2}
For each $i=1,2$, there is a upper bound for the lengths of the
$I$-fibers of $X_i$. \qed
\end{cor}

 The map $f_i:Y_i=S_i\times I\ra M$ is a local isometry but is not an
embedding in general. In particular the center surface
$f_i|:S_i\times \{0\}\ra M$ may not be an embedding, but it
follows from Corollary \ref{strip4} that the map is an embedding
when restricted on each component of $(S_i\times \{0\})\cap {\cal
C}_i$. Hence we may slightly perturb, if necessary, the cusp part of the $
S_i\times \{0\}$ in $Y_i$, keeping it totally geodesic and
transverse to the $I$-fibers, so that the resulting
surface, when restricted to its cusp part,
 will be an embedding under the map $f_i$.
 We still use $S_i$ to denote this surface, and
 we still denote $Y_i$ as $S_i\times I$ and $Y_i^-$ as
$S_i^-\times I$. We call $S_i$ the (topological) center surface of
$Y_i$. Note that $f_i:S_i\ra M$ is quasi-Fuchsian and each
component of $p^{-1}(f_i(S_i))$ is contained in $\g(X_i)$ as a
topological center surface for some $\g\in \G$.

The restriction map  $f_i:(Y_i^-,  \p_p Y_i^-)\ra (M^-,\p M^-)$ is
 a proper map of pairs and  $f_i|:
 (S_i^-,\p S_i^-)\ra (M^-, \p M^-)$ is a
 proper map which is an embedding on $\p S_i^-$ (This property will remain
 valid if we shrink the cusp $C$ of $M$ geometrically).
 We fix an orientation for $S_i$, and let $S_i^-$ and $\p S_i^-$
 have the induced orientation.
 Let $\b_{i,j}$, $j=1,...,n_i$,
denote the components of $\p S_i$ indexed so that their images
$f_i(\b_{i,j})$, $j=1,....,n_i$,  appear consecutively on $\p
M^-$. Let $\D$ be the geometric intersection number between
$f_1(\b_{1,1})$ and $f_2(\b_{2,1})$. Since each $f_i(\b_{i,j})$ is
a Euclidean circle in the Euclidean torus $\p M^-$, each pair of
circles  $f_1(\b_{1,j})$ and $f_2(\b_{2,k})$ have exactly $\D$
intersect points.
 Hence there are a total of $d=n_1n_2\D$
 intersection points between $f_1(\p S_1^-)$
and $f_2(\p S_2^-)$ in the torus $\p M^-$ (all distinct in $\p
M^-$). Let $t_1,...,t_d$ denote these intersection points. The
points $f_i^{-1} \{t_1, ..., t_d\}$ can be indexed as
 $\{t_{i,j,k}, j=1,...,n_i, k= 1, ..., d_i\}$, where $d_i=\D n_{i_*}$
 and $i_*$ is the number such that $\{i,i_*\}=\{1,2\}$.
 We may further assume that
$\{t_{i,j,k}, k = 1, ..., d_i\}$ are contained successively in the
component $\b_{i,j}$, following the orientation of $\b_{i,j}$, for
each $j=1,...,n_i$.

We remark that all the results and notations in this section will
still be valid and consistent if we replace the cusp $C$ by a
smaller one.

\section{The manifold $K_i$}\label{kiyi}

We continue to use the notations established in Section \ref{cqf}.
 The purpose of this section is to construct, for each of $i=1,2$, a manifold $K_i$,
 which, on an intuitive level, corresponds to the intersection
 of $Y_1$ and $Y_2$ in $M$, and which will be used to cut and paste
two immersions.

For each of the points $t_j$, $j=1,...,d$, which was defined at
the end of Section \ref{cqf}, there is a unique embedded  geodesic
ray $R_j$ in $C$, based at $t_j$, perpendicular to $\p C$. We
shall associate to each $R_j$ (thus to $t_j$), a metrically
complete and convex hyperbolic manifold $K_{i,k}$ with a local
isometry, $g_{i,k}$, into $Y_i$ (for each of $i=1,2)$ such that
\newline(1)  the truncated version of $K_{i,k}$, denoted $K_{i,k}^-$ (whose
definition will be given below), is a compact $3$-manifold;
\newline
(2) there is an isometry $h_{k}:K_{1,k}\ra K_{2,k}$ such that
$h_k|: (K_{1,k}^-, \p_pK_{1,k}^-)\ra (K_{2,k}^-,\p_pK_{2,k}^-)$ is
a proper isometry. (The definition of $\p_pK_{i,k}^-$ will be
given below.)

To do this, we first
 choose  points $b_j$, $j=1,...,d$,  in $\p B_\infty$ such that
$p(b_j)=t_j$.  Recall that $p:\mathbb H^3\ra M$ and $p_i:X_i\ra
Y_i$ are fixed universal covering maps. Let $\tilde
S_i=p_i^{-1}(S_i)$. Then $\tilde S_i$ is the (topological) center
surface of $X_i$. Since $\G$ acts transitively on the set
$p^{-1}(t_j)$ for each fixed $j$, there is an element $\g_{i,j}$
of $\G$ such that $\g_{i,j}(\tilde S_i)$ contains the point $b_j$.
Let $X_{i,j}=\g_{i,j}(X_i)$. Then $\g_{i,j}(\tilde S_i)$ is the
center surface of $X_{i,j}=\g_{i,j}(X_i)$, and $X_{i,j}$ is
invariant under the action of the subgroup
$\g_{i,j}\G_i\g_{i,j}^{-1}$.
 Let $W_j= X_{1,j}\cap X_{2,j}$. Then $W_j$ is a
  metrically complete and  strictly convex (thus simply
  connected)
$3$-dimensional submanifold of $\mathbb H^3$ which is invariant
under the action of the subgroup $(\g_{1,j}\G_1\g_{1,j}^{-1})\cap
(\g_{2,j}\G_2\g_{2,j}^{-1})$. Let $Z_{i,j}=\g_{i,j}^{-1}(W_j)$.
Then $Z_{1,j}= X_1 \cap \gamma_{1,j}^{-1}\gamma_{2,j} (X_2)$ is
contained in $X_1$ and is invariant under the action of the
subgroup $\G_{i,j} = \G_1\cap
(\g_{1,j}^{-1}\g_{2,j}\G_2\g_{2,j}^{-1}\g_{1,j})$, and similarly
 $Z_{2,j}=X_2 \cap \gamma_{2,j}^{-1}\gamma_{1,j} (X_1)$ is contained in $X_2$
 and is invariant under the
action of the subgroup $\G_{2,j} = \G_2\cap
(\g_{2,j}^{-1}\g_{1,j}\G_1\g_{1,j}^{-1}\g_{2,j})$.

\begin{lem}\label{gij}
The subgroup $\G_{i,j}$ contains
 no parabolic elements, for any $i = 1,2, j = 1, ..., d$.
\end{lem}

\pf Recall that  $i$ and $i_*$ denote the number $1$ or $2$ such that
$\{i,i_*\}=\{1,2\}$, and that
 $\G_{i,j}=\G_i\cap
(\g_{i,j}^{-1}\g_{i_*,j}\G_{i_*}\g_{i_*,j}^{-1}\g_{i,j})$.
Also recall  $\G_i$, $i=1,2$, are the fundamental groups of
 two embedded, cusped, quasi-Fuchsian surfaces, with
different  boundary slopes. Thus no parabolic element in $\G_i$
is conjugate in $\G$ to any element  in $\G_{i_*}$ (cf. the
proof of Lemma 2.1 in \cite{CL2}). Hence the conclusion of the
lemma follows. \qed

 Recall that ${\cal B}_i$ is the horoball region of $X_i$,
which is the intersection of $X_i$ with the collection of horoballs
in ${\cal B}$ based at parabolic fixed points of $\G_i$. Note that
$\L(X_i)=\L(\G_i)$. We claim that the limit set $\L(Z_{i,j})$ of
$Z_{i,j}$ is equal to the intersection
$\L(\G_i)\cap\L(\g_{i,j}^{-1}\g_{i_*,j}\G_{i_*}\g_{i_*,j}^{-1}\g_{i,j})$.
Indeed, the containment $\L(Z_{i,j})\subset
\L(\G_i)\cap\L(\g_{i,j}^{-1}\g_{i_*,j}\G_{i_*}\g_{i_*,j}^{-1}\g_{i,j})$
is obvious.  For the other containment,
 suppose that $x$ is in $\L(\G_i)\cap\L(\g_{i,j}^{-1}\g_{i_*,j}\G_{i_*}\g_{i_*,j}^{-1}\g_{i,j})$.
   Then there are geodesic rays
 $\a$ and $\a^{\prime}$,
 contained in $H_i$ and $\gamma_{i,j}^{-1}\gamma_{i_*,j} H_{i_*}$ respectively,
 with $x$ as a limit endpoint.  Then far enough along these geodesics, each point
 in $\a$ is within an epsilon-neighborhood
 of $\a^{\prime}$, and vice versa.  Therefore,
 far enough along these geodesics, $\a$ and $\a^{\prime}$ are both contained
 in $X_i \cap  \gamma_{i,j}^{-1}\gamma_{i_*,j} X_{i_*}$,
 and therefore $x$ is a limit point of $Z_{i,j}$.

Since quasi-Fuchsian groups are geometrically finite,
 we may apply  Theorem 3.14 of \cite{MT} (which is originally due to
 Susskind \cite{Su})  to conclude that
$\L(\G_i)\cap\L(\g_{i,j}^{-1}\g_{i_*,j}\G_{i_*}\g_{i_*,j}^{-1}\g_{i,j})
=
\Lambda(\G_{i,j}) \cup P_{i,j}$, where $P_{i,j}$ is the set of
points $\zeta \in \Omega(\Gamma_{i,j})=S^2_\infty -\L(\G_{i,j})$
such that\\ (1) $Stab_{\G_i}(\zeta)$ and
$Stab_{\g_{i,j}^{-1}\g_{i_*,j}\G_{i_*}\g_{i_*,j}^{-1}\g_{i,j}}(\zeta)$
generate a rank two Abelian group, and\\ (2) $Stab_{\G_i}(\zeta)
\cap
Stab_{\g_{i,j}^{-1}\g_{i_*,j}\G_{i_*}\g_{i_*,j}^{-1}\g_{i,j}}(\zeta)
= \{id\}$.\\ Also, $\Lambda(\G_{i,j}) = \Lambda_p(\G_{i,j}) \cup
 \Lambda_c(\G_{i,j})$, where $\L_c$ denotes the set of
 conical limit points and $\L_p$ the set of parabolic limit points
 (see \cite{MT} or \cite{Ra} for their definitions).
By Lemma \ref{gij}, $\G_{i,j}$ contains no parabolic elements, and
thus
 $\Lambda_p(\G_{i,j}) = \emptyset$.
 Thus $\L(\G_i)\cap\L(\g_{i,j}^{-1}\g_{i_*,j}\G_{i_*}\g_{i_*,j}^{-1}\g_{i,j})
 = \Lambda_c(\G_{i,j}) \cup P_{i,j}$.

 Let ${\cal B}_{i,j}$ be the intersection of
$Z_{i,j}$ with the collection of horoballs in ${\cal B}$ based at
points of $P_{i,j}$. We call ${\cal B}_{i,j}$ the horoball region
of $Z_{i,j}$.  Let $Z_{i,j}^-=Z_{i,j}\setminus {\cal B}_{i,j}$,
which is the truncated version of $Z_{i,j}$. We call
$Z_{i,j}^-\cap \p {\cal B}_{i,j}$ the parabolic boundary of
$Z_{i,j}^-$ and denote it by $\p_p Z_{i,j}^-$. Note that
$Z_{i,j}^-$ is locally convex everywhere except on its parabolic
boundary. Each of $Z_{i,j}$,
 ${\cal B}_{i,j}$, $Z_{i,j}^-$ and $\p_p Z_{i,j}^-$ is invariant under the
action of $\G_{i,j}$.

Some members of $\{Z_{i,1},...,Z_{i,d}\}$ maybe the same
submanifold of $X_i$ modulo the action of $\G_i$ on $X_i$, i.e.
some one maybe a translation of another by an element of $\G_i$.

\begin{lem}
 The equality $Z_{1,j} = \g_1 (Z_{1,k})$ holds
 for some $\g_1 \in \G_1$
 if and only if  $Z_{2,j} = \g_2 (Z_{2,k})$ for some $\g_2 \in \G_2$.
\end{lem}

\pf Suppose that $Z_{1,j}=\g_1(Z_{1,k})$ for some element $\g_1\in \G_1$.
 Let
$\g_2=(\g_{2,j}^{-1}\g_{1,j})\g_1(\g_{1,k}^{-1}\g_{2,k})$. Then by our
construction, $\g_2$ maps $Z_{2,k}$ to $Z_{2,j}$.
 Also, $Z_{2,k}$ contains a point in $p^{-1} (t_k )\cap \tilde{S}_2$,
 and $\g_2$ maps this to another point in
 $p^{-1}( t_k) \cap \tilde{S}_2$.
Since the  $\G$-stabilizer of any point in
 $p^{-1}(t_k)$ is trivial and since
 $\G_2$ acts transitively on the
set $p^{-1}(t_{k})\cap \tilde S_2$, the element
$\g_2$ must belong to $\G_2$.
\qed

Let $j_1,...,j_q$ be such that $\{Z_{i,j_1},...,Z_{i,j_q}\}$ is a
maximal set of representatives of $\{Z_{i,1},...,Z_{i,d}\}$ which
are mutually inequivalent under the action of $\G_i$ on $X_i$ for
each $i=1,2$. Note that the set $\{Z_{i,j_1},...,Z_{i,j_q}\}$ is
well defined (independent of the choices for the points $b_j\in
p^{-1}(t_j) \cap \p B_\infty)$, up to translations by elements in
$\G_i$.

\begin{lem}\label{transitive}The subgroup $\G_{i,j_k}$ acts
transitively on $p^{-1}(t_j)\cap Z_{i,j_k}\cap \tilde S_i$, for
each fixed $j, k, i$.\end{lem}

\pf We know that $\Gamma_i$ acts transitively on $p^{-1}(t_j)\cap
\tilde S_i$
 and
$\g_{i,j_k}^{-1}\g_{i_*,j_k}\G_{i_*}\g_{i_*,j_k}^{-1}\g_{i,j_k}$
acts transitively on $p^{-1}(t_j)\cap
\g_{i,j_k}^{-1}\g_{i_*,j_k}(\tilde S_{i_*})$, so given two
distinct points $\tilde{t}$ and
        $\tilde{t}^{\prime}$ in $p^{-1}(t_j) \cap Z_{i,j_k}\cap \tilde S_i$,
        there exists $\gamma \in \Gamma_i$, and $\gamma' \in
    \g_{i,j_k}^{-1}\g_{i_*,j_k}\G_{i_*}\g_{i_*,j_k}^{-1}\g_{i,j_k}$
        such that each of them maps $\tilde{t}$ to $\tilde{t}^{\prime}$.
     But there is a unique element of $\G$ which maps $\tilde{t}$
     to $\tilde{t}^{\prime}$.
     Thus $\gamma = \gamma'$ and so
        $\gamma\in \G_i\cap
(\g_{i,j_k}^{-1}\g_{i_*,j_k}\G_{i_*}\g_{i_*,j_k}^{-1}\g_{i,j_k})=\G_{i,j_k}$.
\qed

 Each of the manifolds $Z_{i,j_k}$, $Z_{i,j_k}^-$,
$\p_pZ_{i,j_k}^-$ and ${\cal B}_{i,j_k}$ is invariant under the
action of the subgroup $$\G_{i,j_k}=\G_i\cap
(\g_{i,j_k}^{-1}\g_{i_*,j_k}\G_{i_*}\g_{i_*,j_k}^{-1}\g_{i,j_k}).$$
 Let $K_{i,k}=Z_{i,j_k}/\G_{i,j_k}$,
 $K_{i,k}^-=Z_{i,j_k}^-/\G_{i,j_k}$,  $\p_pK_{i,k}=\p_p
 Z_{i,j_k}^-/\G_{i,j_k}$, and  ${\cal C}_{i,k}={\cal B}_{i,j_k}/\G_{i,j_k}$.
For each $k=1,...,q$, $K_{1,k}$ and $K_{2,k}$ are isometric,
metrically complete, convex, hyperbolic manifolds.
 The isometry from $K_{1,k}$ to $K_{2,k}$ is the map $h_k$ which makes
the following diagram commute:
$$\displaystyle\begin{array}{ccc}
Z_{1,j_k}&\stackrel{\g_{2,j_k}^{-1}\circ\g_{1,j_k}}{\lra}
&Z_{2,j_k}\\\da&&\da
\\K_{1,k}&\stackrel{h_k}{\lra}&K_{2,k},\end{array}$$
where the vertical maps are the covering maps. Also for each $i$
and $k$, there is a local isometry  $g_{i,k}$ from $K_{i,k}$ into
$Y_i$ which is the restriction of the covering map $X_i/
\G_{i,j_k}\ra Y_i$. Let $K_i$ be the disjoint union of $\{K_{i,k},
k=1,...,q\}$. We have the isometry $h:K_1\ra K_2$ with
$h|_{K_{1,k}}=h_k$. We also have the local isometry $g_i:K_i
\rightarrow Y_i$ with $g_i|_{K_{i,k}}=g_{i,k}$.

\begin{lem}\label{ppz}
The restriction of the covering map $Z_{i,j_k}^-\ra K_{i,k}^-$ to
every component of  $\p_p Z_{i,j_k}^-$ is an isometric embedding,
for each of $i=1,2$ and each of $k=1,...,q$. In fact the
restriction of the covering map $Z_{i,j_k}\ra K_{i,k}$ to every
component of ${\cal B}_{i,j_k}$ is an isometric embedding, for
each of $i=1,2$ and each of $k=1,...,q$.
\end{lem}

\pf It follows from Corollary \ref{strip4} and the transitivity of
the action  of $\G$ on  the set $p^{-1}(t_j)$ (for any fixed $j$)
that every component of $\p_pZ_{i,j_k}^-$ is a Euclidean
parallelogram  in some horosphere. Now the first statement of
 the lemma follows from
the fact that $\G_{i,j_k}$ has no parabolic elements (Lemma
\ref{gij}). The second assertion can be proved similarly.
\qed

 We have just  shown that each component  $\tilde D$
 of $\p_pZ_{i,j_k}^-$ (for any $i, j_k$) is a Euclidean
parallelogram in a horosphere. We define the {\it (topological)
center point} of $\tilde D$ to be the point $\tilde D\cap\tilde
S_i\cap \g_{i,j_k}^{-1}\g_{i_*,j_k}(\tilde S_{i_*})$. The union of
all the center points in $\p_pZ_{i,j_k}^-$ is invariant under the
action of the subgroup $\G_{i,j_k}$.  By Lemma \ref{ppz},  each
component $D$ of $\p_pK_{i,k}$ is the isometric image of a component
$\tilde D$ of $\p_pZ_{i,j_k}^-$ under the covering map $Z_{i,j_k}\ra
K_{i,k}$. We define the (topological) center point of $D$  to be the
image of the center point of $\tilde D$.
 Thus by our construction, for each
$t_j \in f_1(\p S_1^-) \cap f_2(\p S_2^-)$,
 there is a component $D$ of $\p_p K_{i,k}^-$
 (for some $k$) whose center point is mapped to the point $t_j$
 under the map $K_{i,k}\ra Y_i\ra M$. In fact there is a geodesic ray,
  based at the center point, in the cusp part of $K_{i,k}$
 which maps isometrically to the ray $R_j\subset C$, under the map
$K_{i,k}\ra Y_i\ra M$.
 This component of $\p_p K_{i,k}^-$
 is said to be \textit{associated to the ray $R_j$} (thus to the point $t_j$),
 and so is the component
 $K_{i,k}$ of $K_i$.

\begin{lem}\label{finite}
 For each of $i=1,2$, the parabolic boundary of $K_i^-$
 has exactly $d$ components
 (each being a Euclidean parallelogram),
associated to the points $t_j$, $j=1,...,d$, respectively.
\end{lem}

\pf We prove this for $i=1$; the case for $i=2$ can be proved
similarly. By the construction, we see that the parabolic boundary
of $K_1^-$ has at least $d$ components, associated  to the points
$t_1,...,t_d$ respectively. Suppose that there are distinct
components $P_1$ and $P_2$ of the parabolic boundary of $K_1^-$
associated to the same point, say $t_1$.
 Also we may assume that  $K_{1,1}^-$ and $K_{1,k}^-$
are the components of $K_1^-$ containing $P_1$ and $P_2$
respectively. We first show that $k=1$ is impossible. So suppose
that both $P_1$ and $P_2$ are components of $\p_p K_{1,1}^-$. Recall
that $K_{1,1}^-=Z_{1,j_1}^-/\G_{1,j_1}$ and $\G_{1,j_1}=\G_1\cap
\g_{1,j_1}^{-1}\g_{2,j_1}\G_2\g_{2,j_1}^{-1}\g_{1,j_1}$. So the
parabolic boundary of $Z_{1,j_1}^-$ contains two components $\tilde
P_1$ and $\tilde P_2$ which are mapped to $P_1$ and $P_2$,
respectively, under the covering map $Z_{1,j_1}^-\ra K_{1,1}^-$.
 Because  the center points of $\tilde
P_1$ and $\tilde P_2$ are contained in $p^{-1}(t_1)$ and because
$\G_{1,j_1}$ acts on $p^{-1}(t_1)\cap Z_{1,j_1}^-\cap \tilde S_1$
transitively (Lemma \ref{transitive}), there is an element $\g\in
\G_{1,j_1}$ such that $\g(\tilde P_1) =\tilde P_2$. Hence both
$\tilde P_1$ and $\tilde P_2$ are mapped to $P_1$ under the covering
map $Z_{1,j_1}^-\ra K_{1,1}^-$, which gives a contradiction. Now
suppose that $k\ne 1$. Then $Z_{1,j_1}^-=X_1^-\cap
\g_{1,j_1}^{-1}\g_{2,j_1}( X_2^-)$ and $Z_{1,j_k}^-=X_1^-\cap
\g_{1,j_k}^{-1}\g_{2,j_k}(X_2^-)$ are two different submanifolds of
$X_1^-$, and there are two components  $\tilde P_1$ and $\tilde
P_2$, belonging to $\p_p Z_{1,j_1}^-$ and $\p_p Z_{1,j_k}^-$
respectively, which are mapped to $P_1$ and $P_2$ under the covering
maps $Z_{1,j_1}^-\ra K_{1,1}^-$ and $Z_{1,j_k}^-\ra K_{1,k}^-$
respectively.  Since  the center points of  $\tilde P_1$ and $\tilde
P_2$ are contained in $p^{-1}(t_1)$,  and $\G_1$ acts on
$p^{-1}(t_1)\cap \tilde S_1$ transitively, there is an element
$\g\in \G_1$ which maps $\tilde P_1$ to $\tilde P_2$.
 So $\g\g_{1,j_1}^{-1}\g_{2,j_1}( X_2)$
 intersects $X_1$ at
$\tilde P_2$. So $(\g\g_{1,j_1}^{-1}\g_{2,j_1})^{-1}(\tilde P_2)$
and $(\g_{1,j_k}^{-1}\g_{2,j_k})^{-1}(\tilde P_2)$ are both
contained in $X_2$. It follows that
$(\g\g_{1,j_1}^{-1}\g_{2,j_1})(\g_{1,j_k}^{-1}\g_{2,j_K})^{-1}$ is
contained in $\G_2$. Therefore
$(\g\g_{1,j_1}^{-1}\g_{2,j_1})(\g_{1,j_k}^{-1}\g_{2,j_k})^{-1}(X_2)
=X_2$, i.e.
$\g\g_{1,j_1}^{-1}\g_{2,j_1}(X_2)=\g_{1,j_k}^{-1}\g_{2,j_k}(X_2)$.
Hence $\g(Z_{1,j_1}) =X_1\cap
\g\g_{1,j_1}^{-1}\g_{2,j_1}(X_2)=X_1\cap
\g_{1,j_k}^{-1}\g_{2,j_k}(X_2)=Z_{1,j_k}$. Hence $Z_{1,j_1}$ and
$Z_{1,j_k}$ are equivalent under the translations of $\G_1$. This
gives a contradiction to our assumption that these $Z_{1,j_k}$,
$k=1,...,q$, are mutually inequivalent under translations of
elements of $\G_1$.
 \qed

\begin{lem}\label{cpt}  $K_{i,k}^-$ is compact for each $i=1,2$,
 $k=1,...,q$.
\end{lem}

\pf   Recall
 that   $K_{i,k} = Z_{i,j_k}/(\G_{i,j_k})$
 and $\G_{i,j_k}=\G_i\cap
(\g_{i,j_k}^{-1}\g_{i_*,j_k}\G_{i_*}\g_{i_*,j_k}^{-1}\g_{i,j_k})$.
The limit set of $Z_{i,j_k}$ is equal to $ \Lambda_c(\G_{i,j_k})
\cup P_{i,j_k}$ (see the two paragraphs following the proof of Lemma
\ref{gij}).

Since $K_{i,k}$ is convex and metrically complete, between any two
points in $K_{i,k}$ there is a distance minimizing geodesic
connecting them. Fix a point $k_0$ in $K_{i,k}^-$. Consider
$N_{(r, K_{i,k})}(k_0)$, the closed $r$-neighborhood of the point $k_0$
in $K_{i,j}$. Then by the Hopf-Rinow Theorem
(Theorem 1.3.5 of \cite{CEG}) $N_{(r, K_{i,k})}(k_0)$ is a compact
subset of $K_{i,k}$ for any $r>0$.

As a subset of $K_{i,k}$,  $K_{i,k}^-$ is closed. Hence
$K_{i,k}^-\cap N_{(r, K_{i,k})}(k_0)$ is a closed subset of the
compact set $N_{(r, K_{i,k})}(k_0)$ and thus is compact. Therefore
if $K_{i,k}^-$ is not compact, then it is not contained in $N_{(r,
K_{i,k})}(k_0)$ for any fixed $r>0$.
 So we can find a point $k_r$ in $K_{i,k}^-$ with $d(k_0,k_r)>r$
 for any $r>0$. By
Lemma \ref{finite}, the parabolic boundary of  $K_{i,k}^-$ has
finitely many components. Also each component of the parabolic
boundary of  $K_{i,k}^-$ is a compact Euclidean parallelogram.
Hence for all sufficiently large $r>0$, the parabolic boundary of
$K_{i,k}^-$ is contained in  $N_{(r, K_{i,k})}(k_0)$. Hence the
points $k_r$ are not in $\p_p(K_{i,k}^-)$ for all sufficiently
large $r>0$. Let $\a_r$ be the distance minimizing geodesic
segment in $K_{i,k}$ with endpoints $k_0$ and $k_r$.

Let $p_{i,k}:Z_{i,j_k}\ra K_{i,k}$ be the covering map. Pick a
point $z_0$ in $Z_{i,j_k}$ such that $p_{i,k}(z_0)=k_0$. Let
$\tilde \a_r\subset Z_{i,j_k}$ be the lift of $\a_r$ starting  at
$z_0$ (note that the lift is unique). Let $z_r$ be the other
endpoint of $\tilde \a_r$. Note that $\tilde \a_r$ is a distance
minimizing geodesic segment in $Z_{i,j_k}$ with $d(z_0,z)
=d(p_{i,k}(z_0),p_{i,k}(z))$ for any $z\in \tilde \a_r$. Now
consider the sequence of the geodesic segments $\{\tilde
\a_n\}_{n=1}^{\infty}$. As $d(z_0, z_n)=d(k_0,k_n)\ra +\infty$ as
$n\ra +\infty$, there is a subsequence of $\{z_n\}$ which
converges to a point $a$ in $S^2_\infty$. We may assume, for
simplicity in notation, that $\{z_n\}$ itself converges to $a$.
Now $a$ is in the limit set of $Z_{i,j_k}$ and thus $a\in
\Lambda_c(\G_{i,j_k}) \cup P_{i,j_k}$.

Let $\tilde \a$ be the geodesic ray in $\mathbb H^3$ starting at
$z_0$ and approaching $a$ (such ray exists and is unique). Since
$Z_{i,j_k}$ is metrically complete and convex, $\tilde \a$ is
contained in $Z_{i,j_k}$. In fact the sequence $\{\tilde \a_n\}$
is approaching $\tilde \a$ in the sense that every point $x$ in
$\tilde\a$ is the limit of a sequence of points $\{x_n\}$ with
$x_n\in \tilde \a_n$. It follows that   each finite sub-segment of
the projection $p_{i,k}(\tilde{\alpha})$ is a distance-minimizing
segment in $K_{i,k}$.

We first show that the point $a$ is not in $P_{i,j_k}$. For
otherwise if $B_a$ is the  horoball component in $\cal B$ based at
$a$, $\tilde \a$ intersects perpendicularly every horosphere
inside $B_a$ based at $a$. Hence  $\tilde \a\cap B_a$ is a
geodesic ray contained in $B_a\cap Z_{i,j_k}$. But $\tilde\a_n$ is
approaching $\tilde\a$, so for sufficiently large $n$, $z_n$ will
enter into $B_a$. So $z_n$ is not in $Z_{i,j_k}^-$ and thus
$p_{i,k}(z_n)=k_n$ is not contained in $K_{i,k}^-$, which
contradicts our construction of $k_n$.

\begin{figure}[!ht]
{\epsfxsize=3in \centerline{\epsfbox{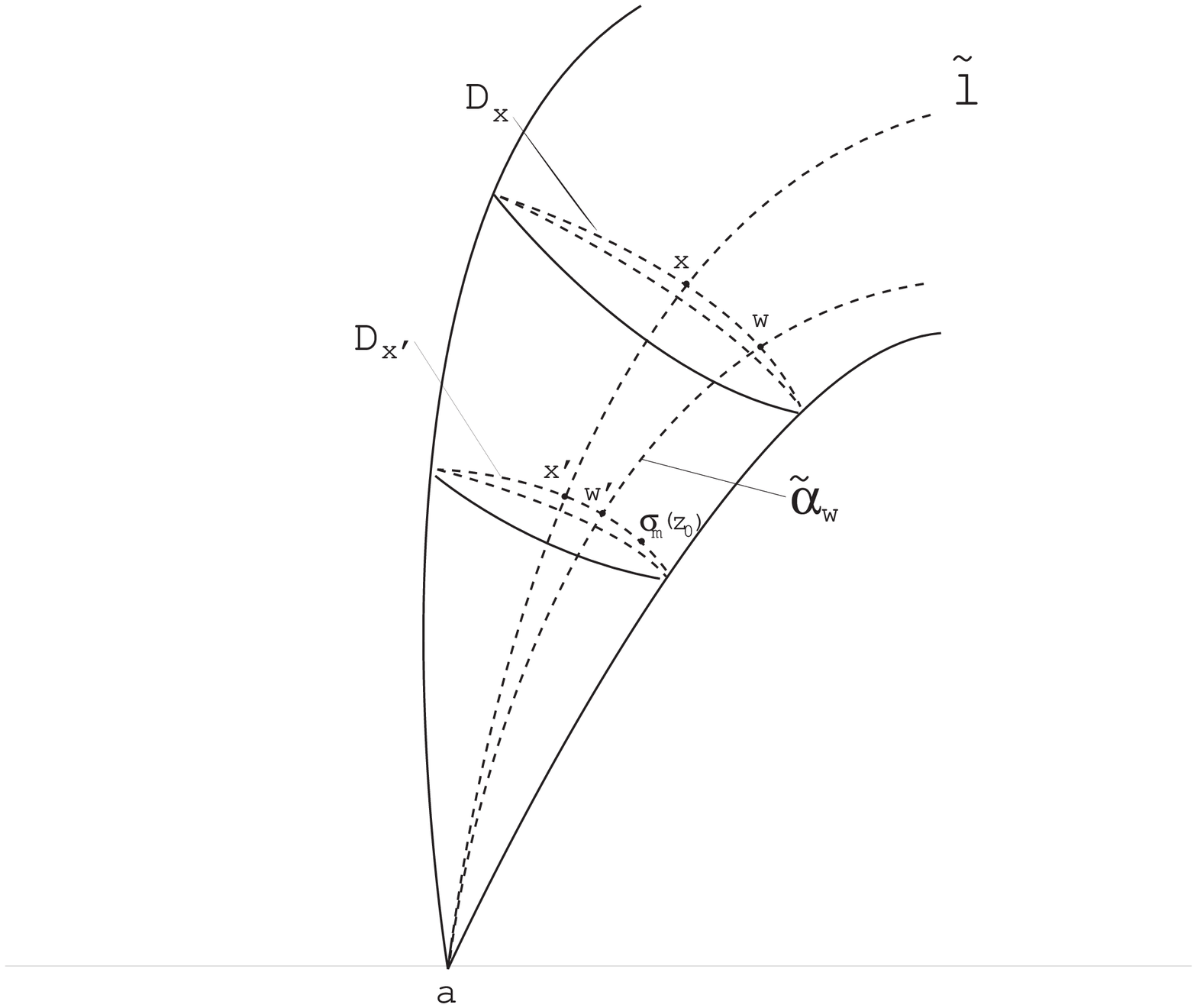}}\hspace{10mm}}
\caption{}\label{conical}
\end{figure}

So $a$ is a conical limit point of $\G_{i,j_k}$. By definition
there is a geodesic ray $\tilde l$ ending at $a$
  and there is a sequence of elements $\s_m$ in $\G_{i,j_k}$ such that
  $\s_m(z_0)$ is contained in  $N_{\e}(\tilde l)$ (which is a fixed $\e$-collared
  neighborhood of
  $\tilde l$ in $\mathbb H^3$) for all sufficiently large $m$,
  and converges  to the point $a$ as $m\ra \infty$.
 Now the ray $\tilde \a$ is contained in $N_\e(\tilde l)$
 except possibly for a finite initial segment.
 Let $\tilde \a_{w}$ be a sub-ray of $\tilde \a$ starting at
 a point $w\in \tilde\a$ such that $\tilde\a_{w}$ is entirely contained
 in $N_\e(\tilde l)$ and $d(z_0,w)>2\e$.
 For any point $x$ on $\tilde l$ let $P_x$ be the hyperbolic plane
 intersecting $\tilde l$ perpendicularly at the point $x$, and let
 $D_x=P_x\cap  N_\e(\tilde l)$.  Then $D_x$ is topologically a disk
separating $N_\e(\tilde l)$ into  two pieces one of which contains
the sub-ray of $\tilde l$ starting at $x$.  Now every point in the
ray $\tilde \a_{w}$ is contained in  $D_x$ for some $x\in \tilde
l$. In particular the endpoint $w$ of $\tilde\a_{w}$ is contained
in some $D_x$. Let $V$ be the component of $N_\e(\tilde l)-D_x$
which contains a sub-ray of $\tilde l$. Since the sequence $(\s_m(z_0))$
approaches $a$, there is some $\s_m(z_0)$ which is contained in $V$.
Let $D_{x'}$ be the disk defined above containing $\s_m(z_0)$.
Then $D_{x'}$ is contained in $V$ and it also intersects a point
$w'$ in $\tilde \a_w$ (cf. Figure \ref{conical}). So $d(\s_m(z_0),
w') \leq d(\s_m(z_0),x')+d(x',w')\leq \e+\e$.  Hence
 $2\e<d(z_0,w)\leq d(z_0,w')=d(p_{i,k}(z_0), p_{i,k}(w'))=
 d(p_{i,k}(\s_m(z_0)),p_{i,k}(w'))\leq d(\s_m(z_0),w')\leq 2\e$
 (here  the first equality follows from the property that $\tilde \a$
 is a distance minimizing curve), giving a contradiction.
 \qed

Let $R$ be a fixed number bigger than  the number $R(\e)$ provided
in Proposition \ref{convh} and also bigger than the upper bound
 provided by Corollary \ref{length2} for
 the lengths of $I$-fibers of $X_i$ (for each of $i=1,2$).
 Consider $N_{(R, X_i)}(Z_{i,j_k})$, the
$R$-collared neighborhood of $Z_{i,j_k}$ in $X_i$. It is a convex
3-submanifold in $X_i\subset \mathbb H^3$ (thus is simply
connected) and is invariant under the action of $\G_{i,j_k}$. We let
$AN_{(R, X_i)}(K_{i,k})$ denote the quotient space
$N_{(R,X_i)}(Z_{i,j_k})/\G_{i,j_k}$, and call it the {\it
abstract $R$-collared neighborhood of $K_{i,k}$ with respect to
$X_i$}. Similarly we can define the truncated version of
$N_{(R,X_i)}(Z_{i,j_k})$ and the truncated version of
$AN_{(R,X_i)}(K_{i,k})$, denoted by $(N_{(R,X_i)}(Z_{i,j_k}))^-$
and $(AN_{(R,X_i)}(K_{i,k}))^-$ respectively. It follows from
Lemma \ref{cpt}  that $(AN_{R,X_i}(K_{i,k}))^-$ is compact.
 We can extend $g_{i,k}:K_{i,k} \ra Y_i$ to a map,
 which we still denote $g_{i,k}$, from $AN_{(R,X_i)}(K_{i,k})$
 to $Y_i$.

 By
construction, $N_{(R,X_i)}(Z_{i,j_k})$ contains all the $I$-fibers
of $X_i$ which meet $Z_{i,j_k}$. Let $Z_{i,j_k}'$ be the
sub-I-bundle of $X_i$ consisting of all the $I$-fibers of $X_i$
which meet $Z_{i,j_k}$. It is easy to see that $Z_{i,j_k}'$ is a
manifold. The manifold  $Z_{i,j_k}'$ is also invariant under the
action of $\G_{i,j_k}$ since $Z_{i,j_k}$ is invariant under the
action of $\G_{i,j_k}$ and since the action of $\G_{i,j_k}\subset
\G_i$ on $X_i$ sends fibers to fibers. Hence
$Z_{i,j_k}'/\G_{i,j_k}=F_{i,k}\times I$ for some surface $F_{i,k}$
(which is non-compact), with the induced $I$-fiber structure. From
the inclusions $K_{i,k}\subset F_{i,k}\times I\subset
AN_{(R,X_i)}(K_{i,k})$ and from the fact that the inclusion map
$K_{i,k}\subset AN_{(R,X_i)}(K_{i,k})$ induces an  isomorphism on
the fundamental groups, we see that the inclusion map $F_{i,k}\times
I\subset AN_{(R,X_i)}(K_{i,k})$ induces a surjective homomorphism on
the fundamental groups.

Note that $\p F_{i,k}\times I$ is precisely the frontier of
$F_{i,k}\times I$ in $AN_{(R,X_i)}(K_{i,k})$.
 Each component of $\p F_{i,k}\times I$ is either an annulus or a strip,
where a strip means $\mathbb R\times I$.

\begin{lem}\label{subbdle1}
Let $A$ be an annulus component of $\p F_{i,k}\times I$. Then
\\(1) $A$ divides $AN_{(R,X_i)}(K_{i,k})$ into two components
$B_1$ and $B_2$.
\\
(2) Suppose $B_1$ is the component whose interior is disjoint from
$F_{i,k}\times I$. Then either $B_1=D\times I$, where $D$ is a
disk, such that $A=\p D\times I$; or $B_1= S^1\times D$, where $D$
is a disk, such that $A=I\times S^1$, where $I$ is an interval
contained in $\p D$.
 \end{lem}

\pf  Since $F_{i,k}\times I$ is a submanifold of
$AN_{(R,X_i)}(K_{i,k})$ and carries the fundamental group of
$AN_{(R,X_i)}(K_{i,k})$, it follows  that $A$ is separating in
$AN_{(R,X_i)}(K_{i,k})$, i.e. we have (1). Part (2) also follows
easily. \qed

Similarly we have

\begin{lem}\label{subbdle2}
Let $E$ be a strip component of $\p F_{i,k}\times I$. Then \\(1)
$E$ divides $AN_{(R,X_i)}(K_{i,k})$ into two components $B_1$ and
$B_2$.
\\
(2) Suppose $B_1$ is the component whose interior is disjoint from
$F_{i,k}\times I$. Then $B_1=\mathbb R\times D$, where $D$ is a
disk,  such that $E=\mathbb R\times I$, where $I$ is an interval
contained in $\p D$. \qed \end{lem}

It follows from Lemmas \ref{subbdle1} and \ref{subbdle2} that the
$I$ bundle structure of $F_{i,k}\times I$ can be extended to one
on $AN_{(R,X_i)}(K_{i,k})$ in an obvious way.

Similarly one can obtain corresponding results in the truncated
setting. Namely $(AN_{(R,X_i)}(K_{i,k}))^-$ has a sub-manifold of
the form $F_{i,k}^-\times I$ which carries the fundamental group
of $(AN_{(R,X_i)}(K_{i,k}))^-$, and the  $I$-bundle structure of
$F_{i,k}^-\times I$ can be extended to one on
$(AN_{(R,X_i)}(K_{i,k}))^-$ in an obvious way, such that the
parabolic boundary of $(AN_{(R,X_i)}(K_{i,k}))^-$ consists of
$I$-fibers. Note that the $I$-fiber structure may not agree
 with the original $I$-fiber structure on $X_i/\G_{i,j}$.

Since $(AN_{(R,X_i)}(K_{i,k}))^-$ is compact for each $i=1,2$,
$k=1,...,q$, and since the horoball region of
$AN_{(R,X_i)}(K_{i,k})$ has a standard shape, we may assume, up to
replacing the cusp $C$ of $M$ by a smaller one, that $g_{i,k}^{-1}
({\cal C}_i)\cap (AN_{(R,X_i)}(K_{i,k}))^- =
\p_p(AN_{(R,X_i)}(K_{i,k}))^-$
 (where ${\cal C}_i$ is the cusp part of $Y_i$) for each of $i=1,2$,
$k=1,....,q$.

 We let
$AN_{(R,X_i)}(K_i)$ denote the disjoint union of
$\{AN_{(R,X_i)}(K_{i,k}); k=1,...,q\}$,
 and let $g_i: AN_{(R,X_i)}(K_i) \ra Y_i$, extending
 the local isometries $g_{i,k}$.
For later use, we record the  following corollary.

\begin{cor}\label{subbdle3}
Suppose that the local isometry $g_i: (AN_{(R,X_i)}(K_i))^-=\amalg_k
(AN_{(R,X_i)}(K_{i,k}))^-\ra Y_i^-$ lifts to an embedding in a
finite cover $\breve Y_i^-$ of $Y_i^-$. Then the $I$-bundle
structure on $\breve Y_i^-$ can be adjusted to one so that the image
of $(AN_{(R,X_i)}(K_i))^-$ is a sub-$I$-bundle in $\breve
Y_i^-$.\qed
\end{cor}

Let  $Fr_{X_i}(N_{(R,X_i)}(Z_{i,j_k}))$ denote the frontier of $
N_{(R,X_i)}(Z_{i,j_k})$ in $X_i$.
 If we define the frontier boundary
$\p_{f}(AN_{(R,X_i)}(K_{i,k}))$ of $AN_{(R,X_i)}(K_{i,k})$ to be
$$Fr_{X_i}(N_{(R,X_i)}(Z_{i,j_k}))/\G_{i,j_k},$$ then
$\p_f(AN_{(R,X_i)}(K_{i,k}))$ is topologically parallel to $\p
F_{i,k}\times I$ by Lemmas \ref{subbdle1} and \ref{subbdle2}. Thus
each component of $\p_f(AN_{(R,X_i)}(K_{i,k}))$ is either an annulus
or a strip. A strip component must enter the cusp region of
$AN_{(R,X_i)}(K_{i,k})$. From the shape of the cusp region of
$AN_{(R,X_i)}(K_{i,k})$ and from Lemma \ref{finite}, we see that
 the
frontier boundary of $AN_{(R,X_i)}(K_i)$ has exactly $d$ strip
components and that the frontier boundary of each component $
AN_{(R,X_i)}(K_{i,k})$ of $ AN_{(R,X_i)}(K_{i})$ has at least two
strip components. We restate this fact in the following corollary
for later use.

\begin{cor}\label{d-strips}(1) $\p_f( AN_{(R,X_i)}(K_{i,k}))$ has at
least two strip components for each $i=1,2$ and $k=1,...,q$. (2)
$\p_f(AN_{(R,X_i)}(K_i))$ has exactly $d$ strip components for
each $i=1,2$.\qed\end{cor}

The following corollary   follows easily from the convexity of
$N_{(R,X_i)}(Z_{i,j_k})$ and from the shape of the parabolic
region of $N_{(R,X_i)}(Z_{i,j_k})$.

\begin{cor}\label{diffcusps}
Every component of $Fr_{X_i}(N_{(R,X_i)}(Z_{i,j_k}))$  has its two
ends contained in two different horoball components of $\cal B$
respectively.\qed\end{cor}

We conclude with some remarks.

\begin{rem} Results and notations in this section will
still be valid if we replace the cusp $C$ by a smaller one.
\end{rem}

\begin{rem}In the construction of $K_{i,k}$ and its local isometry
$g_{i,k}$ into $Y_i$ some choices were made (for instance the
universal cover $Z_{i,j_k}$ in $X_i$ which is defined up to
translation by elements of $\G_i$). But, up to isometry, the construction is
independent of all such choices; i.e. if $g_{i,k}':K_{i,k}'\ra Y_i$ is another
result of this construction, then there is  an isometry
$\phi_{i,k}: K_{i,k}\ra K_{i,k}'$ such that $g_{i,k}=g_{i,k}'\circ
\phi_{i,k}$.
\end{rem}

\section{Constructing  $J_i$}\label{ji}

In Section \ref{kiyi}, we constructed, for each $i=1,2$, the
manifold $AN_{(R,X_i)}(K_i)$, which
is the disjoint union of $\{AN_{(R,X_i)}(K_{i,k}); k=1,...,q\}$,
such that each component of
$AN_{(R,X_i)}(K_i)$ is a metrically complete, convex, hyperbolic
$3$-manifold, and we defined a local isometry
 $g_i:AN_{(R,X_i)}(K_i) \ra Y_i$.
 In this section we construct, for each $i=1,2$, a
{\bf connected},
 metrically complete,  convex, hyperbolic 3-manifold $J_i$
with a local isometry $g_i:J_i \ra Y_i$, such that $J_i$ contains
$AN_{(R,X_i)}(K_i)$ as a hyperbolic submanifold, and
 $J_i\setminus AN_{(R,X_i)}(K_i)$  is a compact $3$-manifold $W_i$
 (which may not be
connected). Obviously we may assume that $q>1$, since otherwise we
may simply take $J_i=AN_{(R,X_i)}(K_{i,1})$.

We continue to use the notations established in early sections. We
have showed (Corollary \ref{d-strips}) that $\p_f(
AN_{(R,X_i)}(K_{i,k}))$ has at least two strip components for each
 $k=1,...,q$ and that $\p_{f}(AN_{(R,X_i)}(K_i))$ has exactly
$d$ strip components. Let $E_{i,k}$ be a fixed strip component of
$\p_{f}AN_{(R,X_i)}(K_{i,k})$ for each fixed $i$ and $k$. Recall
$p_{i,k}:N_{(R,X_i)}(Z_{i,j_k})\ra AN_{(R,X_i)}(K_{i,k})$ is the
universal covering map.  Then each component of
$p^{-1}_{i,k}(E_{i,k})\subset Fr_{X_i}(N_{(R,X_i)}(Z_{i,j_k}))$ is
a strip isometric to $E_{i,k}$ under the map $p_{i,k}$. Let
$\tilde E_{i,k}$ be a fixed component of $p^{-1}_{i,k}(E_{i,k})$.

The required $J_i$ will be constructed by gluing  components of
$AN_{(R,X_i)}(K_i)$ with
 a compact $3$-manifold $W_i$, along an  attaching
region in $\amalg E_{i,k}$. We shall construct the connecting
manifold $W_i$ in $X_i$. The procedure is as follows: find a
suitable translation of $N_{(R,X_i)}(Z_{i,j_k})$ by an element
$\t_{i,k}$ of $\G_i$, then take the convex hull of
$\{\t_{i,k}(N_{(R,X_i)}(Z_{i,j_k})); k=1,...,q\}$ in $X_i$. The
added part in forming the convex hull is the manifold $W_i$, which
will be shown to be compact, and the attaching region of $W_i$
with $\t_{i,k}(N_{(R,X_i)}(Z_{i,j_k}))$ is contained in
$\t_{i,k}(\tilde E_{i,k})$.

If such $W_i$ can be found, we can glue it with each
$AN_{(R,X_i)}(K_{i,k})$ along $E_{i,k}$ using the isometry
$\t_{i,k}(\tilde E_{i,k}) \stackrel {\t_{i,k}^{-1}}{\ra}\tilde
E_{i,k}\stackrel{p_{i,k}} {\ra}E_{i,k}$.  The resulting manifold
$J_i$ is a convex hyperbolic $3$-manifold,  with a local isometry
into $Y_i$, extending the map $g_i: AN_{(R,X_i)}(K_{i}) \ra Y_i$.
 It is easy to
see that the hyperbolic structure in $AN_{(R, X_i)}(K_i)$ and the
hyperbolic structure on $W_i$ match up along their gluing
surfaces, forming a global hyperbolic structure for $J_i$.

We now give the construction of $W_i$,
 beginning with some well known
facts. Let $\g$ be any hyperbolic element of
 $PSL_2(\c)$. The
axis of $\g$ is denoted $A_\g$. Let $a, a'$ be the
two limit points of $A_\g$, which are
the two fixed points of $\g$ in $S_\infty^2$ . Then for any point $x$ in
$\overline{\mathbb H}^3$, the sequence $\g^n(x)$ approaches one of the points
$a,a'$, say $a$, as $n\ra \infty$, and approaches $a'$, as $n\ra
-\infty$. Thus for any fixed closed subset $W$ of
$\overline{\mathbb H}^3$ which is disjoint from $a'$, and for any
fixed open neighborhood $U$ of $a$ in $\overline{\mathbb H}^3$
there is an integer $n$ such that $\g^{n}(W)\subset U$.

\begin{lem}\label{axis}
For any open arc $\a$ in $\L_i=\L(\G_i)$, there exists a
hyperbolic element $\g$ of $\G_i$ such that the two limit points
of $A_\g$ are contained in $\a$.
\end{lem}

\pf Since fixed points of hyperbolic elements of $\G_i$ are
dense in $\L_i$,  there is a hyperbolic element $\d$ in $\G_i$
with at least one of its two fixed points  contained in $\a$.  Now
take a hyperbolic element $\eta$ of $\G_i$ such that the limit
points of $A_\eta $ are disjoint from the limit points of $A_\d$.
By the notes given proceeding the lemma, there is an integer $n$
such that the two limit  points of $\d^n(A_\eta)$ are both
contained in $\a$. Let $\g=\d^n\eta\d^{-n}$, then
$A_\g=\d^n(A_\eta)$. \qed

Each strip $\tilde E_{i,k}$ (defined earlier in this section) has
exactly two limit points in $\L_i$ and each of them is a parabolic
fixed point of $\G_i$ (the two parabolic fixed points are distinct
because of the convexity of $N_{(R,X_i)}(Z_{i,j_k})$). Note that
$\tilde E_{i,k}$ separates $X_i$ into two parts, one of which
contains $N_{(R,X_i)}(Z_{i,j_k})$. Let $U_{i,k}$ be the part whose
interior is disjoint from $N_{(R,X_i)}(Z_{i,j_k})$. Let $\a_{i,k}$
be the  limit set of  $U_{i,k}$.

\begin{lem}\label{uij}
After translations by suitable
 elements of $\G_i$, we may assume
 that $N_{(R,X_i)}(Z_{i,j_2})$  is contained in the
interior of $U_{i,1}$ and that $N_{(R,X_i)}(Z_{i,j_1})$  is
contained in the interior of $U_{i,2}$
\end{lem}

\pf
Let $\overline{X}_i$ denote the closure of $X_i$ in
$\overline{\mathbb H}^3$. Note that $N_{(R,X_i)}(Z_{i,j_k})$ has
the same  limit set as $Z_{i,j_k}$. Let $\g$ be a hyperbolic
element of $\G_i$ whose axis has a limit point, $a$, disjoint from
the limit points of $N_{(R,X_i)}(Z_{i,j_2})$. Then by the notes
given in the paragraph proceeding Lemma \ref{axis}, we may move
$N_{(R,X_i)}(Z_{i,j_2})$  by a power of $\g$ into a small open
neighborhood of  $a'$ (which is the other limit point of $A_\g$)
in $\overline{X}_i$. In particular, we may assume
 that the limit set of this translate does
 not contain the limit set of $U_{i,1}$. Then,
 applying Lemma \ref{axis}, there
 is an element in $\G_i$ which translates
 $N_{(R,X_i)}(Z_{i,j_2})$ into $U_{i,1}$.

 Thus we may assume that
$N_{(R,X_i)}(Z_{i,j_2})$ is contained in the interior of
$U_{i,1}$. If $N_{(R,X_i)}(Z_{i,j_1})$ is contained in the
interior of $U_{i,2}$ already, then we are done. So suppose not.
Then $U_{i,2}$ is contained in the interior of $U_{i,1}$. Let $\g$
be a hyperbolic element of $\G_i$ such that the two limit points
of  $A_\g$ are contained in the interior of the arc $\a_{i,2}$;
such an element exists by Lemma \ref{axis}. Then, after replacing
$N_{(R,X_i)}(Z_{i,j_2})$ by its translate
 under a suitably high power of $\g$, one may check that the
 conclusion of the lemma is satisfied . \qed

By Lemma \ref{uij}, $\tilde E_{i,1}$ and $\tilde E_{i,2}$ co-bound
a connected submanifold of $X_i$, $V_1$, whose interior is
disjoint from both $N_{(R,X_i)}(Z_{i,j_1})$ and
$N_{(R,X_i)}(Z_{i,j_2})$. The limit set of $V_1$ consists of two
disjoint arcs in $\L_i$.  Now if $q>2$, then by a method similar to
 the proof of Lemma \ref{uij}, we may assume, up to translation by a hyperbolic
 element of $\G_i$, that $N_{(R,X_i)}(Z_{i,j_3})$ is in
 in the interior of $V_1$, and that both
$N_{(R,X_i)}(Z_{i,j_1})$ and $N_{(R,X_i)}(Z_{i,j_2})$ are
contained in the interior of $U_{i,3}$. In other
words, the three strips $\tilde E_{i,1}$, $\tilde E_{i,2}$ and
$\tilde E_{i,3}$ co-bound a connected submanifold $V_2$ in $X_i$
such that the interior of $V_2$ is disjoint from
$N_{(R,X_i)}(Z_{i,j_1})$, $N_{(R,X_i)}(Z_{i,j_2})$ and
$N_{(R,X_i)}(Z_{i,j_3})$.

By a simple induction, we may assume that  $N_{(R,X_i)}(Z_{i,j_k}),
k=1,...,q$, are located in $X_i$ in such way that the $q$ strips
$\tilde E_{i,k}$, $k=1,...,q$, co-bound a connected submanifold
$V$  of $X_i$ whose interior is disjoint from
$N_{(R,X_i)}(Z_{i,j_k}), k=1,...,q$. Now we take the convex hull
of the set $\{N_{(R,X_i)}(Z_{i,j_k}), k=1,...,q\}$ in $X_i$, and
let $Z_i$ be the resulting convex manifold. Let $W_i$ be the
complement of the interior of $N_{(R,X_i)}(Z_{i,j_k}), k=1,...,q$,
in $Z_i$. Then, by Lemma \ref{convh2},
 $W_i$ is a compact submanifold of $X_i$.
 This $W_i$ is the desired connecting manifold.
 The attaching region in $\p W$, to be glued to $E_{i,k}$,
 is $W_i\cap \tilde E_{i,k}$.

We still use $g_i$ to denote the  local isometry $J_i\ra Y_i$.
Since $W_i$ is compact, there exists a cusp $C'$ of $M$, smaller
than or equal to $C$, such that $g_i(W_i)$ is disjoint from the
corresponding cusp region ${\cal C'}_i$ of $Y_i$. We may assume
that the cusp $C$ itself already satisfies this condition. Under
this assumption, $W_i$ is disjoint from
$AN_{(R,X_i)}(K_i)\setminus (AN_{(R,X_i)}(K_i))^-$, and  the
components of $(AN_{(R,X_i)}(K_i))^-$ are
 connected together by $W_i$ along  the
frontier boundary of $AN_{(R,X_i)}(K_i)$, forming a connected
compact manifold which we denote by $J_i^-$. The parabolic
boundary $\p_p J_i^-$ is defined to be the parabolic boundary of
$(AN_{(R,X_i)}(K_i))^-$. Then $g_i|: (J_i^-, \p_p J_i^-)\ra
(Y_i^-, \p_pY_i^-)$ is a proper map of pairs.

Each component of $\p_p J_i^-$ can be isometrically  embedded in
$\p B_\infty$ as a Euclidean parallelogram. The convex hull of
such a parallelogram is a convex 3-ball in $B_\infty$ lying
vertically above the parallelogram.   We let $\hat J_i$ denote the
manifold obtained by capping off each of its parabolic
boundary component by a convex $3$-ball as just described. Then
$\hat J_i$ is a connected, compact, convex $3$-manifold with a local
isometry (which we still denote by $g_i$) into $Y_i$.

\section{Constructing  $C_n(J_i^-)$}\label{cni}

From Sections \ref{cqf}, \ref{kiyi} and \ref{ji}, we have the
following setting: for each $i=1,2$, $f_i:Y_i=S_i\times I\ra M$ is
a local isometry; $f_i|: (Y_i^-=S_i^-\times I, \p_p Y_i^-= \p
S_i^-\times I)\ra (M, \p M)$ is a proper map; $f_i|:\p S_i^-\ra \p
M$ is an embedding; $\p S_i^-$ has $n_i$ components $\{\b_{i,j},
j=1,...,n_i\}$ with induced orientation;  $\D$ is the geometric
intersection number between $f_1(\b_{1,1})$ and $f_2(\b_{2,1})$;
there are $d=\D n_1n_2$ intersection points $\{t_1,...,t_d\}$
between $f_1(\p S_1^-)$ and $f_2(\p S_2^-)$ in $\p M$;
$\{t_{i,j,k}, j=1,...,n_i, k = 1, ..., d_i\}$ are the
 points in $f_i^{-1} \{t_1, ..., t_d\}$, where $d_i=\D n_{i_*}$,
 indexed so that  $\{t_{i,j,k}, k = 1, ..., d_i\}$ are contained
successively in the component $\b_{i,j}$ (following the orientation
of $\b_{i,j}$) for each $j=1,...,n_i$; $K_i$ is the disjoint union
of the ``intersection manifolds'' $\{K_{i,j},
 j=1,...,q\}$; the manifold $AN_{(R,X_i)}(K_{i})=\Pi AN_{(R,X_i)}(K_{i,k})$
 is the abstract $R$-collared
 neighborhood of $K_i$ with respect to $X_i$;
 each component $AN_{(R,X_i)}(K_{i,k})$
 is a metrically complete, convex, hyperbolic
$3$-manifold with a local isometry
$g_{i,k}:AN_{(R,X_i)}(K_{i,k})\ra Y_i$; the restriction $g_{i,k}|:
((AN_{(R,X_i)}(K_{i,k}))^-, \p_p(AN_{(R,X_i)}(K_{i,k}))^-)\ra
(Y_i^-,\p_p Y_i^-)$ is a proper map; $g_i:AN_{(R,X_i)}(K_i)\ra
Y_i$ is the local isometry with
$g_i|_{AN_{(R,X_i)}(K_{i,k})}=g_{i,k}$; $J_i$ is a  metrically
complete convex (thus connected) hyperbolic $3$-manifold  with
local isometry $g_i: J_i\ra Y_i$; $g_i|: (J_i^-, \p_p J_i^-)\ra
(Y_i^-, \p_pY_i^-)$ is a proper map; $J_i$ contains
$AN_{(R,X_i)}(K_i)$ as a submanifold;
 and $\p_p J_i^-=\p_p
(AN_{(R,X_i)}(K_i))^-$.

Also recall that there are exactly $d$ components in $\p_p
J_i^-=\p_p (AN_{(R,X_i)}(K_{i}))^-$, one each associated to the
points $t_1,...,t_d$ respectively.  Let $D_{i,j,k}$, $j=1,...,n_i,
k=1,...,d_i$ denote the components of $\p_pJ_i^-$ and let
$b_{i,j,k}$ be the topological center point of $D_{i,j,k}$,
indexed so that $g_i(b_{i,j,k})=t_{i,j,k}$.

 The purpose of this section is to construct,
for each  sufficiently large integer $n$,  a connected, compact,
convex, hyperbolic $3$-manifold $C_n(J_i^-)$, with a local isometry
into $Y_i$, such that $C_n(J_i^-)$  contains $J_i^-$ as a
hyperbolic submanifold. $C_n(J_i^-)$ is obtained by gluing
together $J_i^-$ with $n_i$ ``multi-$1$-handles'' $H_{i,j}^n,
j=1,...,n_i$, along the attaching region $\p_p J_i^-$ (see Fig. 4
 for a preview).
A more precise description of $C_n(J_i^-)$ will be clear after its
construction. The needed properties of $C_n(J_i^-)$ will be
described  in later sections.

Now we proceed to construct the multi-$1$-handle $H_{i,j}^n$ for
each fixed $i\in\{1,2\}$ and each fixed $j\in\{1,...,n_i\}$.
 Let $c_{i,j}$ be a fixed
component of $p^{-1}(f_i(\b_{i,j}))$ in the horizontal
horosphere $\p B_\infty$. The transitivity of the action of  $\G$
implies that  there is element $\d_{i,j}\in \G$ such that
$\d_{i,j}(\tilde S_i)$ contains  $c_{i,j}$. By Corollary
\ref{strip4}, $\d_{i,j}(X_i)\cap \p B_\infty$ is a strip in $\p
B_\infty$ between two parallel Euclidean lines, which contains
$c_{i,j}$ as its (topological) center line. Let $E_{i,j}$ denote
this strip.

Along $c_{i,j}$,  we index the set of points
$p^{-1}(\{f_i(t_{i,j,k}), k=1,...,d_i\})$ as $a_{i,j,k,m}$,
$k=1,...,d_i$, $m\in \mathbb Z$, such that \newline (1) for each
fixed $m$, the points $\{a_{i,j,k,m}, k=1,..., d_i\}$ appear
consecutively along the line $c_{i,j}$ following the orientation
of $c_{i,j}$ (which is induced from that of $\b_{i,j}$);\newline
(2) the point $a_{i,j,d_i,m}$ is  followed immediately by the
point $ a_{i,j, 1, m+1}$, for every $m$;\newline (3)
$p(a_{i,j,k,m})=f_i(t_{i,j,k})$, for all $k=1,...,d_i$ and $m\in
\z$.

 For an arbitrary (fixed) sufficiently large integer $n>0$,
 consider the following $d_i$
points on $c_{i,j}$: $a_{i,j,1,0}$, $a_{i,j,2,n}$, $a_{i,j,3,2n}$,
..., $a_{i,j,d_i,(d_i-1)n}$. Again by transitivity of the action
of $\G$, there are elements $\g_{i,j,1}, \g_{i,j,2},
...,\g_{i,j,d_i}\in \G$ such that $\g_{i,j,1}(\tilde S_{i_*}),
\g_{i,j,2}(\tilde S_{i_*}), ..., \g_{i,j,d_i}(\tilde S_{i_*})$
contain the points $a_{i,j,1,0}$, $a_{i,j,2,n}$,  ...,
$a_{i,j,d_i,(d_i-1)n}$ respectively. Consider the corresponding
translations of $X_{i_*}$: $\g_{i,j,1}( X_{i_*}), \g_{i,j,2}(
X_{i_*}), ..., \g_{i,j,d_i}(X_{i_*})$. Each of
$\d_{i,j}(X_i)\cap\g_{i,j,1}( X_{i_*}),
\d_{i,j}(X_i)\cap\g_{i,j,2}( X_{i_*}), ...,  \d_{i,j}(X_i)\cap
\g_{i,j,d_i}( X_{i_*})$ is a translation of some component in
$\{Z_{i,j_1},..., Z_{i,j_q}\}$. Let $Z_1,...,Z_{d_i}$ denote
$\d_{i,j}(X_i)\cap\g_{i,j,1}( X_{i_*}),
 ...,  \d_{i,j}(X_i)\cap
\g_{i,j,d_i}( X_{i_*})$ respectively, and let $N_R(Z_i) =
N_{(R,\d_{i,j}(X_i))}(Z_i)$.
 Each of $N_R(Z_1),..., N_R(Z_{d_i})$ is a translation of some component in
$\{N_{(R,X_i)}(Z_{i,j_1}),..., N_{(R,X_i)}(Z_{i,j_q})\}$.

\begin{figure}[!ht]
{\epsfxsize=6in \centerline{\epsfbox{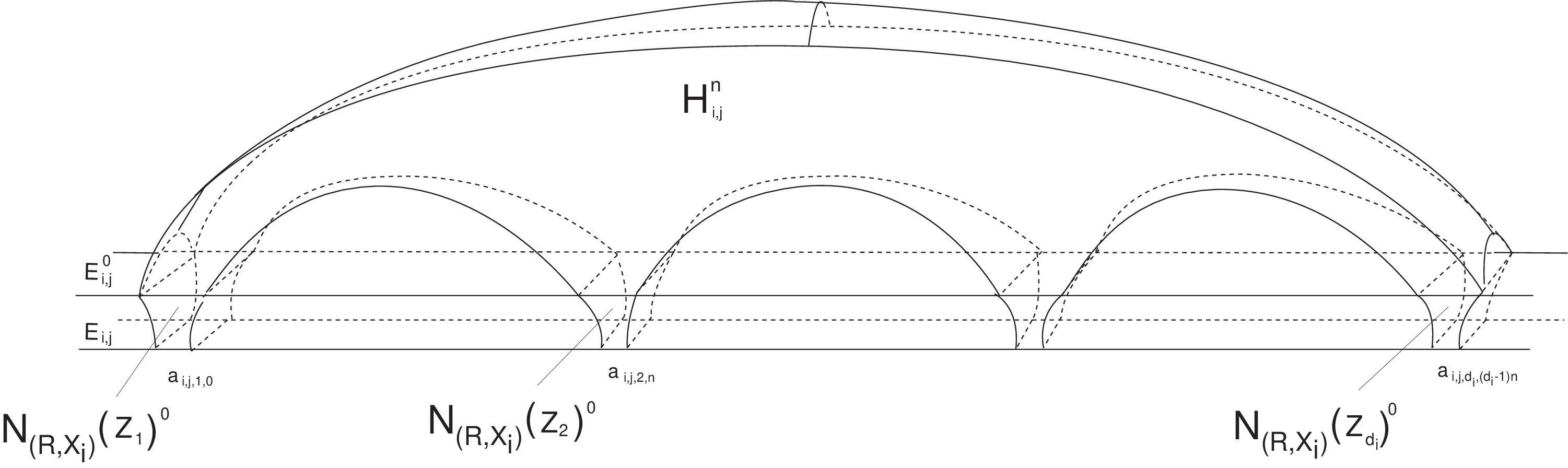}}\hspace{10mm}}
\caption{}\label{cin}
\end{figure}

Let $B_\infty^0$ be a fixed horoball based at $\infty$ which is a
little smaller than $B_\infty$, i.e. its boundary $\p B_\infty^0$
is a little higher than $\p B_\infty$. Let
$E_{i,j}^0=\d_{i,j}(X_i)\cap \p B_\infty^0$, and let
$N_R(Z_1)^0,..., N_R(Z_{d_i})^0$ be the part of
$N_R(Z_1),..., N_R(Z_{d_i})$ between the two
horizontal planes $\p B_\infty^0$ and $\p B_\infty$. Then
 $N_R(Z_1)^0\cap \p B_\infty,
N_R(Z_2)^0\cap \p B_\infty,..., N_R(Z_{d_i})^0\cap
\p B_{\infty}$ are Euclidean parallelograms contained in
$E_{i,j}$, containing the points $a_{i,j,1,0}$, $a_{i,j,2,n}$,
..., $a_{i,j,d_i,(d_i-1)n}$ as their topological  center points,
respectively, and they are isometric to $D_{i,j,k}, k=1, 2...,d_i$,
respectively. As $n$ is sufficiently large,
$N_R(Z_1)^0,..., N_R(Z_{d_i})^0$ are mutually far
apart from each other. We now take the
convex hull of the set $\{N_R(Z_1)^0$,...,
$N_R(Z_{d_i})^0\}$ in $\mathbb H^3$ and let $H_{i,j}^n$ be
the resulting convex manifold. Obviously $H_{i,j}^n$ is contained
in $\d_{i,j}(X_i)\cap B_\infty$.

Let $U_{i,j}$ be the part of $\d_{i,j}(X_i)$ between $E_{i,j}^0$
and $E_{i,j}$. Then $N_R(Z_1)^0,..., N_R(Z_{d_i})^0$
 are all contained in $U_{i,j}$, far apart
from each other. We now show

\begin{lem}\label{cap1}
If $n$ is sufficiently large, then  $H_{i,j}^n\cap  U_{i,j}=\{
N_R(Z_1)^0,..., N_R(Z_{d_i})^0\}$. \end{lem}

\pf  Let $F_k$ be the frontier of $N_R(Z_k)^0$ in
$U_{i,j}$, $1\leq k\leq d_i$. Then $F_k$ is contained in $\p
(N_R(\g_{i,j,k}(X_{i_*})))$. Since $N_R(\g_{i,j,k}(X_{i_*}))$ is
strictly convex and since $\p N_R((\g_{i,j,k}(X_{i_*})))$ is
smooth, then for any point $x\in F_k$, there is a unique geodesic plane
$P_x$ in $\mathbb H^3$ such that
$P_x \cap N_R(\g_{i,j,k}(X_{i_*}))=x$. Obviously $P_x$
 is not a vertical plane. Thus $P_x \cap \p B_\infty$ is a Euclidean circle
in $\p B_\infty$ with finite Euclidean diameter $a_x$. Since $F_k$
is compact, the set of numbers $\{a_x; x\in F_k\}$ has a finite
maximal value $a_k$. Let $a=max\{a_1,...a_{d_i}\}$, and let $c$ be
the maximal Euclidean diameter of the parallelograms
$\{N_R(Z_1)^0\cap E_{i,j}^0, ...,
N_R(Z_{d_i})^0\cap E_{i,j}^0\}$.
 By taking $n$ large enough, we can ensure that
$N_R(Z_1)^0\cap E_{i,j}^0, ..., N_R(Z_{d_i})^0\cap
E_{i,j}^0$ are mutually far apart from each other by Euclidean
distance at least $a+c$.  Then  the convex hull will
satisfy the condition $H_{i,j}^n \cap U_{i,j}=\{
 N_R(Z_1)^0, ..., N_R(Z_{d_i})^0\}$. \qed

 The manifold $H_{i,j}^n$ provided by Lemma \ref{cap1}
 is  the multi-$1$-handle  we
were seeking. (Figure \ref{cin} gives an illustration of
$H_{i,j}^n$ when $d_i=4$). We may assume the choice of $n$ works
in constructing all the multi-$1$-handles $H_{i,j}^n$, $j=1,...,
n_i$,  $i=1,2$.

We now glue the multi-$1$-handles $H_{i,j}^n$, $j=1,...,n_i$,
 to $J_i^-$ along $D_{i,j,k}, k=1,...,d_i$,
  $j=1,..., n_i$, (the gluing isometry should be clear).
 By our explicit construction,
  one can see that the hyperbolic structure
  on $H_{i,j}^n$ and the hyperbolic structure
  on $J_i^-$ match up after the gluing, forming a global
  (convex) hyperbolic structure.
  Thus we obtain  a compact, convex $3$-manifold $C_n(J_i^-)$.
 We also have a local isometry $g_i: C_n(J_i^-) \ra Y_i$,
 extending the local isometry $g_i:J_i^- \ra Y_i$.

\section{Strong separability in the free group}\label{csfg}

In this section, we present our main group theoretical result,
(Theorem \ref{hall}), which, together with the techniques used in
its proof, will have crucial applications in this paper.

 Let $S^-$ be a connected, compact, orientable
surface with genus $g$ and with $b>0$ boundary components. Fix a
point $s$ in $S^-$ as the base point, and let $F=\pi_1(S^-,s)$.
Then $F$ is a free group.  We may choose a free basis of $F$,
$$a_1,b_1,a_2,b_2,...,a_g,b_g,x_1,...,x_{b-1},$$ such that
$$x_1,x_2,...,x_{b-1},x_b=[a_1,b_1][a_2,b_2]\cdots
[a_g,b_g]x_1x_2\cdots x_{b-1}$$ are  represented by embedded loops
in $S^-$ (based at $s$) which are freely isotopic  to the $b$
boundary components of $S^-$ respectively. An element $\g$ of $F$
is \textit{peripheral} iff $\g$ is conjugate to
some power of some $x_i$. We prove the following:

\begin{thm} \label{hall}
Let $H \subset F$ be a finitely generated subgroup
 containing no nontrivial peripheral elements of $F$, and
let $y_1, ..., y_a\in F - H$.  Then there exists a
subgroup $G$ of $F$, with $|F:G| = m < \infty$,
 such that $G$ contains $H$ but does not
contain any elements of $\{y_1, ..., y_a, x_i, x_i^2,
...,x_i^{m-1}: i=1,...,b\}$.
\end{thm}

 In particular, the subgroup $H$ is separable;
 indeed, by M. Hall's Theorem, \textit{every} finitely
 generated subgroup of $F$ is separable.
 However, Theorem \ref{hall} gives much more information about $G$,
 since the number of
 elements to be separated is tied up with the index of the
 subgroup $G$ in $F$.

As an aside, we record a topological consequence of Theorem \ref{hall},
 which may be of independent interest.

\begin{cor}
Let $f:\alpha \ra S^-$ be an immersion of a geodesic loop in a
hyperbolic
 surface $S^-$ with $b > 0$ boundary components.  Then
 $f$ lifts to an embedding in a finite cover $\tilde{S}^- \ra S^-$,
 such that $\tilde{S}^-$ has exactly $b$ boundary components.
\qed\end{cor}

The proof of Theorem \ref{hall} is based on a technique, due
originally to Stallings and developed thoroughly in \cite{KM}, of
using folded graphs. We first need to recall some definitions, and
refer to \cite{KM} for more details. Let $L$ be a free basis for a
free group $F$, and let $L^{-1}$ be the set $\{x^{-1}; x\in L\}$).
An \textit{$L$-labeled directed graph} is a  graph such that each
edge of the graph is oriented, i.e.  with an initial vertex and a
terminal vertex assigned,  and  is  labeled with a unique element of
$L$. Given an $L$-labeled directed graph $\mathcal{G}$,
 we form an $L \cup L^{-1}$-labeled
graph $\widehat{\mathcal{G}}$ as follows:
 for each edge $e$ of $\cal G$-- say with
label $x$, initial vertex $v_1$ and terminal vertex $v_2$-- add a
new edge, denoted $e^{-1}$, with label $x^{-1}$, initial vertex
$v_2$ and terminal vertex $v_1$. The introduction of $\widehat{\cal
G}$ is purely for technique convenience.

An $L$-labeled directed graph $\mathcal{G}$ is said to be
\textit{$L$-regular} if, for every vertex $v$ of $\mathcal{G}$ and
every $x \in L \cup L^{-1}$, there is exactly one edge of
$\widehat{\mathcal{G}}$ with initial vertex $v$ and with label $x$.
An $L$-labeled directed graph $\mathcal{G}$ is called
\textit{folded} if there is no pair of distinct edges $e,
e^{\prime}$ in
 $\widehat{\mathcal{G}}$ with the same initial vertex and the same label.
 Obviously a regular graph is folded.

 If $\cal G$ is folded, then every  reduced  path
 (i.e. path containing no subpath of the form $e, e^{-1}$)
   in $\widehat{\mathcal{G}}$ determines a unique
freely reduced word in $L \cup L^{-1}$, and thus a unique element of
$F$. If we fix a vertex $v_0 \in \mathcal{G}$, then the set of all
elements of $F$ corresponding to the set of  reduced loops in
$\widehat{\mathcal{G}}$ based at $v_0$ is a subgroup of $F$, denoted
$L(\mathcal{G}, v_0)$. A proof of the following lemma is contained
in \cite{KM}.

\begin{lem} \label{km}
If $\mathcal{G}$ is a  finite and $L$-regular  graph, then
$L(\mathcal{G},v_0)$
 is a finite-index subgroup of $F$, and its index in $F$
 is equal to the number of vertices in $\mathcal{G}$.
\end{lem}

An example of an $L$ directed graph ${\cal G}_0$ is the wedge of
$|L|$ circles each given some fixed orientation and labeled with the
labels of $L$, one each. If we denote the vertex by $v_0$, then
$L({\cal G}_0, v_0)=F$. The point of Lemma \ref{km} is that a graph
${\cal G}$ as given in the lemma is naturally a finite sheeted
covering of the graph  ${\cal G}_0$ with degree equal to the number
of vertices of ${\cal G}$.

We now proceed to prove Theorem \ref{hall}. From now on in this
section $L$ denotes the free basis of the free group
$F=\pi_1(S^-,s)$ given at the beginning of this section.
 We may certainly
assume that $F$ is not a cyclic group, and thus we have either
$g>0$, or $g=0$ and $b>2$. Elements of $F$ will be considered as
words in letters in $L\cup L^{-1}$. It follows directly from the
proof of Hall's Theorem in \cite{KM} that  there is a connected,
finite, folded, $L$-labeled directed graph $\mathcal{G}_0$, with
base vertex $v_0$, such that $L(\mathcal{G}_0,v_0) = H$, and the
words
 $y_1, ..., y_a$ are representable by non-closed paths in
 $\widehat{\mathcal{G}}_0$ with the base vertex $v_0$ as their initial vertex.
Also, no loop in $\widehat{\mathcal{G}}_0$ (based at any vertex)
 represents a non-zero
 power of any $x_i$, for otherwise $H$ would contain nontrivial peripheral elements.
Note that ${\cal G}_0$ is the quotient of the minimal $H$-invariant
subtree of the Cayley graph of $F$ with respect to the given
generators.

We need some more definitions.
Suppose that $\mathcal{G}$ is a
 finite, connected, $L$-labeled directed graph.
 For each $i=1,...,b$, we call a  path in
 $\widehat{\mathcal{G}}$ an $x_i$-\textit{path}
 if it represents  a subword of the word $x_i^k$ for some
 non-negative integer $k$.
 A single vertex of the graph is also considered as an $x_i$-path,
corresponding to the empty subword of $x_i$. An $x_i$-path is
called an $x_i$-\textit{loop} if it is a loop representing the word $x_i^k$
for some positive integer $k$.
 An $x_i$-path is
called {\it maximal} if  it is not contained in any other
$x_i$-path besides itself. Now suppose further that ${\cal G}$ is
folded and $\widehat{\cal G}$ contains no $x_i$-loops, for all
$i=1,...,b$.
 Then for each $i=1,...,b$,  every $x_i$-path is contained
in a unique maximal $x_i$-path with finite length (where the
length of a path is the number of edges that the path contains).
If $ i<b$, then any maximal $x_i$-path is an embedded path, and
any two different maximal $x_i$-paths are disjoint. For a maximal
$x_b$-path, every oriented edge in the path appears only once
 in the path but the path may
cross  itself at some common vertices. Any two different maximal
$x_b$-paths have disjoint oriented  edges but may cross each other
at some common vertices.
 It follows that
there are only finitely many maximal $x_i$-paths in $\widehat{\cal G}$,
 which we denote
 by $C_{i,j}$, $i=1,...,b$, $j=1,...,m_i$. For a maximal $x_i$-path
 $C_{i,j}$,
 its initial (respectively terminal) vertex  is missing an incoming
 (respectively outgoing) edge
 whose label is  the predecessor (respectively successor) to
 the first (respectively the last) label
 of $C_{i,j}$, where $C_{i,j}$ is considered as a subword of the word
 $x_i^k$ (for some $k\geq 0$).
 We shall call these two missing labels the \textit{initial} and
 \textit{terminal} missing
 labels of $C_{i,j}$ respectively.
 Of course if $i<b$, then the initial or terminal missing label for
 every $C_{i,j}$ is always $x_i$.
 Note that  for a maximal $x_b$-path $C_{b,j}$,
  if the first label  of $C_{b,j}$ is the  letter $a_1$,
  then the initial missing label of $C_{b,j}$ is the  letter  $x_{b-1}$
  if $b>1$ or the latter $b_g^{-1}$ if $b=1$;
    and similarly if the last label  of $C_{b,j}$ is the  letter  $b_g^{-1}$,
  then the terminal  missing label of $C_{b,j}$ is the  letter $x_1$ if $b>1$
  or the letter $a_1$ if $b=1$.

\begin{lem}\label{pair1}
Let ${\cal G}$ be a finite, connected, $L$-labeled, directed graph
such that ${\cal G}$ is folded and such that $\widehat{\cal G}$
contains no $x_i$-loops for any $i=1,...,b$. Then $x\in L$ is the
initial missing label of some maximal $x_i$-path $C_{i,j}$ if and
only if $x$ is the terminal missing label for some maximal
$x_i$-path $C_{i,j'}$.
\end{lem}

\pf Suppose that the number of vertices of ${\cal G}$ is $m$. Let
$k$ be the number of existing directed edges of ${\cal G}$ with
label $x$. Then $m-k$ is equal to the number of initial missing
edges of $\widehat{\cal G}$ with label $x$ and is also equal to
the number of terminal missing edges of $\widehat{\cal G}$ with
label $x$. The lemma follows. \qed

We shall let
  $$L_*= \{a_1,b_1,...,a_g,b_g\}.$$
The proof of the following lemma is obvious.

\begin{lem}\label{pair2}
Let ${\cal G}$ be a finite connected $L$-labeled directed graph
such that ${\cal G}$ is folded and such that $\widehat{\cal G}$
contains no $x_i$-loops for any $i=1,...,b$.
 If $x\in L_*\cup L_*^{-1}$ is the initial or terminal
missing label of some $C_{b,j}$ at a vertex $v$, then $x^{-1}$
must also be the terminal or initial missing label of some
$C_{b,j'}$ respectively at the same vertex $v$.\qed
\end{lem}

 By Lemma \ref{km}, it is enough to show that the graph
 $\mathcal{G}_0$ embeds in a finite
 $L$-regular graph $\mathcal{G}_*$
 such that $\widehat{\mathcal{G}}_*$ contains no $x_i$-loop
 (based at any vertex)
 representing the word $x_i^k$ for any $i=1,...,b$ and
 $k=1,...,m_*-1$, where $m_*$ is the number of vertices of ${\cal G}_*$.
 Indeed, assuming such ${\cal G}_*$ is found,  we have
 \newline
 (1) $G=L({\cal G}_*, v_0)$ is an index $m_*$ subgroup of $F$ (Lemma \ref{km});
  \newline
  (2) $G$ contains $H$ as a subgroup but does not contain any of
 the elements $y_1,...,y_a$ (because ${\cal G}_0$ is an embedded subgraph
 of  ${\cal G}_*$ and $y_1,...,y_a$ are represented by non-closed paths
 with initial vertex $v_0$);
and \newline
 (3)  $G$ does not contain
 any of the elements $x_i^k$, $i=1,...,b$, $k=1,...,m_*-1$ (because $\widehat{\cal G}_*$
 has no $x_i$-loop
 representing the word $x_i^k$ for any $i=1,...,b$ and
 any $k=1,...,m_*-1$).

In the rest of this section we  show that such a graph ${\cal G}_*$
exists.

\noindent{\bf Definition}. Let ${\cal G}$ be a  finite, connected,
$L$-labeled, directed graph such that ${\cal G}$ is folded and such that
$\widehat{\cal G}$ contains no $x_i$-loops for any $i=1,...,b$. A
graph  ${\cal G}'$ is called a {\it good extension} of ${\cal G}$
if  \\ (1) ${\cal G}'$ is  a finite, connected, $L$-labeled, directed
graph; \\ (2)  ${\cal G}'$ contains ${\cal G}$ as an embedded
subgraph;\\ (3) ${\cal G}'$ is folded;\\ (4) $\widehat{\cal G}'$
contains no $x_i$-loops for all $i=1,...,b$.

 \noindent{\bf Definition}. Let ${\cal G}$
be a  finite, connected, $L$-labeled, directed graph such that ${\cal
G}$ is folded and such that $\widehat{\cal G}$ contains no $x_i$-loops
for any $i=1,...,b$. A graph  ${\cal G}'$ is called a {\it perfect
extension} of ${\cal G}$ if  \\ (1) ${\cal G}'$ is a finite,
connected, $L$-labeled, directed graph; \\ (2)  $ {\cal G}'$
contains ${\cal G}$ as an embedded subgraph;\\ (3) ${\cal G}'$ is
$L$-regular;\\ (4)  $\widehat{\cal G}'$ contains no loop
representing the word $x_i^k$ for any $i=1,...,b$ and
any $k=1,...,m-1$, where $m$ is the number of vertices of ${\cal G}'$.

We shall describe a canonical procedure for constructing  a finite
sequence of  graphs ${\cal G}_0,{\cal G}_1, ...,{\cal G}_n$ such
that\\ (1) ${\cal G}_0$ is the graph given above;\\ (2) ${\cal
G}_{p+1}$ is a good extension of ${\cal G}_p$ for each $p=1,...,
n-2$ (if $n>1$);
\\(3) ${\cal G}_n$ is a perfect extension of ${\cal G}_{n-1}$.
\\Obviously if such a sequence of graphs can be constructed, then
 ${\cal G}_n$ will be the graph which we seek.

We divide our discussion  into three cases: $b=1$, $b=2$ and
$b>2$. We need the following definitions for all the three cases.

\noindent{\bf Definitions.} For $g>0$, let ${\cal G}$ be a finite,
connected, $L$-labeled, directed graph such that ${\cal G}$ is
folded and such that $\widehat{\cal G}$ contains no $x_i$-loops for
any $i=1,...,b$.
 Let $x \in L_* \cup L_*^{-1}$.
\\A maximal $x_b$-path $C_{b,j}$ of $\widehat{\cal G}$ is called a
{\it type $I$} maximal $x_b$-path \textit{with missing label} $x$ if
\\ (1) $x$ is the initial missing label of $C_{b,j}$ and $x^{-1}$
is the terminal missing label of the same path $C_{b,j}$; and\\
(2) the initial vertex and the terminal vertex of $C_{b,j}$ are
the same vertex.\\
 A
maximal $x_b$-path $C_{b,j}$ of $\widehat{\cal G}$ is called a
{\it type $II$} maximal $x_b$-path \textit{with missing label} $x$
 (again we assume that $g>0$) if  $x$ is the
initial missing label of $C_{b,j}$ and  is also the terminal
missing label of the same path.

\begin{figure}
{\epsfxsize=5in \centerline{\epsfbox{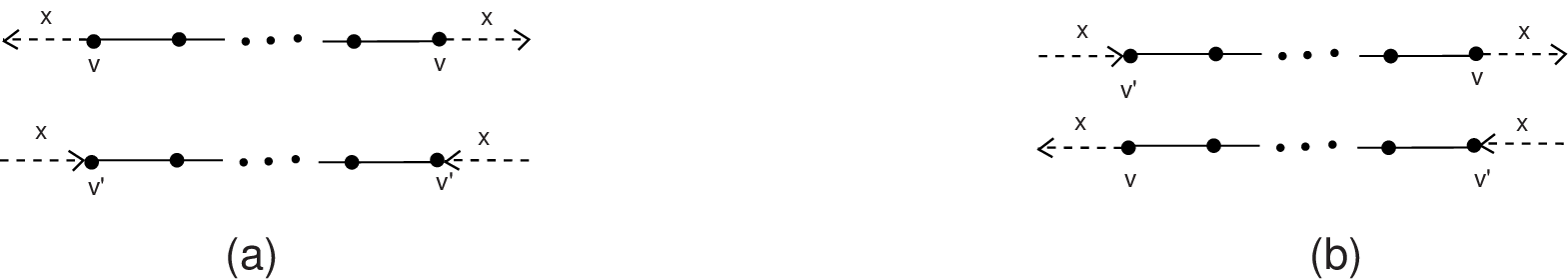}}\hspace{10mm}}
\caption{}\label{fg2}
\end{figure}

 Figure \ref{fg2} (a) illustrates  a pair of type I maximal
$x_b$-paths, and  Figure \ref{fg2} (b) shows a pair of type II
maximal $x_b$-paths. In these figures, a missing label is
represented
 by a dotted, labeled edge;
  the initial missing label is given at the
left end of a path and the terminal one at the right end.

\noindent{\bf Case 1}. $b=1$.

In this case, we have $g>0$, free basis
$L=L_*=\{a_1,b_1,...,a_g,b_g\}$, and a surface $S^-$
 with a single boundary component, which is freely isotopic to a loop
representing the commutator $x_b=x_1=[a_1,b_1]\cdots [a_g,b_g]$.

\begin{figure}[!ht]
{\epsfxsize=4in \centerline{\epsfbox{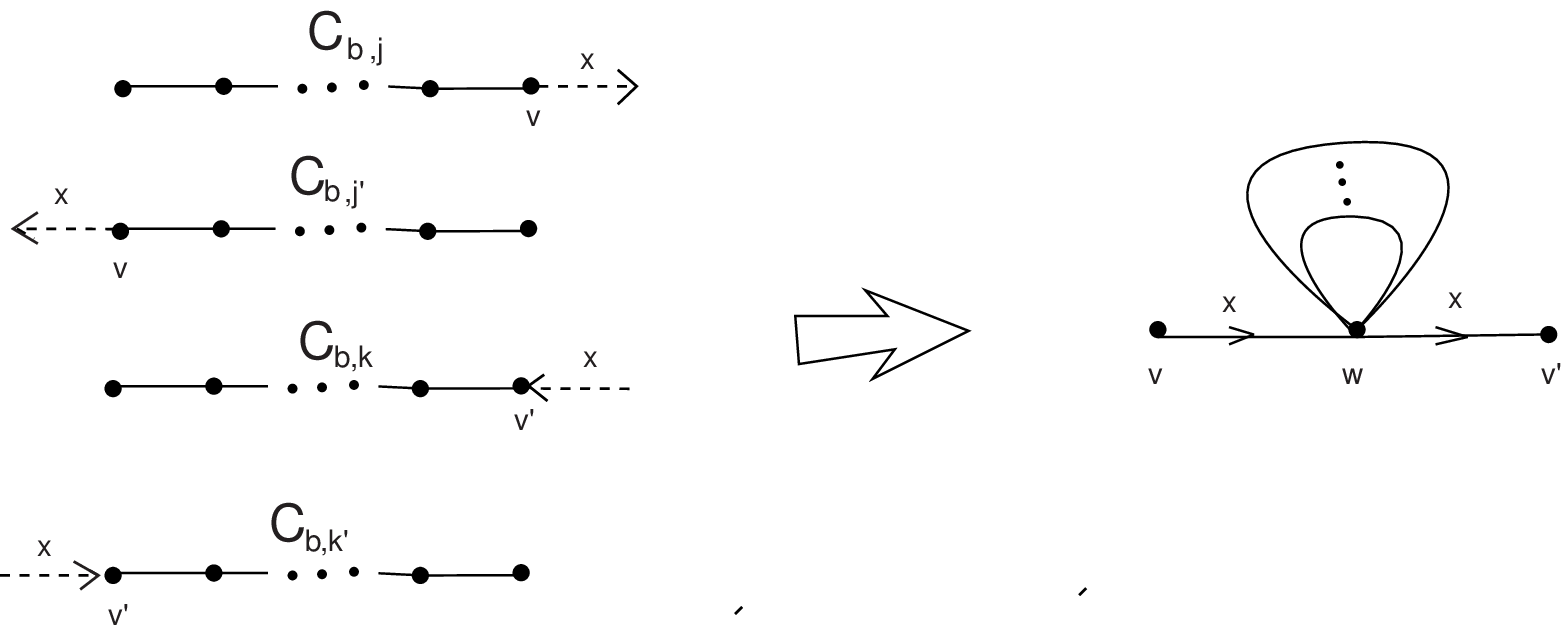}}\hspace{10mm}}
\caption{}\label{fig1}
\end{figure}

 Start with the graph ${\cal G}_0$. Since $\widehat{\cal G}_0$  has no $x_b$-loops,
 every maximal $x_b$-path in $\widehat{\cal G}_0$ has
 both initial and terminal missing labels.
  Suppose that $x$ is a missing label.
Then $x\in L_*\cup L_*^{-1}$. By Lemmas \ref{pair1}
and \ref{pair2}, we have maximal $x_b$-paths (possibly
non-distinct) $C_{b,j}, C_{b,j^{\prime}}, C_{b,k}, C_{b,k^{\prime}}$,
  with missing labels as illustrated on the left hand side
in Figure \ref{fig1}.   Note that although we draw these paths
separately, they may actually share some common vertices. Also a
path we draw may not be simply connected, i.e. some of its vertices
maybe the same vertex. Note also that instead of drawing an oriented
edge with label $x^{-1}$, we often draw, equivalently, an edge with
the opposite orientation and with label $x$ (i.e. every edge we draw
shall be considered as two directed edges with opposite orientations
and inverse labels). Also, a pair of paths in the figure may in fact
be non-distinct (e.g. maybe $C_{b,j} = C_{b,j^{\prime}}$).  However,
the four paths given on the left hand side  in Figure \ref{fig1}
must satisfy
 $C_{b,j}\ne C_{b,k}$ and $C_{b,j'}\ne C_{b,k'}$.

\noindent {\bf Operation 1}. Suppose that the four maximal
$x_b$-paths given in Figure \ref{fig1} satisfy: $C_{b,j}\ne
C_{b,j'}$, $C_{b,k}\ne C_{b,k'}$, and either $C_{b,j}\ne C_{b,k'}$ or
$C_{b,j'}\ne C_{b,k}$. Then we may perform the following operation:
add a new vertex $w$; for each letter in $L-\{x\}$ add a new  edge
with both its initial and terminal vertices at $w$ and with that
letter as the label; add a new  edge with label $x$, initial vertex
 $v$ and terminal vertex $w$; and add a
new  edge with label $x$, initial vertex $w$ and
terminal vertex $v'$. A piece of the resulting graph is
shown on the right hand side in Figure \ref{fig1}, where there are
exactly $2g-1$ single-edge-loops at the new vertex $w$, with labels
one each from $L-\{x\}$.

We claim that the resulting graph is a good extension of ${\cal
G}_0$, and that the number of maximal $x_b$-paths in the new graph is
reduced. Indeed, the new graph  is obviously finite, connected,
$L$-labeled and contains ${\cal G}_0$ as an embedded subgraph,
 and one can check
 that it is folded. One can also check that the number of maximal
$x_b$-paths in the new graph is reduced. Namely the
 maximal $x_b$-paths $C_{b,j}$ and $C_{b,j'}$ are joined
into a single maximal $x_b$-path, as are $C_{b,k}$ and $C_{b,k'}$.
 No extra maximal $x_b$-paths are created, since all the added edges
are used in the two new maximal $x_b$-paths. In particular no
$x_b$-loops are created.

\noindent {\bf Operation 2}. Suppose that the four maximal
$x_b$-paths given in Figure \ref{fig1} satisfy:
 $C_{b,j} \neq C_{b,k^{\prime}}, C_{b,j^{\prime}} \neq C_{b,k}$
 and either $C_{b,j} \neq C_{b,j^{\prime}}$ or $C_{b,k} \neq C_{b,k^{\prime}}$.
 Then we may perform the following
operation: add a new edge $e$ with the label $x$, initial vertex
$v$ and terminal vertex $v^{\prime}$.
 Again one can easily check that the resulting graph from an operation 2
 is a good extension of the old graph, and the number of maximal
$x_b$-paths in the new graph is reduced.

For each pair $\{x,x^{-1}\}$, we perform Operations 1 and 2
 as many times as possible.  Since each operations
reduces the number of maximal $x_b$-paths, this process will
terminate in a graph ${\cal G}_1$ for which neither
Operation 1 nor Operation 2 may be applied for any pair
$\{x,x^{-1}\}$.

\begin{lem}\label{3cases}(1) ${\cal G}_1$ is a good extension
 of ${\cal G}_0$;\\
 (2) In $\widehat{\cal G}_1$, for each letter pair $\{x,x^{-1}\}$, we have either

\noindent I. every  maximal $x_b$-path with $x$ or $x^{-1}$ as an
initial or terminal missing label is of Type I;

\noindent II. there are only two  maximal $x_b$-paths which have
$x$ or $x^{-1}$ as a missing label; moreover one of the two
paths is a type II maximal $x_b$-path with $x$ as its missing
label and the other path is a type II maximal $x_b$-path with
 $x^{-1}$ as a missing label; or

\noindent III. there is no maximal $x_b$-path with $x$ or
$x^{-1}$ as an initial or terminal missing label.
\end{lem}

\pf Part (1) of the lemma holds since we have checked that
Operations 1 and 2 always yield good extensions.

To show part (2) of the lemma,  fix a letter
 pair $\{x, x^{-1}\}$, and suppose that case III does not
happen. Then as discussed above, in $\widehat{\cal G}_1$ we have
maximal $x_b$-paths $C_{b,j}$,
$C_{b,j'}$, $C_{b,k}$ and $C_{b,k'}$ (possibly non-distinct) as shown in Figure
\ref{fig1}. If $C_{b,j}=C_{b,j'}$, then we must have
$C_{b,k}=C_{b,k'}$, and vice versa. For otherwise Operation 2
would apply. It follows that if $C_{b,j} = C_{b,j'}$
 or  $C_{b,k} = C_{b,k'}$ are the same path, then
any other such pair of paths must be the same, i.e.  Case
I holds. So we may assume that for any two pairs
 of maximal $x_b$-paths
 $(C_{b,j}, C_{b,j'})$ and $(C_{b,k}, C_{b,k'})$ as shown in Figure \ref{fig1},
 we have   $C_{b,j}\ne C_{b,j'}$ and $C_{b,k}\ne C_{b,k'}$.
 But since Operation 1 does not apply to them, we must have
$C_{b,j}= C_{b,k'}$ and $C_{b,j'}=C_{b,k}$. Thus we have a pair of
Type II maximal $x_b$-paths with $x$ and $x^{-1}$ as their missing
labels, respectively. Again since Operation 1 does not apply, there
is at most one such pair, i.e. Case II occurs.
 \qed

\begin{figure}
{\epsfxsize=5in \centerline{\epsfbox{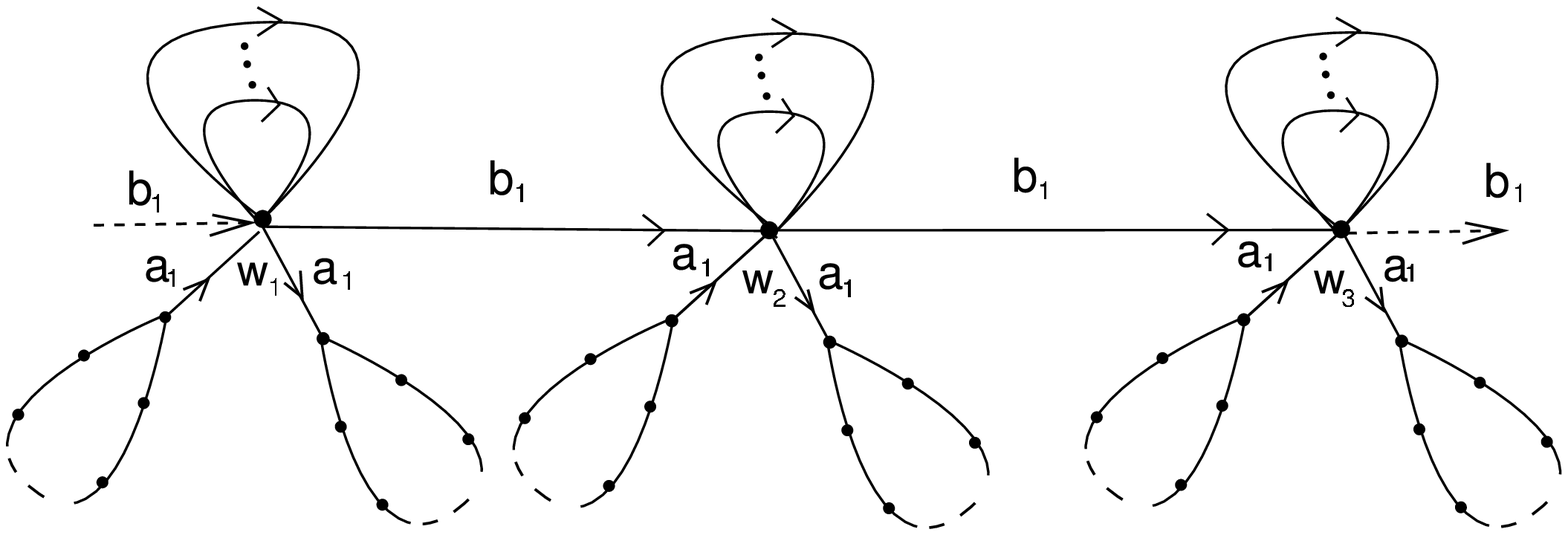}}\hspace{10mm}}
\caption{}\label{fg1}
\end{figure}

\noindent {\bf Operation 3}.
 Suppose that $\widehat{\cal G}_1$ satisfies Property I of
 Lemma \ref{3cases}, for some letter pair $\{x, x^{-1} \}$.
 Let $k$ be the total number of pairs of Type I maximal $x_b$-paths
 with missing labels in $\{x, x^{-1}\}$.
 Then we may perform the operation shown in
Figure \ref{fg1} to change all of them into a single pair of Type
II maximal $x_b$-paths with a letter different from $x$ or
$x^{-1}$ as the missing label. For illustration, in Figure
\ref{fg1}  we assume that $k=3$ and $x=a_1$. We
add three new vertices $w_1,w_2,w_3$, two new edges with label
$b_1$, and new loops at each $w_i$ having labels one each from
$L_*-\{a_1,b_1\}$. The three pairs of type I maximal $x_b$-paths
with missing label $a_1$ and $a_1^{-1}$ become one pair of type II
maximal $x_b$-paths with missing label $b_1$ and $b_1^{-1}$. All
the added new edges are used in this pair of new maximal
$x_b$-paths.  If the missing labels for the type I pairs are $b_1$
and $b_1^{-1}$ instead of $a_1$ and $a_1^{-1}$, then in Figure
\ref{fg1} we exchange the letters $a_1$ and $b_1$.
 In general, for
arbitrary $k$ and arbitrary missing label pair
$\{x,x^{-1}\}$, it should be clear how to make a similar
operation.

 Again one can check that, if Operation 3 is applied to a good
 extension,  then
 the resulting graph is a good extension. Also,
 after Operation 3 is applied, the $k$ pairs of Type I paths are
 replaced with a single pair of Type II paths, whose
 missing labels are different from $\{x, x^{-1} \}$.
 If $\{x, x^{-1} \} = \{a_i,a_i^{-1}\}$, then the
 new label pair is $\{b_i,b_i^{-1}\}$, and vice versa.
 Note that after applying Operation 3, the total
 number of maximal $x_b$-paths cannot increase,
 although it may stay the same.

 Now, starting with the graph ${\cal G}_1$,
  we perform Operation 3   as many times as possible.
 We are left with a graph, ${\cal G}_2$,
 whose maximal $x_b$-paths are all of Type II.
 Suppose ${\cal G}_2$ has more than one
 pair of Type II maximal $x_b$-paths with
 missing labels in $\{x, x^{-1} \}$.
 Then we may apply Operation 2.
 The effect is to replace a pair of Type II
 paths with a single Type II path.
 Therefore, after performing Operation 2
 repeatedly, we arrive at a graph ${\cal G}_3$,
 such that, for every missing label  pair $\{x,x^{-1}\}$,
 there is exactly one pair of Type II paths
 corresponding to that pair, and there are no other maximal $x_b$-paths
 with $x$ or $x^{-1}$ as missing label.
  Note that ${\cal G}_3$ cannot be $L$-regular,
 since our operations thus far never create an
 $x_b$-loop.

\begin{figure}[!ht]
{\epsfxsize=5in \centerline{\epsfbox{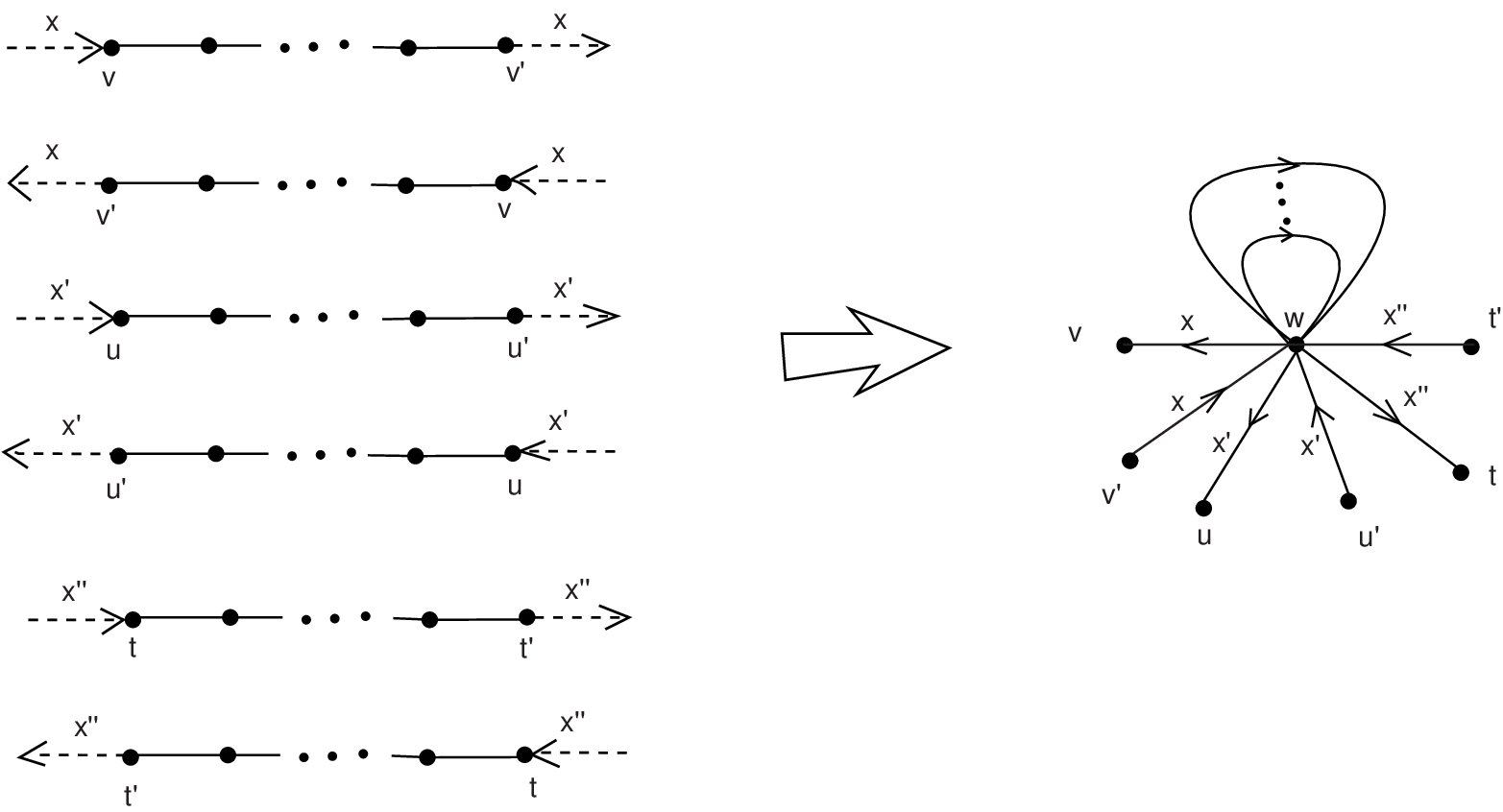}}\hspace{10mm}}
\caption{}\label{fig3}
\end{figure}

\noindent{\bf Operation 4}. We add a single
 vertex $w$, and then add the appropriate edges  to
 make the graph folded, with no missing labels
 from $L_*$.
  This is illustrated in Figure \ref{fig3} in the case
 where there are exactly three pairs of maximal $x_b$-paths,
 with missing labels
  $\{x,x^{-1}\}$, $\{x',x^{'-1}\}$, $\{x'',x^{''-1}\}$.
 The loops at the added  vertex $w$ have labels  from
 the remaining labels in $L_*$.

The resulting graph ${\cal G}_4$ is what we wanted, i.e. it is a
perfect extension of ${\cal G}_3$. Indeed, it is easy to check
that the graph is finite, connected, and folded, with ${\cal G}_3$
as an embedded subgraph; also, in the current case,
 $L_* = L$, so  ${\cal G}_4$ is $L$-regular.
 So we only need to check that $\widehat{\cal G}_4$
has no $x_b$-loops representing the word $x_b^k$ for any
$k=1,...,m-1$, where $m$ is the number of vertices of ${\cal
G}_4$. To see this holds, refer to Figure \ref{fig3}, and trace
out an $x_b$-loop, starting at any vertex of ${\cal G}_4$. One
sees that an $x_b$-loop does not occur until every edge of
$\widehat{\cal G}_4$ has been used exactly once (in particular
each maximal $x_b$-path $C_{b,j}$ in $\widehat{\cal G}_3$ has been
traced out exactly once). This $x_b$-loop
represents the word $x_b^m$, since ${\cal G}_4$ is $L$-regular and
contains $m$ vertices.
 This completes the proof of Theorem \ref{hall} when $b=1$.

\noindent{\bf Case 2}. $b=2$.

In this case,  the free basis $L=\{a_1,b_1,...,a_g,b_g, x_1\}$ and
$x_b=x_2=[a_1,b_1]\cdots [a_g,b_g]x_1$ (note that $g>0$ in this
case too). Recall $L_*=\{a_1,b_1,...,a_g,b_g\}$.

The first step in this  case is to construct a good
extension graph ${\cal G}_1$ of ${\cal G}_0$ such that
$\widehat{\cal G}_1$ has no missing label pairs from $L_*\cup
L_*^{-1}$. The procedure for constructing such ${\cal G}_1$
 is the
same as that given in Case 1. Only in the current case, after
applying an Operation 1, or 3, or $4$, we may increase the number
of maximal $x_b$-paths with missing label $x_1$ and may also
 increase the number of maximal $x_1$-paths.
But no $x_1$-loops will be created  because during each of these
operations no new edge with label $x_1$ is added.

\begin{figure}[!ht]
{\epsfxsize=5in \centerline{\epsfbox{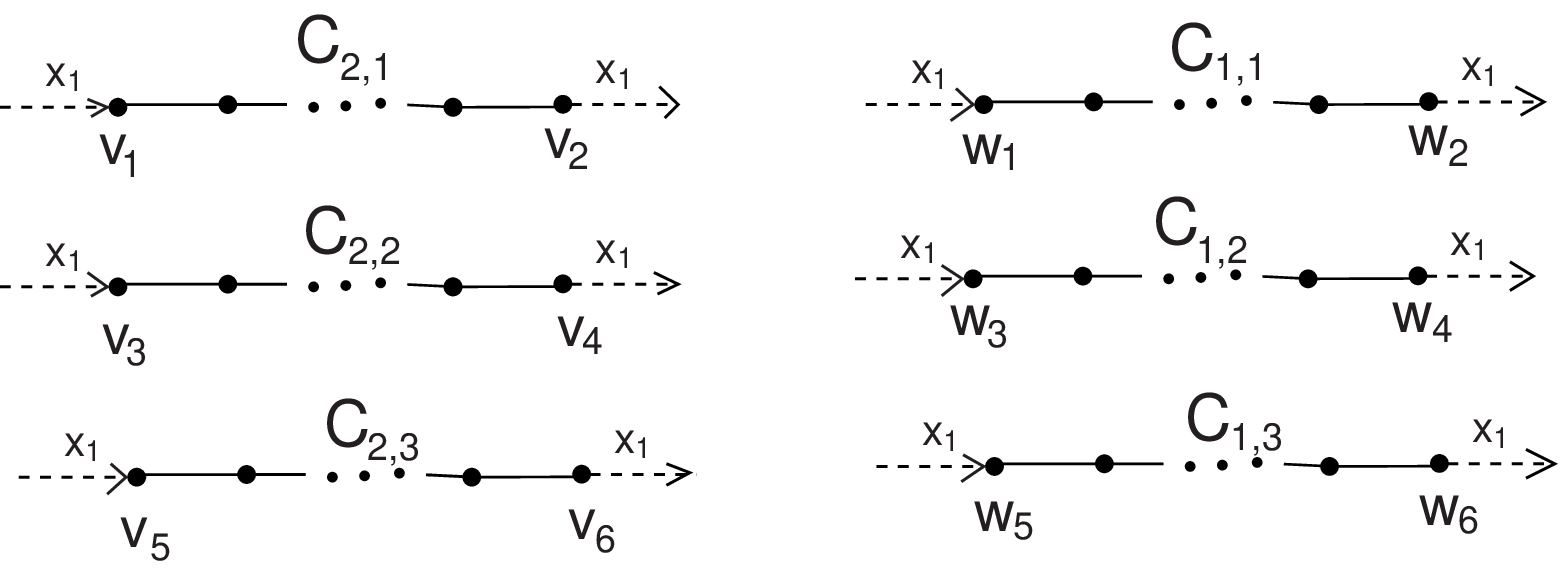}}\hspace{10mm}}
\caption{}\label{fig4}
\end{figure}

\begin{figure}[!ht]
{\epsfxsize=1.5in \centerline{\epsfbox{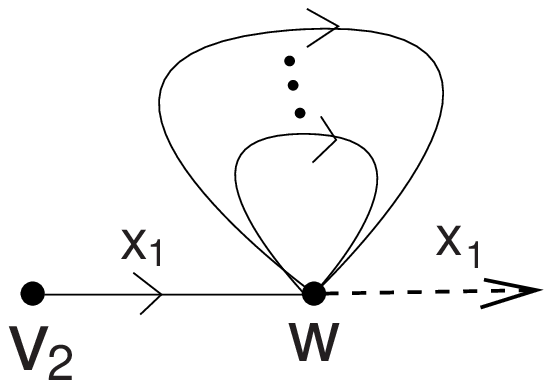}}\hspace{10mm}}
\caption{}\label{fig5}
\end{figure}

Let ${\cal G}_1$ be the resulting good extension graph after all
missing edges with labels from $L_*$ are eliminated.  So every
maximal $x_{b}$-path $C_{b,j}$ in $\widehat{\cal G}_1$ has both
its initial and terminal missing labels being $x_1$. In the
current case, we also need to consider maximal $x_1$-paths. Their
missing labels are of course always $x_1$. Note also that in
$\widehat{\cal G}_1$ the number of maximal $x_b$-paths is the same
as the number of maximal $x_1$-paths (since $b=2$ and since $x_1$ is now
the only missing label from $L$). Suppose that there are at least
three maximal $x_b$-paths in $\widehat{\cal G}_1$. We illustrate
three such paths in Figure \ref{fig4} left hand side. Thus there
are at least three maximal $x_1$-paths, which are shown in Figure
\ref{fig4} right hand side. Since the graph ${\cal G}_1$ is
connected and folded, the vertices shown in Figure \ref{fig4}
satisfy $v_1\ne v_3\ne v_5\ne v_1$, $v_2\ne v_4\ne v_6\ne v_2$.
Also  $C_{1,1}$, $C_{1,2}$, $C_{1,3}$ are mutually disjoint
embedded paths. By adding to $v_2$ a subgraph which is shown in
Figure \ref{fig5} (in the figure, the loops have labels one each
from $L_*$), we may assume that $v_2\ne v_3$ and $v_2\ne v_5$. We
may also assume that $v_2=w_2$, and $v_3\ne w_1$. Now we simply
add an edge with label $x_1$ to ${\cal G}_1$ pointing from $v_2$
to $v_3$. Then the new graph is folded with one less number of
maximal $x_b$-paths, and with no $x_i$-loops created, $i=1,2$
(since $v_3$ is not an end vertex in the path $C_{1,1}$).

\begin{figure}[!ht]
{\epsfxsize=5in \centerline{\epsfbox{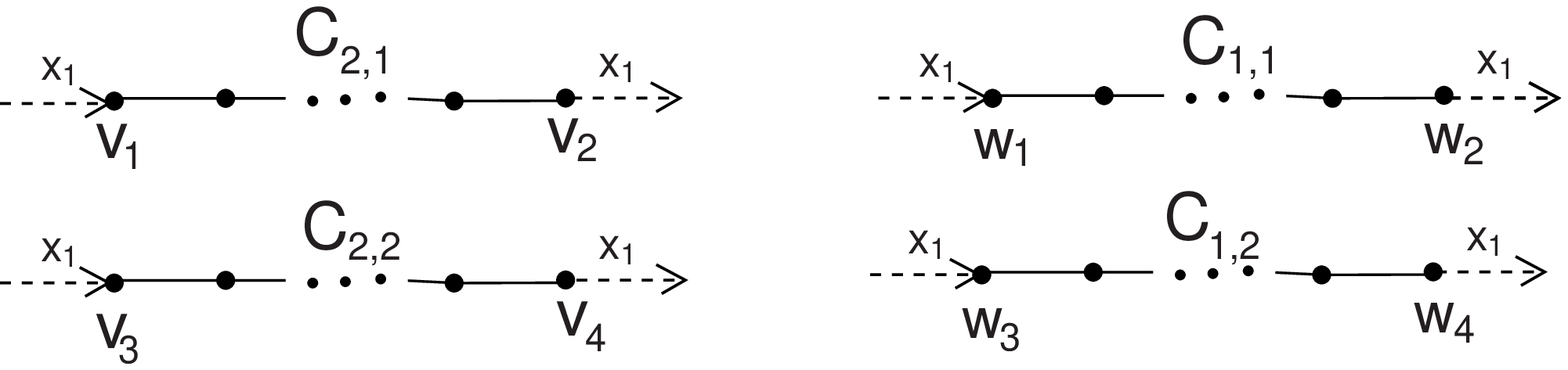}}\hspace{10mm}}
\caption{}\label{fig6}
\end{figure}

\begin{figure}[!ht]
{\epsfxsize=3in \centerline{\epsfbox{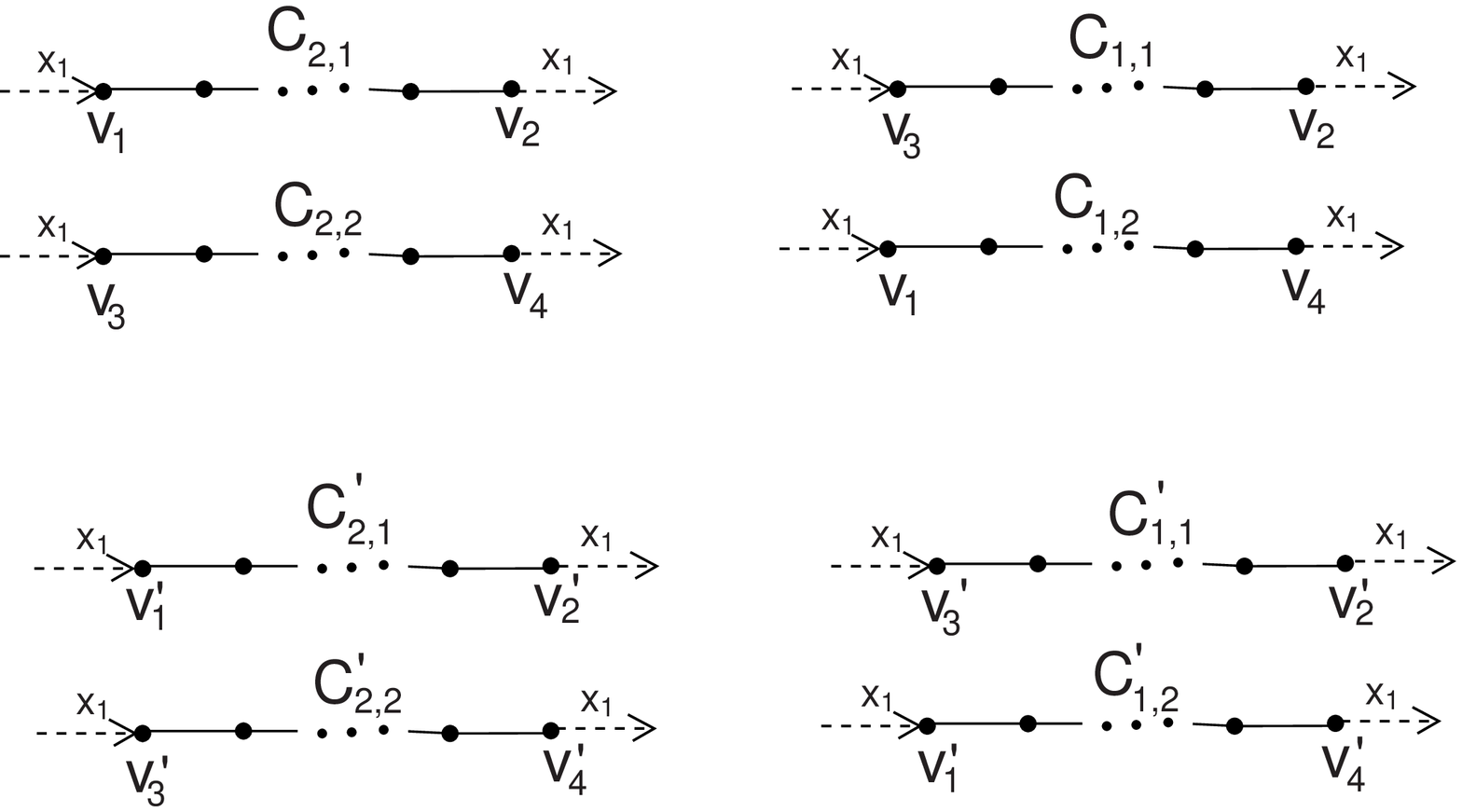}}\hspace{10mm}}
\caption{}\label{fig7}
\end{figure}

So we may assume we have a graph, still denoted  ${\cal G}_1$,
such that $\widehat{\cal G}_1$  has at most
 two maximal $x_b$-paths.
 If there is only one such path, we simply add  a single edge
 with the missing label $x_1$. The resulting graph  is
 what we requested.
 So we may assume that there are exactly two such paths, as shown in
 Figure \ref{fig6}.
 Again we
have $v_1\ne v_3$, $v_2\ne v_4$ and may assume $v_2\ne v_3$ and
$v_2=w_2$. If $v_3\ne w_1$, then we simply add an edge
  with label
$x_1$ pointing from $v_2$ to $v_3$. So we may assume that
$v_3=w_1$. In such case we cannot add an edge from $v_2$ to $v_3$
with label $x_1$ since that will create a graph which has an
$x_1$-loop but is not yet $L$-regular. The initial and terminal
vertices of these paths are illustrated in the first two rows of
Figure \ref{fig7}. Now we take an identical copy ${\cal G}_1'$ of
the graph ${\cal G}_1$. The corresponding maximal $x_i$-paths in
$\widehat{\cal G}_1'$ are denoted by $C_{i,j}'$. The union of
these paths from the two graphs are shown in Figure \ref{fig7}. We
now connect these two graphs together as follows: add an edge from
$v_2$ to $v_1'$, add an edge  from $v_4$ to $v_3'$, add an edge
from $v_2'$ to $v_3$, and add an edge from $v_4'$ to $v_1$, all
with label $x_1$. Then one can easily check that the resulting
graph  is what we requested. This proves the theorem when $b=2$.

\noindent {\bf Case 3}. $b\geq 3$.

Again we first eliminate all missing labels belonging to
$L_*=\{a_1,b_1, ...,a_g,b_g\}$, with a similar method as in Case
2. Namely during the relevant operations, no new edges with label
$x_i$, $i=1,...,b-1$, are added. Note that the resulting graph
${\cal G}_1$ must have missing $x_i$-labels for each $i=1,...,b-1$
(since no $x_i$-loops were created). Let $C_{b,j}$ be a maximal
$x_b$-path in $\widehat{\cal G}_1$. Its initial missing label
might not be the same as its terminal missing label. If so, we
call such path a maximal $x_b$-path \textit{with mixed missing labels}.
Note that for each $i=1,...,b-1$, the number of initial missing
$x_i$ labels is equal to the number of terminal missing $x_i$
labels. It follows that if the graph $\widehat{\cal G}_1$ has a
maximal $x_b$-path $C_{b,j}$ with mixed missing labels and with
$x_1$, say,  as the terminal missing label, then there must be
another maximal $x_b$-path $C_{b,j'}$ with mixed missing labels
and with $x_1$ as the initial missing label, and vice versa.

\begin{figure}[!ht]
{\epsfxsize=1.5in \centerline{\epsfbox{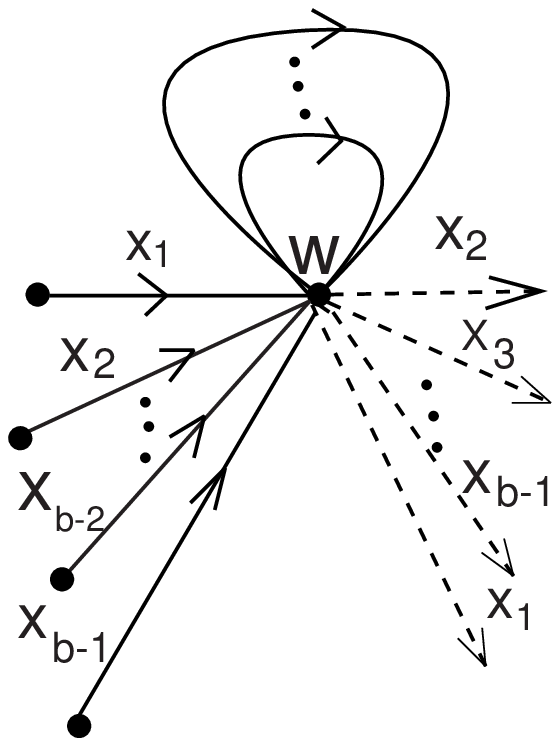}}\hspace{10mm}}
\caption{}\label{fig8}
\end{figure}

Take $b-1$ maximal $x_b$-paths  with mutually distinct terminal
missing labels. Such $b-1$ paths must exist as we already noted.
 By adding a subgraph as shown in Figure \ref{fig8}
to the $b-1$ terminal vertices of the  maximal $x_b$-paths (in
Figure \ref{fig8}, loops at the vertex $w$ have labels from $L_*$)
we may assume that the terminal vertex is different from the
initial vertex for every one of these $b-1$ maximal $x_b$-paths,
and we may also assume that the graph has a maximal
$x_b$-path with mixed missing labels.  Note that the operation
given in Figure \ref{fig8} results in a good extension,
 and does not change the number of maximal
$x_b$-paths.  Together with the notes given in the preceding
paragraph, we see that for any given label $x_i$, $i=1,...,b-1$,
we may change our graph with the operation given in Figure
\ref{fig8} so that the resulting graph has a maximal $x_b$-path
with mixed missing labels and with $x_i$ as the terminal missing
label, without increasing the total number of maximal $x_b$-paths.

\begin{figure}[!ht]
{\epsfxsize=5in \centerline{\epsfbox{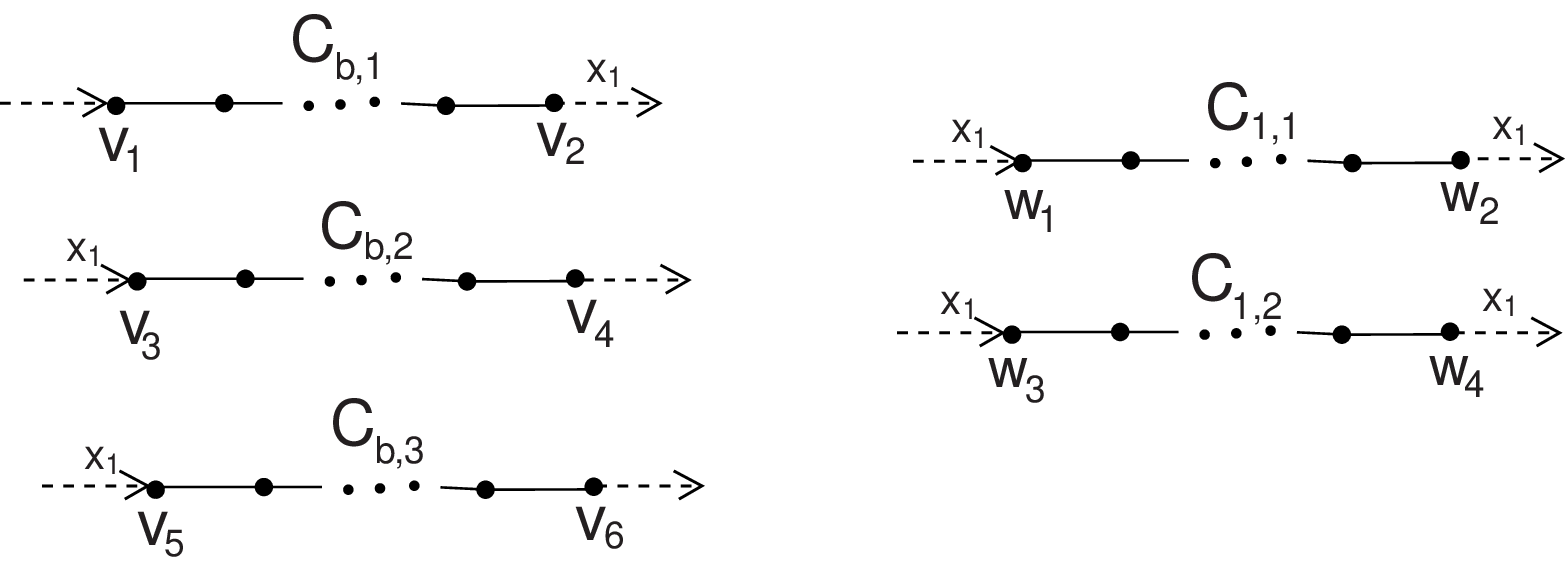}}\hspace{10mm}}
\caption{}\label{fig9}
\end{figure}

Now suppose that there is a maximal $x_b$-path $C_{b,1}$ with
mixed missing labels and  with $x_1$, say,  as terminal missing
label, and suppose that there are at least two maximal $x_1$-paths
in the graph. Then there are at least two maximal $x_b$-paths
$C_{b,2}$ and $C_{b,3}$ which have $x_1$ as their initial missing
label. The situation is illustrated in Figure \ref{fig9}. We may
assume that $v_2=w_2$. We may assume one of $v_3$ and $v_5$, say
$v_3$, is different from $w_1$ since the graph is folded. Now we
add an edge connecting $v_2$ to $v_3$ with the label $x_1$. Then
 no $x_b$-loop is created since $C_{b,1}$ has mixed
missing labels. Also no $x_1$-loop is created since $v_2\ne v_3$
and $v_3\ne w_1$. But the number of maximal  $x_b$-paths is
reduced.

\begin{figure}[!ht]
{\epsfxsize=6in \centerline{\epsfbox{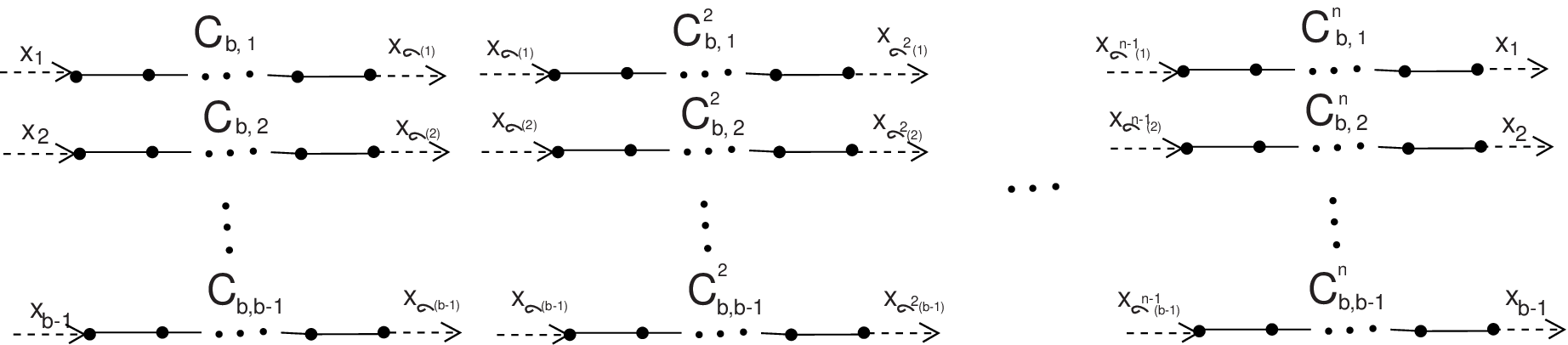}}\hspace{10mm}}
\caption{}\label{fig10}
\end{figure}

Repeating such operation, we may assume that our graph has exactly
one maximal $x_i$-path for each $i=1,2,...,b-1$. Thus there are
exactly $b-1$ maximal $x_b$-paths, $C_{b,i}$, $i=1,...,b-1$. We
may assume that the initial missing label of $C_{b,i}$ is $x_i$
for $i=1,...,b-1$. Then the terminal missing label of $C_{b,i}$ is
$x_{\s(i)}$ for some permutation $\s$ of the set
$\{1,2,...,b-1\}$. Suppose that $\s$ has order $n$. We take $n$
copies of the graph and connect them as indicated in Figure
\ref{fig10}. The resulting graph again has exactly $b-1$ maximal
$x_b$-paths, but the permutation of their missing labels is the
identity now. By adding the subgraph in Figure \ref{fig8} to the
graph at the right side, the permutation becomes a cyclic
permutation $i\ra i+1$, $i=1,...,b-1$ (defined mod $b-1$). Now we
simply fill in  $b-1$ edges at the obvious places with the right
labels. The resulting graph is what we wanted.
 The  proof of Theorem \ref{hall} is finally completed.

\begin{rem}\label{regular}
{\rm Note that the arguments in this section actually show that if
${\cal G}_\#$ is a finite, $L$-labeled, directed, folded graph with
base vertex $v_0$, with corresponding
 subgroup $G_\# = L({\cal G}_\#,v_0) \subset F=\pi_1(S^-,s)$), such that
\\
$\bullet\;$ ${\cal G}_\#$ does not contain any loop representing
the word $x_i^j$  for any  $i=1,...,b$, $j\in \z-\{0\}$, and
\\
$\bullet\;$ $y_1,...,y_r$ are some fixed, non-closed paths
 based at $v_0$ in ${\cal G}_\#$,\\
 then there is a finite, connected, $L$-regular
graph ${\cal G}_*$  such that
\\
$\bullet\;$ ${\cal G}_*$ contains ${\cal G}_\#$ as an embedded
subgraph,
 and thus in particular $y_1,...,y_a$ remain non-closed paths
 based at $v_0$ in ${\cal G}_*$, and
\\
$\bullet\;$ ${\cal G}_*$ contains no loops representing the word
$x_i^j$, for each $i=1,...,b$, $j=1,...,m_*-1$, where $m_*$ is the
number of vertices of ${\cal G}_*$. \\That is, the graph ${\cal
G}_*$ is a perfect extension of ${\cal G}_\#$.
 \\
 In terms of groups, $L({\cal G}_*,v_0)$ represents a
subgroup $G_*$ of $F$ of index $m_*$ such that
\\
$\bullet\;$ $G_*$ contains $G_\#$ as a subgroup;
\\
$\bullet\;$ $ G_*$ does not contain any of the elements $x_i^j$,
$i=1,...,b$, $j=1,...,m_*-1$;
\\
$\bullet\;$ $ G_*$ does not contain any of the elements
$y_1,...,y_r$ (considered as words in the generators in $L\cup
L^{-1}$).}
\end{rem}

This remark will be used later in Section \ref{wraps}.

\section{Lifting immersions to embeddings}\label{lifting}

We recall some of the notations from earlier sections.  Each of $i,
i_* $ denotes a number $1$ or $2$ such that $\{i,i_*\}=\{1,2\}$.
We have the universal covering maps $p:\mathbb H^3\ra M$, $p|:
\mathbb H^3\setminus {\cal B}\ra M^-$, $p_i: X_i=\tilde S_i\times
I\ra Y_i=S_i\times I$, and $p_i|:X_i^-=\tilde S_i^-\times I\ra
Y_i^-=S_i^-\times I$. We have a local isometry $f_i:Y_i=S_i\times
I\ra M$ such that $f_i|: (Y_i^-=S_i^-\times I, \p_p Y_i^-= \p
S_i^-\times I)\ra (M, \p M)$ is a proper map and such that $f_i:\p
S_i^-\ra \p M$ is an embedding.
 The surface $S_i^-$ has $n_i$
boundary components $\{\b_{i,j}, j=1,...,n_i\}$ with induced
orientation. There are $d=\D n_1n_2$ intersection points $\{t_j,
j=1,...,d\}$ between $f_1(\p S_1^-)$ and $f_2(\p S_2^-)$ in $\p
M$. The points in $f_i^{-1} \{t_1, ..., t_d\}$
 are $\{t_{i,j,k}, j=1,...,n_i, k = 1, ..., d_i\}$, where $d_i=\D n_{i_*}$,
 indexed so that  $\{t_{i,j,k}, k = 1, ..., d_i\}$ are contained
successively in the component $\b_{i,j}$ (following the orientation
of $\b_{i,j}$) for each $j=1,...,n_i$.
 We constructed a metrically complete, convex,
 hyperbolic $3$-manifold $J_i$ with a local isometry
$g_i: J_i\ra Y_i$ such that $g_i|: (J_i^-, \p_p J_i^-)\ra (Y_i^-,
\p_pY_i^-)$ is a proper map.
 The parabolic boundary of $J_i^-$,
$\p_p J_i^-$, has exactly $d$ components
 $\{D_{i,j,k}, j=1,...,n_i, k = 1, ..., d_i\}$,
 and the topological center point of $D_{i,j,k}$ is denoted $b_{i,j,k}$.
 We have
 $g_i(b_{i,j,k})=t_{i,j,k}$.
For each sufficiently large integer $n>0$, we constructed  a
compact, convex, hyperbolic $3$-manifold $C_n(J_i^-)$ which
contains $J_i^-$ as an embedded  submanifold, and a
local isometry $g_i:C_n(J_i^-)\ra Y_i$, extending the
 map $g_i:J_i^- \ra Y_i$.

\noindent{\bf Base point convention}. From now on, we will fix
$t_1$ as a basepoint for each
 of $M$, $M^-$ and $\p M$, and $t_{i,1,1}$ will
 be the base point for each of $S_i$, $S_i^-$, $Y_i$ and $Y_i^-$.
 After re-ordering $\{t_1,...,t_d\}$, we may assume that
 $f_i(t_{i,1,1}) = t_1$ for $i =1, 2$, and that
 the point $\tilde{b}=\tilde S_1\cap\tilde S_2\cap
\p B_\infty$ is in $p^{-1}(t_1)$. The point $\tilde{b}$ will be the base point
for each of $\mathbb H^3$, $X_i$, $X_i^-$, $\tilde S_i$, $\tilde
S_i^-$, and $\p B_{\infty}$. The  point $b_{i,1,1}$  will be the
base point for each of $J_i$, $J_i^-$, $\hat J_i$ and
$C_n(J_i^-)$. Under these choices of base points, we can identify,
 as in Section \ref{property}, each of
$\pi_1(M,t_1)$ and $\pi_1(M^-,t_1)$ with the group $\G$; identify
each of $\pi_1(S_i, t_{i,1,1})$, $\pi_1(S_i^-, t_{i,1,1})$,
$\pi_1(Y_i, t_{i,1,1})$, $\pi_1(Y_i^-, t_{i,1,1})$ with the
quasi-Fuchsian group $\G_i \subset \G$; and identify $\pi_1(\p M, t_1)$  with
the stabilizer of $\infty$ in $\G$. Under such identifications,
the induced map $f_i^*:\pi_1(S_i,t_{i,1,1})=\G_i \ra \pi_1(M,
t_1)=\G$ is the inclusion homomorphism, and  each of the inclusion
maps $(S_i, t_{i,1,1})\subset  (Y_i, t_{i,1,1})$, $(S_i^-,
t_{i,1,1})\subset (Y_i^-, t_{i,1,1})$, $(S_i^-, t_{i,1,1})\subset
(S_i, t_{i,1,1})$, $(Y_i^-, t_{i,1,1})\subset (Y_i, t_{i,1,1})$
and $(M^-, t_1)\subset (M, t_1)$ induces the identity isomorphism
on the fundamental groups.

\begin{figure}[!ht]
{\epsfxsize=3in \centerline{\epsfbox{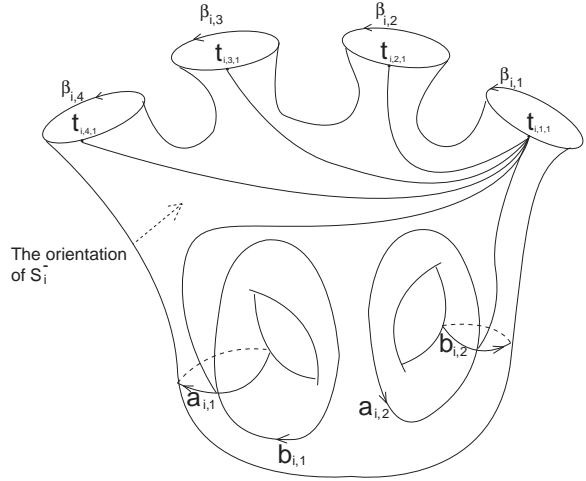}}\hspace{10mm}}
\caption{Choice of generators for $\pi_1(S_i^-,
t_{i,1,1})$}\label{generators}
\end{figure}

\noindent {\bf Choice of a free basis for $\G_i$}. Recall that
$n_i$ is the number of boundary components of  the truncated
surface $S_i^-$. Let $g_i$ be the genus of $S_i^-$.  As in Section
\ref{csfg}, the group $\G_i=\pi_1(S_i^-, t_{i,1,1})=\pi_1(Y_i,
t_{i,1,1})=\pi_1(S_i, t_{i,1,1})=\pi_1(Y_i, t_{i,1,1})$ has a set
of generators
$$X=\{a_{i,1},b_{i,1},...,a_{i,g_i},b_{i,g_i},x_{i,1},...,x_{i,n_i-1}\}$$
such that the elements
$$x_{i,1},x_{i,2},...,x_{i,n_i-1},x_{i,n_i}=[a_{i,1},b_{i,1}][a_{i,2},b_{i,2}]\cdots
[a_{i,g_i},b_{i,g_i}]x_{i,1}x_{i,2}\cdots x_{i,n_i-1}$$
 have
representative loops  freely homotopic to  the $n_i$ boundary
components of $S_i^-$ respectively. In the current case,
 we pick representative loops based at the point
$t_{i,1,1}$ for the elements
$a_{i,1},b_{i,1},...,a_{i,g_i},b_{i,g_i},x_{i,1},...,x_{i,n_{i-1}}$ as
shown in Figure \ref{generators}. For instance $x_{i,2}$ is represented by the
loop which goes along the given arc
 from $t_{i,1,1}$ to $t_{i,2,1}$, then goes around $\b_{i,2}$
 once following  the given orientation and then
 comes back to $t_{i,1,1}$ along the given arc from $t_{i,2,1}$ to $t_{i,1,1}$.
 The representative for $x_{i,n_i}$ is obtained similarly,
 except that in this case, we choose a loop which disagrees with
 the orientation of $\b_{i,n_i}$.
  Then it is easy to see that
$$x_{i,n_i}=[a_{i,1},b_{i,1}][a_{i,2},b_{i,2}]\cdots
[a_{i,g_i},b_{i,g_i}]x_{i,1}x_{i,2}\cdots x_{i,n_i-1}$$ is
in fact satisfied.

\noindent {\bf Choice of generators for
$\pi_1(C_n(J_i^-),b_{i,1,1})$}. Recall from the construction of
$C_n(J_i^-)$ that $C_n(J_i^-)$ is obtained  by gluing together
$J_i^-$ and $n_i$ multi-$1$-handles $H_{i,1}^n,...,H^n_{i,n_i}$
along the parabolic regions $D_{i,j,k}$,
$j=1,...,n_i,k=1,...,d_i$, of $J_i^-$. Recall that $b_{i,j,k}$ is
the center point of $D_{i,j,k}$
 which maps to $t_{i,j,k}$.

 Let $\alpha_{i,j,k}\subset J_i^-$ be a fixed,
 oriented path from $b_{i,1,1}$ to $b_{i,j,k}$,
 $j=1,...,n_i,k=1,...,d_i$ ($\a_{i,1,1}$ is the constant path).
 For $j=1,..., n_i$, $1 \leq k \leq d_i-1$,
 let $\delta_{i,j,k}(n)$ be the oriented geodesic arc in $H_{i,j}^n$ from
 $b_{i,j,k}$ to $b_{i,j,k+1}$.
 For $1 \leq k \leq d_i-1$,
 let $z_{i,j,k}(n)$ be the loop
 $\alpha_{i,j,k}\cdot\delta_{i,j,k}\cdot\overline{\alpha_{i,j,k+1}} $,
 where the symbol ``$\cdot$'' denotes path concatenation
 (sometimes omitted), and $\overline{\alpha_{i,j,k}}$
 denotes the reverse of $\alpha_{i,j,k}$.
 Also we always write path (in particular loop) concatenation from left to right.
 We also consider $z_{i,j,k}(n)$ as an element of $\pi_1( C_n (J_i^-),
 b_{i,1,1})$.
 Fix a set of generators
$w_{i,1},...,w_{i,\ell_i}$ for $\pi_1(J_i^-, b_{i,1,1})$. Then it's
not hard to see, by recalling the structure of $H_{i,j}^n$, that
$\pi_1(C_n(J_i^-),b_{i,1,1})$ is generated by the set of elements
$$w_{i,1},...,w_{i,\ell_i},
 z_{i,j,k}(n), 1 \leq j \leq n_i, 1 \leq k \leq d_i-1 .$$
 In fact
$$\pi_1(C_n(J_i^-),b_{i,1,1}) =\pi_1(J_i^-,b_{i,1,1})* <
z_{i,j,k}(n)| 1 \leq j \leq n_i, 1 \leq k \leq d_i-1>,$$ where $*$
denotes the free product, and $ < z_{i,j,k}(n)| 1 \leq j \leq n_i,
1 \leq k \leq d_i-1>$ is the free group freely generated by the
$z_{i,j,k}(n)$'s.

By Lemma \ref{inj}, the local isometry $g_i: (C_n(J_i^-),
b_{i,1,1})\ra (Y_i,t_{i,1,1})$ induces an injective homomorphism
$g_i^*:\pi_1( C_n(J_i^-), b_{i,1,1})\ra
\G_i=\pi_1(Y_i,t_{i,1,1})$. If $\a$ is an oriented arc in
$C_n(J_i^-)$, we use $\a^*$ to denote the oriented  arc
$g_i\circ\a$ in $Y_i$. We use $\g^*$ to denote the image of an
element $\g$ of $\pi_1( C_n(J_i^-), b_{i,1,1})$ under the map
$g_i^*$. Then $g_i^*(\pi_1( C_n(J_i^-), b_{i,1,1}))$ is
generated by the set of elements  $$w_{i,1}^*,...,w_{i,\ell_i}^*,
 z_{i,j,k}(n)^*, 1 \leq j \leq n_i, 1 \leq k \leq d_i-1.$$

Now consider the images of these generators
 in $Y_i$. The oriented
path $\alpha_{i,j,k}^*$ in $Y_i^-$ runs from $t_{i,1,1}$
to $t_{i,j,k}$.
 For $j=1,...,n_i$, $1 \leq k \leq d_i-1$,
 let $\eta_{i,j,k}$ be an oriented arc in $\beta_{i,j}$ from
 $t_{i,j,k}$ to $t_{i,j,k+1}$ following the orientation of $\b_{i,j}$, and
 let $\s_{i,j,k}\subset Y_i^-$ be the loop
 $\alpha_{i,j,k}^*\cdot\eta_{i,j,k}\cdot\overline{\alpha_{i,j,k+1}}^*$.
  Let $\s_{i,j,0}$ be the constant path based at $t_{i,1,1}$.
 Let $x_{i,j}^{\prime}$ be the loop
 $\alpha_{i,j,1}^*\cdot\beta_{i,j}\cdot\overline{\alpha_{i,j,1}}^*$,
 where  $\b_{i,j}$ is considered an oriented loop
 starting and ending at the point $t_{i,j,1}$.

\begin{lem}
\label{lk} Considered as an element in $\G_i$,
  $$z_{i,j,k}(n)^* =(\overline{\s_{i,j,k-1}}\cdots\overline{\s_{i,j,0}})(x_{i,j}')^n
 (\s_{i,j,0}\cdots\s_{i,j,k}),
$$ for each $j=1,...,n_i, k=1,...,d_i-1$.
\end{lem}

\pf From the construction of $H_{i,j}^n$, we see that
 the arc $\d_{i,j,k}(n)^*$
is isotopic in $Y_i$, with the endpoints fixed, to an
arc which: starts from the point $t_{i,j, k}$, goes around
$\b_{i,j}$ exactly $n$ times (following the orientation of
$\b_{i,j}$), then continues along
$\b_{i,j}$ until it reaches the point $t_{i,j,k+1}$. Now it is
easy to check that the loop  $z_{i,j,k}(n)^*$ is homotopic in $Y_i$,
fixing the base point $t_{i,1,1}$, to the loop
$(\overline{\s_{i,j,k-1}}\cdots\overline{\s_{i,j,0}})(x_{i,j}')^n
 (\s_{i,j,0}\cdots\s_{i,j,k})$.
This proves the lemma. \qed

\begin{rem}{\rm
The elements $\s_{i,j,k}$ $j=1,...,n_i, k=1,...,d_i-1$, are
independent of the integer $n$.}
\end{rem}

\noindent{\bf Definition.} Suppose that
 $\breve p_i:\breve{\b}_{i,j}\ra
\b_{i,j}$ is a covering
 map, and let $\breve \b_{i,j}$ have orentation induced from $\b_{i,j}$.
 Let $\a\subset \tilde \b_{i,j}$ be an embedded, connected, compact arc
 with orientation induced from $\breve \b_{i,j}$, whose
 initial point is in $\breve p_i^{-1}(t_{i,j,k})$
 and whose terminal point is in
 $\breve p_i^{-1}(t_{i,j,k+1})$ (here $k+1$ is defined mod $d_i$).
 We say that $\a$ has {\it wrapping number} $n$ if
 there are exactly
$n$ distinct points of $\breve p_i^{-1}(t_{i,j,k})$ which are
contained in the interior of $\a$.

 In   the next section we show the following

\begin{pro}\label{adjust}For each $i=1,2$ and $n\geq 0$,
there is a finite cover $\breve Y_i=\breve S_i\times I$ of
$Y_i=S_i\times I$, having the following properties:
\newline
$(1)$ $\p_p \breve Y_i^-=\p \breve S_i^-\times I$ has the same
number of components as  $\p_pY_i^-=\p S_i^-\times I$;
\newline $(2)$ the map $g_i: J_i^-\ra Y_i^-$ lifts to an embedding
$\breve g_i: J_i^-\ra \breve Y_i^-$;
\newline
$(3)$ the points
$\breve g_i (b_{i,1,1}),...,\breve g_i(b_{i,n_i,d_i})$ are evenly
spaced; i.e. there is an integer $N_i>n$ such that each  of the $n_id_i$
components of $\p \breve S_i^-\setminus\{\breve g_i
(b_{i,1,1}),...,\breve g_i(b_{i,n_i,d_i})\}$ has wrapping
number equal to the integer $N_i$.
\end{pro}

\section{Adjusting the wrapping numbers}\label{wraps}

In this section we prove Proposition \ref{adjust}.

Recall from Section \ref{ji} that $\hat J_i$ is a
connected, compact, convex, hyperbolic 3-manifold
obtained from $J_i^-$ by capping off each component of $\p_p
J_i^-$ with a compact, convex $3$-ball, and that
$\pi_1(J_i, b_{i,1,1})=\pi_1(J_i^-,b_{i,1,1})=\pi_1(\hat
J_i,b_{i,1,1})$.
 Also, $\hat J_i$ is a
submanifold of $C_n(J_i^-)$, so by Lemma \ref{inj}, $\pi_1(\hat
J_i,b_{i,1,1})$ can be considered as a subgroup of
$\pi_1(C_n(J_i^-),b_{i,1,1})$.

By Proposition \ref{sep1} there is a set of elements $y_{i,1},...,
y_{i,r_i}$ in $\G_i-g_i^*(\pi_1(\hat J_i,b_{i,1,1}))$ such that,
if $G_i$ is a finite index subgroup of $\G_i$ which separates
$g_i^*(\pi_1(\hat J_i,b_{i,1,1}))$ from $y_{i,1},...,y_{i,r_i}$,
then the local isometry $g_i:\hat J_i\ra Y_i$ lifts to an
embedding $\breve g_i$ in the finite cover $\breve Y_i$
corresponding to $G_i$ .

To prove Proposition \ref{adjust},  we shall construct a finite
index subgroup $G_i$ of $\G_i$, of sufficiently large  index $m_i$,
such that
\newline (i) $m_i=N_id_i+1$ for some integer $N_i>n$;\newline
(ii) $G_i$ contains the elements $w_{i,1}^*,...,w_{i,\ell_i}^*$ (defined
 in Section \ref{lifting}), and
thus contains the subgroup $g_i^*(\pi_1(\hat J_i,b_{i,1,1}))$;
\newline
(iii) $G_i$ contains the elements $z_{i,j,k}(N_i)^*$,
$j=1,...,n_i$, $k=1,...,d_i-1$;
\newline (iv) $G_i$ does not contain any of
$x_{i,j}^l$, $j=1,...,n_i$, and $l=1,...,m_i-1$;
\newline
(v) $G_i$ does not contain any of $y_{i,1},...,y_{i,r_i}$.

\begin{pro}\label{same}
Assuming such a subgroup $G_i$ can be found, then the corresponding
finite cover $\breve Y_i=X_i/G_i$ of $Y_i$  will satisfy all the
properties  given in Proposition \ref{adjust}.
\end{pro}

\pf  Let $\breve q_i:X_i\ra \breve Y_i$ and $\breve p_i:\breve
Y_i\ra Y_i$ be the covering maps. Properties (ii) and (v) imply
that the map $g_i:(\hat J_i, b_{i,1,1})\ra (Y_i, t_{i,1,1})$ lifts
to an embedding $\breve g_i:(\hat J_i, b_{i,1,1})\ra (\breve Y_i,
\breve q_i(\tilde{b}))$. By our  choice for the cusp $C$ of $M$
(which determines the cusp region ${\cal C}_i$ for $Y_i$), the
restriction of $\breve g_i$ on $J_i^-$ gives a proper embedding
$\breve g_i:(J_i^-, \p J_i^-)\ra (\breve Y_i^-, \p \breve Y_i^-)$,
i.e. we have (2) of Proposition \ref{adjust}.

We claim that condition (iv) implies $\p_p\breve Y_i^-$ has
 the same number of components as $\p_pY_i^-$, i.e. we
have (1) of Proposition \ref{adjust}. To see this,  recall that
each of $x_{i,j}$, $j=1,...,n_i$, has a representative loop
(Figure \ref{generators}) which is homotopic, with the base point
$t_{i,1,1}$ fixed, to an embedded loop, $x_{i,j}''$, in $S_i^-$,
 such that $x_{i,j}''$ is parallel to
$\b_{i,j}$  in $S_i^-$. Since $G_i$ does not contain any of the
elements $x_{i, j}^l$, $j=1,...,n_i$, $l=1,..., m_i-1$, then $\breve
p_i^{-1}(x_{i,j}'')$ is a single, embedded loop in $\breve S_i^-$
for each fixed $j=1,...,n_i$. Hence $\breve p_i^{-1}(\b_{i,j})$ is
a single component of  $\p\breve S_i^-$ for each $j=1,...,n_i$.
This proves the claim.

We now show that condition (i) and (iii) imply (3) of Proposition
\ref{adjust}. Namely we want to show that  the set of points
$\breve g_i (b_{i,1,1}),...,\breve g_i(b_{i,n_i,d_i})$ are evenly
spaced in $\p\breve S_i^-$ so that each  of the $n_id_i$
components of $\p \breve S_i^-\setminus\{\breve g_i
(b_{i,1,1}),...,\breve g_i(b_{i,n_i,d_i})\}$ has wrapping
number equal to $N_i$.

Consider the manifold $C_{N_i}(J_i^-)$. As noted in the previous
section, the subgroup  $$g_{i,N_i}^*(\pi_1(C_{N_i}(J_i^-),
b_{i,1,1}))\subset \G_i$$ is generated by the elements
$w_{i,1}^*,...,w_{i,\ell_i}^*$ and $z_{i,j,k}(N_i)^*$, $j=1,...,n_i,
k=1,...,d_i-1$.
 Hence the
group $g_{i,N_i}^*(\pi_1(C_{N_i}(J_i^-), b_{i,1,1}))$ is contained
in $G_i$ by conditions (ii) and (iii).
 Therefore the map
 $g_i:(C_{N_i}(J_i^-),b_{i,1,1})\ra (Y_i, t_{i,1,1})$ lifts
 to a map
 $\breve g_i:(C_{N_i}(J_i^-), b_{i,1,1})\ra
 (\breve Y_i, \breve q_i(\tilde{b}))$,
 i.e. $\breve p_i\circ\breve  g_i=g_i$.

 Let $\breve\b_{i,j}$ be the component of
 $\p \breve S_i^-$ which covers $\b_{i,j}$, $j=1,...,n_i$.
 Then by (1) and condition (i),  $\breve p_i: \breve \b_{i,j}\ra \b_{i,j}$
 is an $N_id_i+1$-fold cyclic covering, for each $j=1,...,n_i$.
 For each fixed $j=1,...,n_i$, the set of points
 $\{\breve g_{i,N_i}(b_{i,j,k}), k=1, ..., d_i\}$ divides
 $\breve \b_{i,j}$ into $d_i$ segments.

 Recall the notations established in Section \ref{lifting}.
 Consider the multi-handle $H_{i,j}^{N_i} \subset C_{N_i}(J_i^-)$
 containing the points $b_{i,j,1}, ... b_{i,j,d_i}$, and the
 geodesic arcs $\d_{i,j,k}(N_i)\subset H_{i,j}^{N_i}$, $k=1,...,d_i-1$.
 By our construction the immersed arc $g_{i,N_i}:\d_{i,j,k}(N_i)\ra S_i$
 is homotopic, with end points fixed, to the arc
 in $\b_{i,j}$ which starts at the point $t_{i,j,k}$,
 wraps $N_i$ times around $\b_{i,j}$
  and then continues to the point
 $t_{i,j,k+1}$, following the orientation of $\b_{i,j}$.
 This latter (immersed) arc lifts to an embedded arc in
$\breve \b_{i,j}$ connecting $\breve g_{i, N_i}(b_{i,j,k})$ and
$\breve g_{i, N_i}(b_{i,j,k+1})$, since  $\breve \b_{i,j}$ is an
$N_id_i+1$-fold cyclic cover of $\b_{i,j}$.
 Now it is easy to see that the conclusion of
(3) follows.
 \qed

To find the required subgroup $G_i$ of $\G_i$, we apply again the
graph technique used in Section \ref{csfg}. We shall use
terminologies established there without recalling them again. From
now on all elements of $\G_i$ will be considered as words in
letters from $L_i\cup L_i^{-1}$, where $$L_i=
\{a_{i,1},b_{i,1},...,a_{i,g_i},b_{i,g_i},x_{i,1},...,x_{i,n_i-1}\}$$
is the free basis of $\G_i$ given in Section \ref{lifting}. For
simplicity a word $w$ in letters of $L_i\cup L_i^{-1}$ shall also
be considered as a path in a $L_i$-labeled directed graph, and the
context will make it clear which is meant.

From Section \ref{csfg} we know that to find the required subgroup
$G_i$ of $\G_i$, it suffices to find a finite, connected,
$L_i$-labeled, directed graph ${\cal G}_{i}$ (with a fixed base
vertex $v_{i,0}$) with the following properties:
\newline
  (0)  ${\cal G}_{i}$ is  $L_i$-regular;\newline
(1) $m_{i}=N_id_i+1$ for some integer $N_i>n$,
 where $m_{i}$ is the number of vertices of ${\cal G}_{i}$; \newline
 (2) each of the words $w_{i,1}^*,..., w_{i,\ell_i}^*$ is  representable
 by a loop, based at $v_{i,0}$, in ${\cal G}_{i}$; \newline
(3) ${\cal G}_{i}$  contains a closed loop, based at $v_{i,0}$,
representing the word $z_{i,j,k}(N_i)^*$, for each $j=1,...,n_i,
k=1,...,d_i-1$;\newline
 (4)  ${\cal G}_{i}$ contains no
 loop representing the word $x_{i,j}^l$ for any $j=1,...,n_i$ and
 $l=1,...,m_{i}-1$;\newline
  (5) each of the words $y_{i,1},...,y_{i,r_i}$ is  representable by a
non-closed path, based at  $v_{i,0}$, in ${\cal G}_{i}$.

\noindent If such a graph can be found, then the
subgroup of $\G_i$ represented by $L({\cal G}_i,v_{i,0})$ will
satisfies all the requirements (i)-(v) set for $G_i$. Indeed,
Properties (0) and (1)  of ${\cal G}_{i}$ imply Property (i) of
$G_i$ (Lemma \ref{km}), and Properties (2)-(5) of ${\cal G}_{i}$
imply Properties (ii)-(v) of $G_i$ respectively. The task of the
rest of this section is  to construct such a graph ${\cal G}_i$.

If ${\cal G}$ is an $L_i$-labeled directed graph,
 then ${\cal G}^f$ will denote the folded graph resulting from folding
 ${\cal G}$ (see \cite{KM} for the folding operation).
 Note that if ${\cal G}$ is an $L_i$-labeled directed graph, and
 ${\cal G}'$ is a graph obtained from ${\cal G}$ by performing some
 folding operations on ${\cal G}$, then there is a uniquely associated
quotient map $q: {\cal G}\ra {\cal G}'$. In particular there is a
uniquely associated  quotient map from ${\cal G}$ to ${\cal G}^f$.

Let $n$ be a  large integer such that the manifold $C_n(J_i^-)$ is
convex, for each of $i=1,2$. Hence the local isometry
$g_i:C_n(J_i^-)\ra Y_i$ induces an injective homomorphism
$g_i^*:\pi_1(C_n(J_i^-), b_{i,1,1})\ra \pi_1(Y_i,t_{i,1,1})$.
Recall that the subgroup $g_i^*(\pi_1(C_n(J_i^-),
b_{i,1,1})\subset \G_i$ is generated by the elements $$w_{i,1}^*,
..., w_{i,\ell_i}^*, z_{i,j,k}(n)^*, 1 \leq j \leq n_i,
 1 \leq k \leq d_i-1.$$
 Let
${\cal G}_{i,0}(n)$ be the connected, finite, $L_i$-labeled, directed
graph which
 results from taking a disjoint union of embedded \textit{loops}-- representing
the reduced versions of the words $$w_{i,1}^*, ..., w_{i,\ell_i}^*,
z_{i,j,k}(n)^*, 1 \leq j \leq n_i,
 1 \leq k \leq d_i-1,$$ respectively--
 and non-closed embedded \textit{paths}-- representing the reduced
 versions of the
 words $$y_{i,1}, ..., y_{i,r_i}$$ respectively--
 and then identifying their  base vertices (their initial vertices)  to
 a common vertex $v_{i,0}$.
 Then obviously
$L({\cal G}_{i,0}(n), v_{i,0})=L({\cal G}_{i,0}(n)^f, v_{i,0})=
g_i^*(\pi_1(C_n(J_i^-), b_{i,1,1}))$.

 We may consider a graph ${\cal G}$  as metric space, by making each
 edge isometric to the interval $[0,1]$, and taking the induced
 path metric.  If $x \in {\cal G}$ and $s \in \mathbb{R}$,
 then $N_s(x)$ denotes the $s$-neighborhood of $x$ in $\cal G$.

\begin{lem}
\label{free} There is an integer $s>0$, independent
 of $n$, such that, when $n$ is large,
  the natural quotient map
 $f:{\cal G}_{i,0}(n) \ra {\cal G}_{i,0}(n)^f$ is
 an embedding on ${\cal G}_{i,0} \setminus N_s(v_{i,0})$,
and each of $f (y_{i,1}), ..., f (y_{i,r_i})$ is still  a
non-closed path
  in ${\cal G}_{i,0}(n)^f$.
\end{lem}

\pf We give an explicit construction of
 ${\cal G}_{i,0}(n)^f$, building it in steps.

Let ${\cal G}_{i,1}$ be the connected, finite, $L_i$-labeled,
directed graph which
 results from taking a disjoint union of embedded loops-- representing
the reduced versions of the words $w_{i,1}^*, ..., w_{i,\ell_i}^*$
  respectively-- and non-closed embedded paths--
 representing the reduced versions of the
 words $y_{i,1}, ..., y_{i,r_i}$  respectively--
 and then identifying their base points to
 a common vertex $v_{i,0}$.
 Then obviously $L({\cal G}_{i,1}, v_{i,0})$ represents the subgroup
$g_i^*(\pi_1( J_i^-, b_{i,1,1}))\subset \G_i$. Since the folding
operation does not change the group that the graph represents,
 $L({\cal G}_{i,1}^f, v_{i,0})=g_i^*(\pi_1(J_i^-,b_{i,1,1}))$.
 By assumption, none of the elements
$y_{i,1},...,y_{i,r_i}$ belong to the subgroup $g_i^*(\pi_1(J_i^-,
b_{i,1,1}))$, so $y_{i,1},...,y_{i,r_i}$   are still non-closed
paths in ${\cal G}_{i,1}^f$ based at $v_{i,0}$.

 Recall from Lemma \ref{lk}
 that
$$z_{i,j,k}(n)^*
=(\overline{\s_{i,j,k-1}}\cdots\overline{\s_{i,j,0}})(x_{i,j}')^n
 (\s_{i,j,0}\cdots\s_{i,j,k})$$
 for $k=1,...,d_i-1$,
where $\s_{i,j,k}$ and $x_{i,j}'$ were defined in Section 8. Note
that $x_{i,j}^{\prime}$ is conjugate to $x_{i,j}$ in $\G_i$.
 Let $\t_{i,j}$ be an element of $\Gamma_i$ such that
 $x_{i,j}^{\prime} = \t_{i,j} x_{i,j}\t_{i,j}^{-1}$.
 Let ${\cal G}_{i,2}$ be the connected graph which
 results from taking the disjoint union of  ${\cal G}_{i,1}^f$
 and non-closed embedded paths representing  the reduced version of
 the words  $\overline{\s_{i,j,k-1}}\cdots\overline{\s_{i,j,0}}\t_{i,j}$,
 $1 \leq j \leq n_i,
1\leq k \leq d_i$, respectively, and then identifying their base
vertices into a single base vertex which we still denote by
$v_{i,0}$. Then obviously we have $L({\cal G}_{i,2}^f, v_{i,0}) =
L({\cal G}_{i,2}, v_{i,0})
    = L({\cal G}_{i,1}^f, v_{i,0})=g_i^*(\pi_1(J_i^-,b_{i,1,1}))$.

 Let
 $v_{i,j,k}$ be the terminal vertex
 of the path $\overline{\s_{i,j,k-1}}\cdots\overline{\s_{i,j,0}}\t_{i,j}$
 in ${\cal G}_{i,2}^f$, for each $j=1,...,n_i$ and
 $k=1,...,d_i$.
For each  $j=1,...,n_i-1$ (when $n_i>1$) and
 $k=1,...,d_i$, let $q_{i,j,k}$ be the maximal
 $x_{i,j}$-path in $\widehat{{\cal G}_{i,2}^f}$ (a maximal $x_{i,j}$-path was
 defined in Section
 \ref{csfg}) which contains the vertex $v_{i,j,k}$.
 For $j=n_i$, and each
 $k=1,...,d_i$, let $q_{i,n_i,k}$  be the
 maximal $x_{i,n_i}$-path  in $\widehat{{\cal G}_{i,2}^f}$
determined by \newline (1) if there is a directed edge of
$\widehat{{\cal G}_{i,2}^f}$ with
 $v_{i,j,k}$ as its initial vertex
 and with the first letter of the word $x_{i,n_i}$ as its label, then $q_{i,n_i,k}$
 contains that edge;
\newline
(2) if the edge described in (1) does not exists, then $v_{i,j,k}$
is the terminal vertex of  $q_{i,n_i,k}$ and  the first letter of
the word $x_{i,n_i}$ is the
 terminal missing label of $q_{i,n_i,k}$.

\noindent Note that each $q_{i,j,k}$ is uniquely determined.
 Also no  $q_{i,j,k}$ can
 be an $x_{i,j}$-loop, since the group
 $L({\cal G}_{i,2}^f, v_{i,0})=g_i^*(\pi_1(J_i^-,b_{i,1,1}))$
      does not contain non-trivial peripheral elements of $\G_i$.
  Let  $v_{i,j,k}^-$ and $v_{i,j,k}^+$ be the initial
and terminal
 vertices of $q_{i,j,k}$ respectively.
Note that if $j<n_i$ and  $q_{i,j,k}$ is not a constant
 path, then $v_{i,j,k}^-$ and $v_{i,j,k}^+$ must be
distinct vertices; however
  $v_{i,n_i,k}^-$ and $v_{i,n_i,k}^+$ may possibly be the
 same vertex, even if $q_{i,j,k}$ is a non-constant path.

For each  $j=1,...,n_i$  and
 $k=1,...,d_i$, let $q_{i,j,k}^-$ be the embedded subpath of $q_{i,j,k}$ with
$v_{i,j,k}^-$ as the initial vertex and with $v_{i,j,k}$ as the
terminal vertex, and  let $q_{i,j,k}^+$ be the embedded subpath of
$q_{i,j,k}$ with $v_{i,j,k}$ as the initial vertex and with
$v_{i,j,k}^+$ as the terminal vertex.

 Note
that the set $\{Length(q_{i,j,k}): i=1,2, j=1,...,n_i,
k=1,...,d_i\}$ is independent of $n$, and thus is bounded. So we
may assume that for each $i=1,2$, $$n>10+max\{2Length(q_{i,j,k}):
i=1,2, j=1,...,n_i, k=1,...,d_i\}.$$
 Now for each  $j=1,...,n_i$  and
 $k=1,...,d_i-1$,  we make a new  non-closed embedded path $\Theta_{i,j,k}(n)$
representing the word $x_{i,j}^{n}$, and we add it to the graph ${\cal
G}_{i,2}^f$, by identifying the initial vertex of
$\Theta_{i,j,k}(n)$ with $v_{i,j,k}$  and the terminal vertex
  with $v_{i,j,k+1}$.
 In the resulting graph there are some obvious
 places one can perform the folding operation:
 for each  $j=1,...,n_i$ and $k=1,...,d_i-1$, the path
$q_{i,j,k}^+$ can be completely folded into the added new path
$\Theta_{i,j,k}(n)$, and likewise the path $q_{i,j,k+1}^-$ can be
completely folded into $\Theta_{i,j,k}(n)$. Let ${\cal
G}_{i,3}(n)$ be the resulting graph after performing these specific
 folding operations  for each  $j=1,...,n_i$ and
$k=1,...,d_i-1$.

 From the explicit construction, it is clear that ${\cal G}_{i,3}(n)$
 has the following properties:
 \newline
 (1) ${\cal G}_{i,3}(n)$ is a connected, finite, $L_i$-labeled,
 directed graph;
\newline
(2) ${\cal G}_{i,3}(n)$ contains loops, based at $v_{i,0}$,
representing the word $z_{i,j,k}(n)^*$ for each $j=1,...,n_i$,
$k=1,...,d_i-1$;
 \newline
 (3) ${\cal G}_{i,3}(n)$  contains ${\cal G}_{i,2}^f$ as an embedded
 subgraph;
\newline
(4)  ${\cal G}_{i,3}(n)$ is obtained from ${\cal
G}_{i,0}(n)$ by a sequence of folds.

It follows from Property (3) that the paths in ${\cal G}_{i,2}^f$
representing the words $y_{i,1}, ..., y_{i,r_i}$ remain each
non-closed
  in ${\cal G}_{i,3}(n)$, and it follows from
 Property  (4) that $L({\cal G}_{i,3}(n),v_{i,0})=L({\cal
G}_{i,0}(n), v_{i,0})=g_i^*(\pi_1(C_n(J_i^-), b_{i,1,1}))$. So
$\widehat{{\cal G}_{i,3}(n)}$ cannot have $x_{i,j}$-loops for any
$j$.

Now we consider the remaining folding operations on
${\cal G}_{i,3}(n)$ that  need to be done, in order to get the
folded graph ${\cal G}_{i,3}(n)^f$.

For each  $j=1,...,n_i$ and $k=1,...,d_i-1$, let
$\Theta_{i,j,k}(n)' =\Theta_{i,j,k}(n)\setminus (q_{i,j,k}^+\cup
q_{i,j,k+1}^-)$. Then by our construction each
$\Theta_{i,j,k}(n)'$ is an embedded $x_{i,j}$-path with
$v_{i,j,k}^+$ as its initial vertex and with $v_{i,j,k+1}^-$ as
the terminal vertex, and contains a subpath representing the word
$x_{i,j}^{10}$. Also all these paths
$\Theta_{i,j,k}(n)'$,$j=1,...,n_i$ and $k=1,...,d_i-1$, are
mutually disjoint in their interior, and their disjoint union is
equal to ${\cal G}_{i,3}(n)\setminus {\cal G}_{i,2}^f$.

For each fixed $j=1,...,n_1$, there is an $x_{i,j}$-path in ${\cal
G}_{i,3}(n)$ with $v_{i,j,1}^-$ as the initial vertex and with
$v_{i,j,d_i}^+$ as the terminal vertex, containing all the
vertices $v_{i,j,k}^{\pm}$, $k=1,...,d_i$, and containing all the
paths $\Theta_{i,j,k}(n)$, $k=1,...,d_i-1$. Since $\widehat{{\cal
G}_{i,3}(n)}$ has no $x_{i,j}$-loops, we see immediately that
when  $j<n_i$, all the vertices
 $v_{i,j,k}^{\pm}$, $k=1,...,d_i$, are mutually distinct.

We know that:\\ (1)  each vertex  $v_{i,j,k}^{\pm}$ is a initial
or terminal vertex of a maximal $x_{i,j}$-path in ${\cal
G}_{i,2}^f$;\\ (2) the graph ${\cal G}_{i,2}^f$ is an embedded,
folded subgraph
 of ${\cal G}_{i,3}(n)$;\\
 (3)  each $\Theta_{i,j,k}(n)'$ is an embedded
path in ${\cal G}_{i,3}(n)$;\\ (4)  for each
fixed $j<n_i$, all the vertices
 $v_{i,j,k}^{\pm}$, $k=1,...,d_i$, are mutually distinct.\\
\\
 It follows that the only remaining folds are at the vertices
 $v_{i,n_i,k}^{\pm}$, where possibly a single edge from
$\Theta_{i,n_i,k}(n)'$ may be folded to a single edge
 from $\Theta_{i,j,k_*}(n)'$, for
 some $1\leq j<n_i$ and some $1\leq k_*\leq d_i-1$. At such a
vertex there is at most one edge from $\Theta_{i,n_i,k}(n)'$ which
may be folded with one $x_{i,j}$-edge of $\Theta_{i,j,k_*}(n)'$.
 Thus ${\cal G}_{i,3}(n)^f$
 is obtained from ${\cal G}_{i,3}(n)$ by performing at most
 $2d_i$ folds  (which occur at some of the vertices
$v_{i,n_i,k}^{\pm}$, $k=1,...,d_i$),
  and every non-closed, reduced path
 in  ${\cal G}_{i,3}(n)$ which is based at $v_{i,0}$ will
 remain non-closed in ${\cal G}_{i,3}(n)^f$.
 In particular, the paths representing the words $y_{i,1}, ... y_{i,r_i}$ are
 each non-closed in
 ${\cal G}_{i,3}(n)^f = {\cal G}_{i,0}(n)^f$.

Let $f_3:{\cal G}_{i,3}(n)\ra {\cal G}_{i,3}(n)^f$ be the natual
map. Then by the construction, we see that  if $s_1$ is greater
than
 $2d_i+Diameter({\cal G}_{i,2}^f)$, then the  map
 $f_3: {\cal G}_{i,3}(n) \ra {\cal G}_{i,3}(n)^f$
 is an embedding on ${\cal G}_{i,3}(n) - N_{s_1}(v_{i,0})$.
Since ${\cal G}_{i,3}(n)$ is a partial folding of ${\cal
G}_{i,0}(n)$, there is a quotient map
 $g:{\cal G}_{i,0}(n)\ra {\cal G}_{i,3}(n)$. Letting $s$ be the diameter of
$g^{-1}(N_{s_1}(v_{i,0}))$, then $g$ is an embedding on ${\cal
G}_{i,0}(n) - N_{s}(v_{i,0})$. Since the map $f: {\cal G}_{i,0}(n)
\ra {\cal G}_{i,0}(n)^f={\cal G}_{i,3}(n)^f$ is the composition of the
maps $g$ and $f_3$,
 we see that $f$ is an embedding on ${\cal G}_{i,0} - N_{s}(v_{i,0})$.
 Obviously the number $s$ is independent of $n$.
The proof of Lemma \ref{free} is now complete.\qed

 Let $s$ be the constant integer guaranteed by Lemma \ref{free}.
 We may assume that $s$ is large enough so that
 $N_s(v_{i,0})$ in ${\cal G}_{i,0}(n)$ contains the loops
 $w_{i,1}^*,..., w_{i,\ell_i}^*$, the paths $y_1, ..., y_{r_i}$
 and the paths representing the words
  $\overline{\s_{i,j,k-1}}\cdots\overline{\s_{i,j,0}}\t_{i,j}$,
 $j=1,...,n_i$, $k=1,...,d_i-1$.
 (The choice of $s$ given in the proof of Lemma \ref{free}
 actually already satisfies this  requirement.)
 We may assume further that  $n$ is large enough
 so that
 the components of ${\cal G}_{i,0}(n)^f\setminus f(N_{v_{i,0}}(s))$
can be denoted by   $\Phi_{i,j,k}(n)$, $1 \leq j \leq n_i$, $1
\leq k \leq
 d_i-1$,
 such that  $\Phi_{i,j,k}(n)$ is an embedded subpath in
 $\Theta_{i,j,k}(n)'$ (and thus is a $x_{i,j}$-path) containing
 a sufficiently large power of $x_{i,j}$.
 This is clearly possible from the proof of Lemma \ref{free}.

 The next step is to modify the graph  ${\cal G}_{i,0}(n)^f$,
 by inserting copies of a certain graph $\Omega$,
 pictured in Figure \ref{fg11}, and then performing folding
 operations, to obtain a  graph  (the graph ${\cal
  G}_{i,4}(n)$given below)
  which contains loops, based at the base vertex $v_{i,0}$,
  representing the words  $$w_{i,1}^*, ...,
w_{i,\ell_i}^*, z_{i,j,k}(n+1)^*, 1 \leq j \leq n_i,
 1 \leq k \leq d_i-1,$$ respectively,
 and which contains
 non-closed paths, based at $v_{i,0}$, representing the
 words $$y_{i,1}, ..., y_{i,r_i}$$ respectively.
 Then from this graph we can go two steps further to
  find the required graph (The graph ${\cal G}_{i,6}(n)$ given
afterwards). The method for constructing ${\cal G}_{i,4}(n)$
  breaks into three cases, i.e.
 \newline(a) when $n_i$ is even, \newline(b) when $n_i>1$ is odd, and \newline(c) when $n_i=1$.
\newline In Figure \ref{fg11}, single edge loops at a vertex have label
one
 each from $L_i^*=\{a_{i,1}, b_{i,1},..., a_{i,g_i}, b_{i,g_i}\}$.
The edges in part (a) and (b)  connecting two adjacent vertices
are $x_{i,j}$-edges, $j=1,2,...,n_i-1$, (precisely $n_i-1$ edges).
In part (a) of the figure, an $x_{i,j}$-edge points from the left
vertex to the right vertex iff
 $j$ is odd, and in part (b) of the figure,
 an $x_{i,j}$-edge points from  left to
right iff
 $j$ is $1$ or an even number.
The edges in part (c) connecting the left two vertices and
pointing from left to right are labeled $a_{i,j}$ and $b_{i,j}$,
$j=1,2,...,g_i$, respectively, while the edges connecting the left
two vertices but pointing from right to left are labeled
$b_{i,1}$, $a_{i,j}$, $b_{i,j}$, $j=2,...,g_i$, respectively. The
right half of (c) is an identical copy of the left half of (c).

\begin{figure}[!ht]
{\epsfxsize=4.5in \centerline{\epsfbox{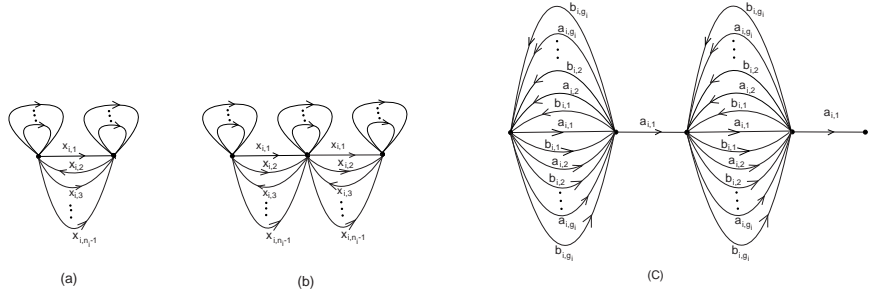}}\hspace{10mm}}
\caption{The graph $\Omega$ when
 (a) $n_i$ is even, (b) $n_i>1$ is odd, (c) $n_i=1$.}\label{fg11}
\end{figure}

 \noindent
 {\bf Case (a):}  $n_i$ is even.

 We shall insert $d_i-1$ copies of
 the graph $\Omega$ (Figure \ref{fg11} part (a)),  denoted
 $\O_k, k=1,...,d_i-1$, as follows.
 For each $1 \leq k \leq d_i-1$, we define a subset of vertices
 ${\cal{U}}_{i,k} = \{ u_{i,j,k} : 1 \leq j \leq n_i \}\subset
 {\cal G}_{i,0}(n)^f$ where,\\
 -- if $j \leq n_i-1$, then
 $u_{i,j,k}$ is a vertex in   $\Phi_{i,j,k}(n)$,
 such that there are at least three edges before it and
 after it in the directed (and thereby ordered) path $\Phi_{i,j,k}(n)$, and\\
 -- $u_{i,n_i,k}$ is the initial vertex of an edge labeled $x_1$
 in  $\Phi_{i,n_i,k}(n)$ such that
  there are at least three edges with label $x_1$ before it and
 after it  in the path $\Phi_{i,n_i,k}(n)$.

 \noindent
 Then cut ${\cal{G}}_{i,0}(n)^f$ at the vertices of ${\cal{U}}_{i,k}$,
 $k=1,...,d_i-1$,
 and for each $k$, insert the graph $\O_k$, which is a copy of
 the graph $\Omega$
   shown in Figure \ref{fg11} (a).
 That is, we\\
 (1) Form a cut graph  ${\cal{G}}_{i,0}(n)^f_c={\cal{G}}_{i,0}(n)^f
 \setminus \{U_{i,k}; k=1,...,d_i-1\}$, whose
 vertex set is obtained from
 the vertex set of ${\cal{G}}_{i,0}(n)^f$ by replacing each
 $u_{i,j,k} \in {\cal{U}}_{i,k}$ with a pair
 of vertices  $u_{i,j,k}^{\pm}$. More precisely
the point $u_{i,j,k}$ cuts the path $\Phi_{i,j,k}$ into two
components; $u_{i,j,k}^+$ is the terminal vertex of one
component, and $u_{i,j,k}^-$ is the initial vertex of the other
 component. If each pair $\{u_{i,j,k}^+,
u_{i,j,k}^-\}$ is identified into a single vertex,
 then the resulting graph is ${\cal{G}}_{i,0}(n)^f$.
\\(2) For each fixed $k=1,...,d_i-1$, we identify the
 vertex set $\{u_{i,j,k}^{\pm}, j=1,...,n_i\}$ of
 ${\cal{G}}_{i,0}(n)^f_c$
 with  the  vertices  of  $\Omega_k$ as follows:\\
-- if $j<n_i$, identify $u_{i,j,k}^+$ with the left vertex of
$\O_k$ if $j$ is odd and to the right vertex if $j$ is even, and
identify $u_{i,j,k}^-$ with the right vertex of $\O_k$ if $j$ is
odd and to the left vertex if $j$ is even,\\ -- identify
$u_{i,n_i,k}^+$ with the left vertex of $\O_k$  and identify
$u_{i,n_i,k}^-$ with the right vertex of $\O_k$.

\begin{figure}[!ht]
{\epsfxsize=4in \centerline{\epsfbox{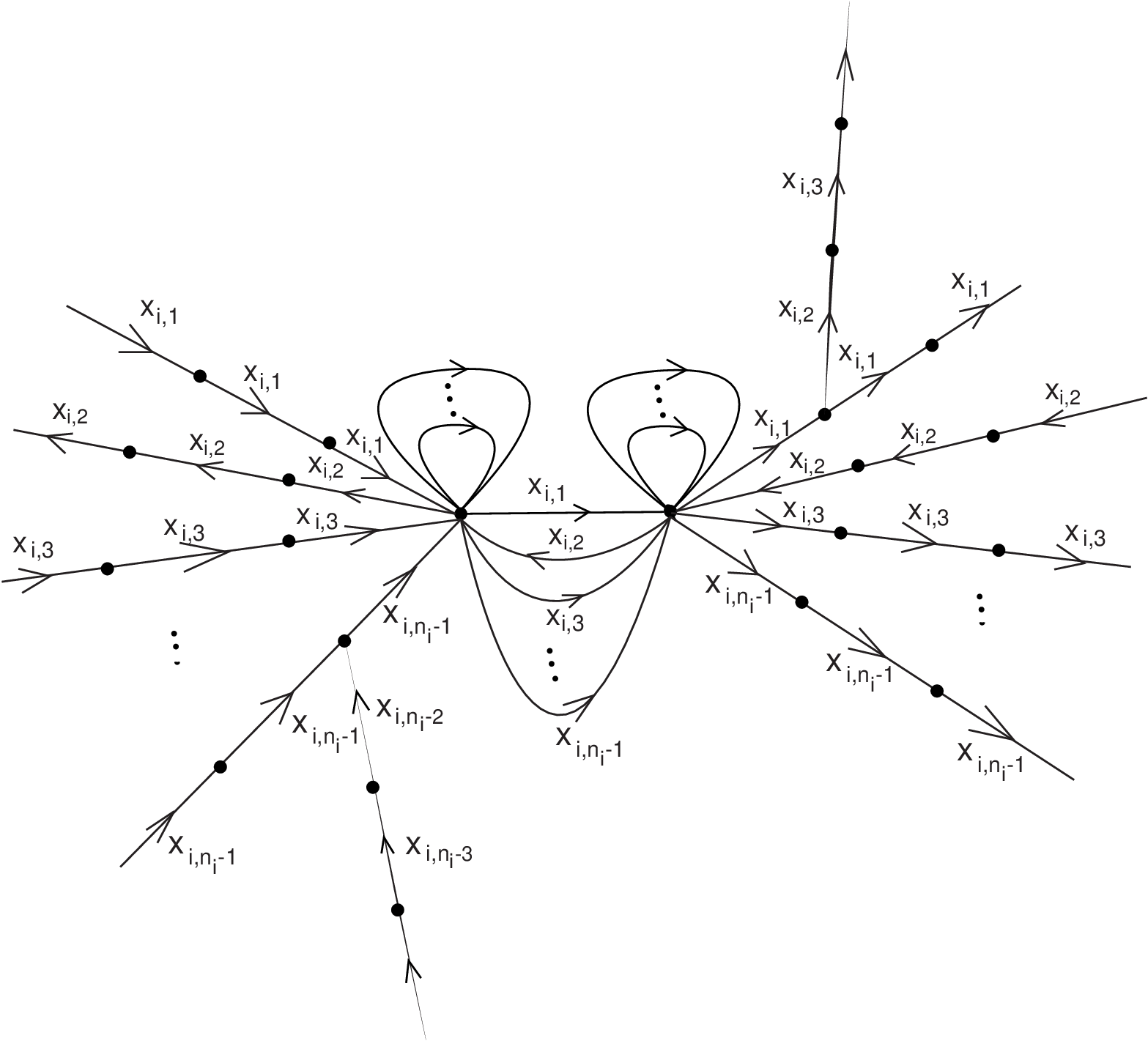}}\hspace{10mm}}
\caption{}\label{fg13}
\end{figure}

 The resulting graph is not
folded, but becomes folded graph after the following obvious
folding operation around each inserted  $\O_k$: \\ -- fold the
subpath $x_{i,n_i-1}a_{i,1}b_{i,1}a_{i,1}^{-1}b_{i,1}^{-1}\cdots
a_{i,g_i} b_{i,g_i}a_{i,g_i}^{-1}b_{i,g_i}^{-1}$  whose terminal
vertex is the vertex $u_{i,n_i,k}^+$ with the loops of $\O_k$  at
the left vertex of $\O_k$ and then with the $x_{i,n_i-1}$-edge of
${\cal G}_{i,4}(n)$ whose terminal vertex is the left vertex of
$\O_k$, and\\ -- fold the two $x_1$-edges whose initial vertices
are the right vertex of $\O_k$.\\ The resulting folded graph
${\cal G}_{i,4}(n)$ around the inserted $\O_k$ is shown in Figure
\ref{fg13}.
 By  our construction we see that ${\cal G}_{i,4}(n)$ is a folded,
 $L_i$-labeled, directed graph, with  no
 $x_{i,j}$-loops, with each of the words
 $w_{i,1}^*,...,w_{i, \ell_i}^*$ still representable
 by a loop based at $v_{i,0}$, and
 with each of the words $y_{i,1},..., y_{i,r_i}$
 still representable by a non-closed path based at $v_{i,0}$.
 Also we see that the graph ${\cal G}_{i,4}(n)$ contains
 loops based $v_{i,0}$
 representing the words $z_{i,j,k}(n+1)^*$, for
 any $j=1,...,n_i$, $k=1,...,d_i-1$.

\begin{figure}[!ht]
{\epsfxsize=4in \centerline{\epsfbox{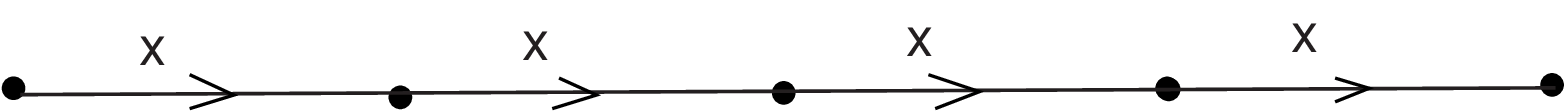}}\hspace{10mm}}
\caption{}\label{fg12}
\end{figure}

The graph  ${\cal G}_{i,4}(n)$ is not $L_i$-regular yet since it
does not contain any $x_{i,j}$-loops. So it
must contain a missing label. Let $x\in L_i$ be a missing label at
a vertex $v$ of ${\cal G}_{i,4}(n)$. Let $\a$ be a finite directed
graph consisting of a single path of edges all labeled with $x$, as
shown in Figure \ref{fg12}. We identify the left end vertex of
$\a$ to the vertex $v$ of ${\cal G}_{i,4}(n)$. The resulting
 graph ${\cal G}_{i,5}(n)$ is obviously still folded, contains
 ${\cal G}_{i,4}(n)$ as an embedded subgraph, and contains
 no $x_{i,j}$-loops for any $j=1,...,n_i$.
 By choosing a long enough path $\alpha$,
 we may assume that the number of vertices of ${\cal G}_{i,5}(n)$
 is bigger than $d_i(n+1)+1$.

Now by Remark  \ref{regular},
 we can obtain an $L_i$-regular graph ${\cal G}_{i,6}(n)$ such that
 \newline
 (1)  ${\cal G}_{i,5}(n)$ is an embedded subgraph of ${\cal G}_{i,6}(n)$;
 thus in particular  in ${\cal G}_{i,6}(n)$
 each of the words
 $w_{i,1}^*,...,w_{i, \ell_i}^*$, $z_{i,j,k}(n+1)^*$,
$j=1,...,n_i$, $k=1,...,d_i-1$ is representable
 by a loop based at $v_{i,0}$, and each of the words $y_{i,1},..., y_{i,r_i}$
 is representable by a non-closed path based at $v_{i,0}$;
 \newline
 (2) ${\cal G}_{i,6}(n)$ contains no
 loops representing the word $x_{i,j}^l$ for any $j=1,...,n_i$,
  $l=1,...,m_i^*-1$, where $m_i^*$ is the number of
 vertices of ${\cal G}_{i,6}(n)$.

  Note that $m_i^*$ is some integer larger than $d_i(n+1)+1$.
 Let $N_i=m_i^*-(d_i-1)(n+1)-1$. Then $N_i>(n+1)$.

\begin{figure}[!ht]
{\epsfxsize=3in \centerline{\epsfbox{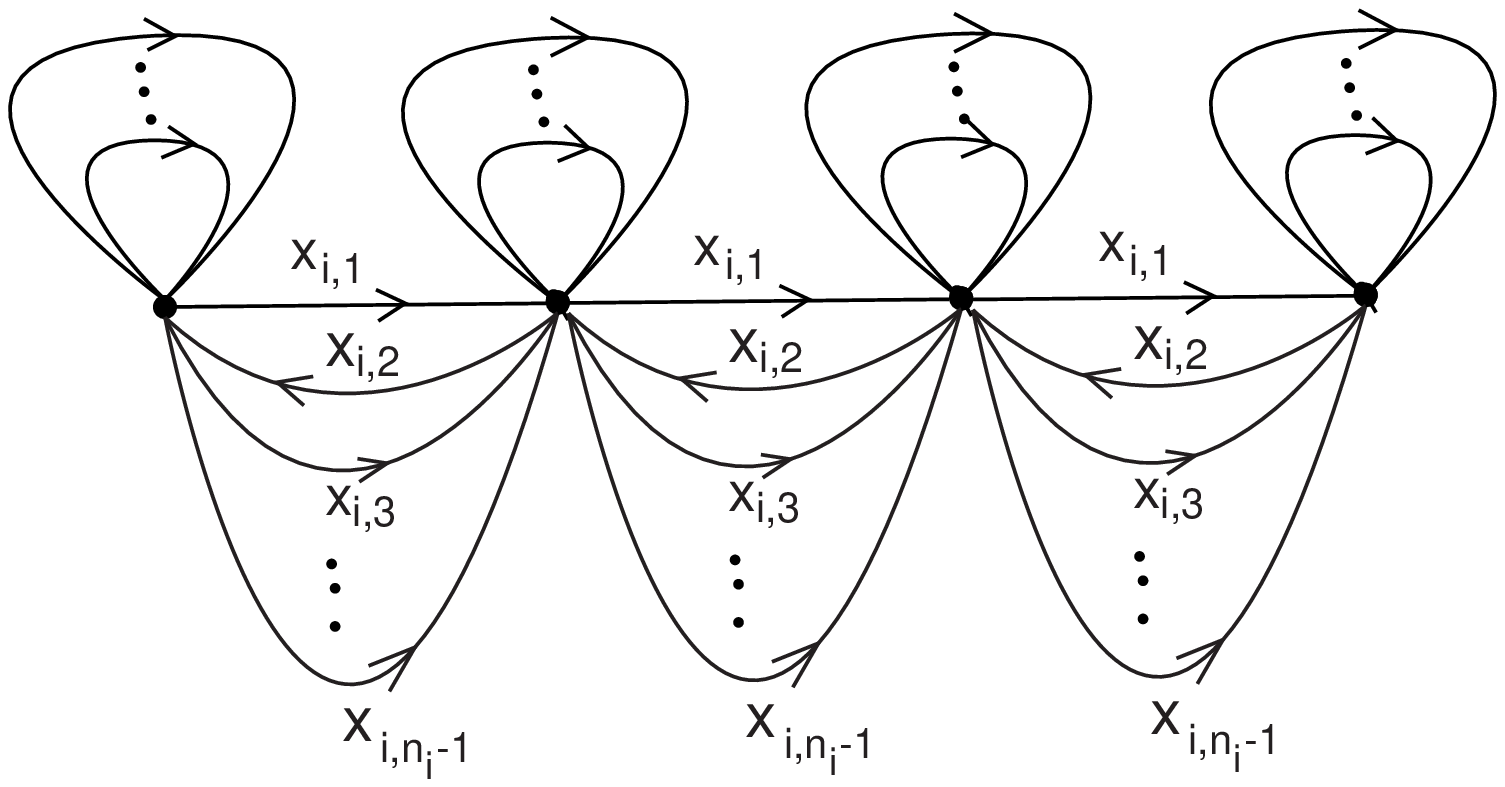}}\hspace{10mm}}
\caption{}\label{fg15}
\end{figure}

 During the transformation from ${\cal G}_{i,4}(n)$ to ${\cal G}_{i,6}(n)$,
  the subgraph of ${\cal G}_{i,4}(n)$  consisting of the
  edges which intersect the subgraph
 $\Omega_k$ (for each fixed $k=1,...,d_i-1$) remained unchanged
 since ${\cal G}_{i,4}(n)$ was locally $L_i$-regular
 already at the two vertices of
 $\O_k$.
 Now we replace $\O_k$, for each of $k=1,...d_i-1$, by
 a graph similar to $\Omega$ but with $N_i-n+1 \geq 3$
 vertices (Figure \ref{fg15}
 illustrates such a graph with four vertices).
 Then the resulting graph ${\cal G}_{i, 7}(n)$ has the following  properties.
\newline
 (1)  ${\cal G}_{i,7}(n)$ is  $L_i$-regular;\newline
 (2) each of the words $y_{i,1},...,y_{i,r_i}$ is still representable by a
non-closed path based at $v_{i,0}$ in ${\cal G}_{i,7}(n)$,\newline
 (3) each of the words $w_{i,1}^*,..., w_{i,\ell_i}^*$ is still representable
 by a loop based at $v_{i,0}$ in ${\cal G}_{i,7}(n)$, \newline
 (4)  ${\cal G}_{i,7}(n)$ contains no
 loops representing the word $x_{i,j}^l$ for each $j=1,...,n_i$ and
 each $l=1,...,m_i-1$,
 where $m_i$ is the number of vertices of ${\cal G}_{i,7}$,\newline
 (5) ${\cal G}_{i,7}(n)$  contains a closed loop based at $v_{i,0}$
representing the word $z_{i,j,k}(N_i)^*$, for each $j=1,...,n_i,
k=1,...,d_i-1$, and \newline
 (6) $m_i$, the number of vertices of
${\cal G}_{i,7}(n)$, is equal to $N_id_i+1$. \newline
 Properties (1)-(5) are obvious by the construction, while property (6)
follows by a simple calculation. Indeed
\begin{eqnarray*}
m_i &=&  m_i^* + (N_i-n+1-2)(d_i-1)\\
 &=&[N_i + (d_i-1)(n+1)+1] + (N_i-(n+1))(d_i-1)\\
 &=& N_id_i+1.
\end{eqnarray*}

\noindent {\bf Case (b)} $n_i>1$ is odd.

We modify the graph ${\cal G}_{i,0}(n)^f$
 as follows.
For each of $k=1,...,d_i-1$, we define a subset of vertices
 ${\cal{U}}_{i,k} = \{ u_{i,j,k} : 1 \leq j \leq n_i \}
 \cup \{  u_{i,n_i,k}^{\prime} \}\subset {\cal G}_{i,0}(n)^f$,
 where\\
 -- if $j \leq n_i-1$, then
 $u_{i,j,k}$ is a vertex in   $\Phi_{i,j,k}(n)$,
 such that there are at least three edges after it and
 at least three edges before it
in the directed path $\Phi_{i,j,k}(n)$; \\ -- $u_{i,n_i,k}$ is the initial
vertex of an edge labeled $x_1$
 in  $\Phi_{i,n_i,k}(n)$ such that
  there are at least three  edges with label $x_1$ before it
  in the directed path $\Phi_{i,j,k}(n)$;
  and\\
  --$u_{i,n_i,k}^{\prime}$  is the initial vertex of an edge
with label $x_2$ in  $\Phi_{i,j,k}(n)$ which appears after the
vertex $u_{i,j,k}$ in the directed path $\Phi_{i,j,k}$.  We also insist
that  $\Phi_{i,j,k}(n)$ contains
 at least three  edges with label $x_1$ between
$u_{i,n_i,k}$ and $u_{i,n_i,k}^{\prime}$
 and at least three  edges with label $x_1$ after $u_{i,n_i,k}^{\prime}$.

 Then cut ${\cal{G}}_{i,0}(n)^f$ at the vertices of ${\cal{U}}_{i,k}$,
 $k=1,...,d_i-1$,
 and for each $k$, insert the graph $\O_k$, which is a copy of the graph $\Omega$
   shown in Figure \ref{fg11} (b).
 That is, we\\
 (1) Form a cut graph  ${\cal{G}}_{i,0}(n)^f_c={\cal{G}}_{i,0}(n)^f\setminus
 \{{\cal U}_{i,k}; k=1,...,d_i-1\}$,
 defined as in Case a, with obvious modifications,
 i.e. we have similarly defined pairs of vertices
 $u_{i,j,k}^{\pm}$, $u_{i,n_i,k}^{'\pm}$ for
 ${\cal{G}}_{i,0}(n)^f_c$ such that
 if each such $\pm$ pair of vertices  are identified,
 then the resulting graph is the original ${\cal{G}}_{i,0}(n)^f$.
\\(2) For each fixed $k=1,...,d_i-1$, we identify the
 vertex set $\{u_{i,j,k}^{\pm}, u_{i,n_i,k}^{'\pm}, j=1,...,n_i\}$ of
 ${\cal{G}}_{i,0}(n)^f_c$
 with  the left and right-most vertices  of  $\Omega_k$ as follows:\\
-- if $j<n_i$, and $j=1$ or $j$ is even, then
 identify $u_{i,j,k}^+$ with the left-most vertex of
$\O_k$ and $u_{i,j,k}^-$ with the right-most vertex.\\
-- if $j<n_i$, $j \neq 1$ and $j$ is odd, then identify $u_{i,j,k}^+$
with the right-most vertex of $\O_k$ and $u_{i,j,k}^-$
 with the left-most vertex,\\
-- identify $u_{i,n_i,k}^+$ with the left-most vertex of $\O_k$
and identify $u_{i,n_i,k}^-$ with the right-most vertex of $\O_k$,
\\
-- identify $u_{i,n_i,k}^{'+}$ with the left-most vertex of $\O_k$
and identify $u_{i,n_i,k}^{'-}$ with the right-most vertex of
$\O_k$.

\begin{figure}[!ht]
{\epsfxsize=4in \centerline{\epsfbox{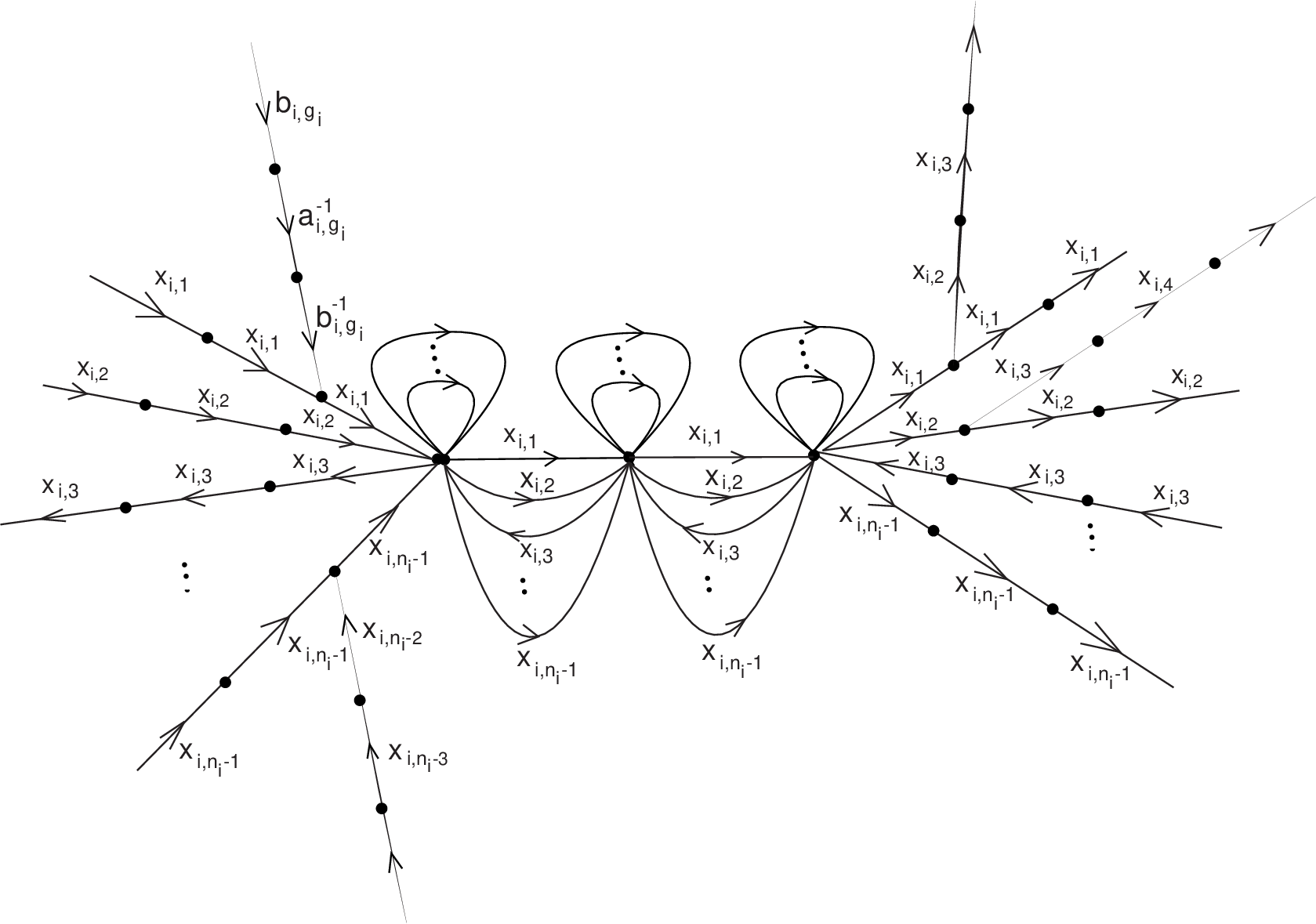}}\hspace{10mm}}
\caption{}\label{fg14}
\end{figure}

 The resulting graph is not
folded, but  becomes folded graph after the following
folding operations are performed around each inserted  $\O_k$: \\ -- fold the
path $x_{i,n_i-1}a_{i,1}b_{i,1}a_{i,1}^{-1}b_{i,1}^{-1}\cdots
a_{i,g_i} b_{i,g_i}a_{i,g_i}^{-1}b_{i,g_i}^{-1}$  whose terminal
vertex is  the vertex $u_{i,n_i,k}^+$ with the loops of $\O_k$ at
the left-most vertex of $\O_k$ and then with the
$x_{i,n_i-1}$-edge of ${\cal G}_{i,4}(n)$ whose terminal  vertex
is the left-most vertex of $\O_k$,
\\ -- fold the two $x_{i,1}$-edges whose initial vertices are the right-most
 vertex of $\O_k$,
\\ -- fold the two $x_{i,1}$-edges whose terminal vertices are the left-most
vertex of $\O_k$,
\\ -- fold the two $x_{i,2}$-edges whose initial vertices are the right-most
vertex of $\O_k$.\\ The resulting folded graph ${\cal
G}_{i,4}(n)^f_0$ around the inserted $\O_k$ is shown in Figure
\ref{fg14}.
 By  our construction we see that ${\cal G}_{i,4}(n)^f$ is a folded,
 $L_i$-labeled, directed graph, with  no
 $x_{i,j}$-loops, with each of the words
 $w_{i,1}^*,...,w_{i, \ell_i}^*$ still representable
 by a loop based at $v_{i,0}$, and
 with each of the words $y_{i,1},..., y_{i,r_i}$
 still representable by a non-closed path based at $v_{i,0}$.
 Also we see that the graph ${\cal G}_{i,4}(n)$ contains
 loops based $v_{i,0}$
 representing the words $z_{i,j,k}(n+2)^*$, for
 all $j=1,...,n_i$, $k=1,...,d_i-1$.

\begin{figure}[!ht]
{\epsfxsize=3in \centerline{\epsfbox{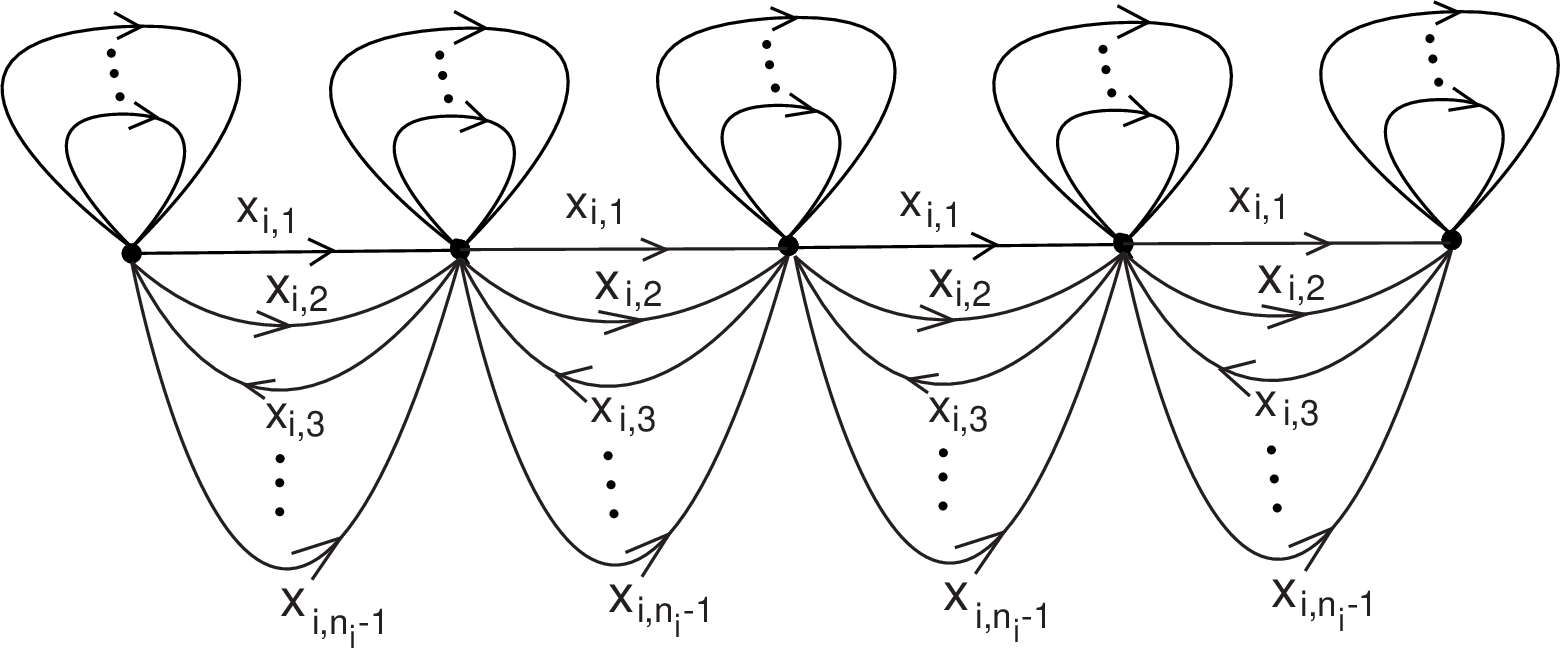}}\hspace{10mm}}
\caption{}\label{fg16}
\end{figure}

We then define ${\cal G}_{i,5}(n)$ and  ${\cal G}_{i,6}(n)$
in a similar manner as Case a; here we may assume that
 ${\cal G}_{i,5}(n)$ has at least $(d_i-1)(n+2)-1$ vertices.
 Let $m_i^*$ be the number of vertices of
${\cal G}_{i,6}$, and let $N_i= m_i^*-(d_i-1)(n+2)-1$.
  To form  ${\cal G}_{i,7}(n)$, we  replace each
 graph $\Omega_k$, $k=1,...,d_i-1$  in ${\cal G}_{i,6}(n)$
 with a graph similar to Figure \ref{fg11}(b) but with
 $1+N_i- n$ vertices.
 In the current case, we need $1+N_i-n$ to be an odd integer
 in order for the construction to work.
 (Figure \ref{fg16} illustrates such a graph with five vertices).
 This is made possible by the following

\begin{lem}
$N_i - n$ is even.
\end{lem}

\pf Since $n_i$ is odd, then the Euler characteristic
$\chi(S_i^-)$ of $S_i^-$ is odd.
 Let $\widehat{S}_i^-$ be the cover of $S_i^-$ corresponding
 to the subgroup  $L({\cal G}_{i,6}(n), v_{i,0})$ of $\G_i$.
 Due to the property (2) of the graph ${\cal G}_{i,6}$,
 $\widehat{S}_i^-$ also has $n_i$ boundary components
 (cf. the second paragraph in the proof of Proposition \ref{same}).
 So  $\chi(\widehat{S}_i^-)$ is also odd.
 Therefore the degree of the cover, which is $m_i^*$, must be odd.
 We have that $m_i^* = N_i + (d_i-1)(n+2) +1$.
 Since $n_i$ is odd, $n_{i_*}$ is even (see Section \ref{cqf}).
 Thus $d_i-1=\D n_{i_*}-1$ is odd. Thus $N_i$ and $n$ are both
 even or both odd. So $N_i - n$ is even. \qed

 The rest of the argument proceeds by obvious analogy with
 the case where $n_i$ is even.
 That is, the graph ${\cal G}_{i,7}(n)$ is  a graph with the
  properties listed as (1)-(6) in Case a.
 Indeed, Properties (1)-(5) are immediate.
 To verify Property (6), we let $m_i$ be the number of
 vertices of ${\cal G}_{i,7}$, and then we have:
\begin{eqnarray*}
m_i &=& m_i^*+ (1+N_i-n-3)(d_i-1)\\
 &=&  N_i + (d_i-1)(n+2) +1+ (N_i-n-2)(d_i-1)\\
 &=&  N_id_i+1.
\end{eqnarray*}

\noindent {\bf Case (c):}  $n_i=1$.

We modify the graph ${\cal G}_{i,0}(n)^f$
 as follows.
For each of $k=1,...,d_i-1$, we pick a pair  vertices $\{u_{i,k},
u_{i,k}'\}$ in $\Phi_{i,1,k}$ as follows:\\
  --$u_{i,k}$  is the terminal vertex of an edge
with label $a_{i,1}$ in  $\Phi_{i,j,k}(n)$ such that
 there are at least three  edges with label $a_{i,1}$ before $u_{i,k}$ in
the directed path $\Phi_{i,j,k}(n)$; and \\
 --$u_{i,k}^{\prime}$  is the terminal vertex of an edge
with label $b_{i,1}$ which appears after the
vertex $u_{i,k}$.  We also insist that
 there are at least three  edges with label $b_{i,1}$ between
$u_{i,k}$ and $u_{i,k}^{\prime}$
 and that there are at least three  edges with label $b_{i,1}$
 after $u_{i,k}^{\prime}$ in
the path $\Phi_{i,j,k}(n)$.

\noindent
 Then cut the graph
${\cal{G}}_{i,0}(n)^f$ at all the pairs
 of vertices $\{u_{i,k}, u_{i,k}'\}$, $k=1,...,d_i-1$,
and for each $k$, insert the graph $\O_k$-- which is a copy of
the graph $\Omega$
   shown in Figure \ref{fg11} (c)-- as follows.
 Form a cut graph  ${\cal{G}}_{i,0}(n)^f_c={\cal{G}}_{i,0}(n)^f\setminus
 \{u_{i,k}, u_{i,k}'; k=1,...,d_i-1\}$,
 and let
 $u_{i,k}^{\pm}$, $u_{i,k}^{'\pm}$ be the corresponding vertices
 for
 ${\cal{G}}_{i,0}(n)^f_c$.
 For each fixed $k=1,...,d_i-1$, we identify the
 vertex $u_{i,k}^+$ with the left-most vertex of $\O_k$, identify
$u_{i,k}^-$ with the right-most vertex of $\O_k$,
 identify $u_{i,k}^{'+}$ with the right-most vertex of $\O_k$
and identify $u_{i,k}^{'-}$ with the left-most  vertex of
$\O_k$.

\begin{figure}[!ht]
{\epsfxsize=5in \centerline{\epsfbox{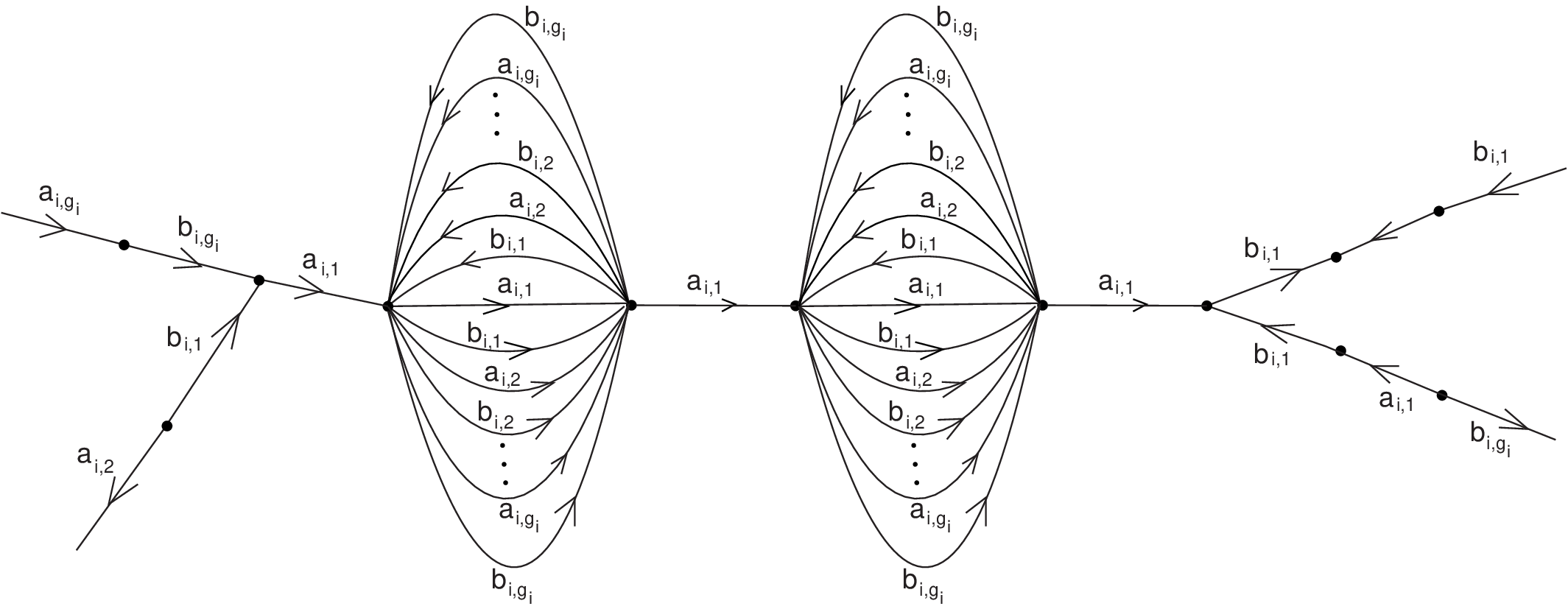}}\hspace{10mm}}
\caption{}\label{fg17}
\end{figure}

 The resulting graph is not
folded, but  becomes folded graph after a single folding operation
around each inserted  $\O_k$:
 fold the two
$a_{i,1}$-edges whose terminal  vertices are the right-most vertex
of $\O_k$. The resulting folded graph ${\cal G}_{i,4}(n)^f_0$
around the inserted $\O_k$ is shown in Figure \ref{fg17}.
 By  our construction we see that ${\cal G}_{i,4}(n)^f$ is a folded
 $L_i$-labeled directed graph, with  no
 $x_{i,1}$-loops, with each of the words
 $w_{i,1}^*,...,w_{i, \ell_i}^*$ still representable
 by a loop based at $v_{i,0}$, and
 with each of the words $y_{i,1},..., y_{i,r_i}$
 still representable by a non-closed path based at $v_{i,0}$.
 Also we see that the graph ${\cal G}_{i,4}(n)$ contains
 loops based at $v_{i,0}$
 representing the words $z_{i,1,k}(n+4)^*$, for
 all $k=1,...,d_i-1$.

As in the previous case, we get ${\cal G}_{i,5}(n)$ and
${\cal G}_{i,6}(n)$. In the current case, $N_i= m_i^*-
(d_i-1)(n+4)-1$, which is assumed larger than $n+4$ (since $m_i^*$
can be assumed arbitrary large).
  To form  ${\cal G}_{i,7}(n)$, we replace the left half
 (with three vertices) of
   $\Omega_k$, for each  $k=1,...,d_i-1$,
 with a graph similar to Figure \ref{fg11}(c) but with
 $N_i- n-1$ vertices.
 In the current case, we also need $1+N_i-n$ to be an odd integer
 in order for the construction to work.
 This is true, and  can be proved as in Case (b).
It is easy to see that ${\cal G}_{i,7}(n)$ has all the
Properties (1)-(5).  To verify Property (6), we have:
\begin{eqnarray*}
m_i &=& m_i^*+(N_i-n-1-3)(d_i-1)\\
    &=& N_i+(d_i-1)(n+4)+1+ (N_i-n-4)(d_i-1)\\
 &=& N_id_i+1
\end{eqnarray*}

\section{$HS$-manifolds}\label{hsmf}

We call a compact, connected, orientable $3$-manifold $W$ with boundary is
 an \textit{$HS$-manifold} if it has the form $W=H\cup (S\times I)$,
where
\newline
(i) each component of $H$  is a  handlebody of genus at least one;
\newline
(ii) each component of $S$ is a compact orientable surface with
boundary;
\newline
(iii)  $H\cap (S\times I)=\p S\times I$;
\newline
(iv) each component of $H\cap (S\times I)=\p S\times I$ is an
annulus in $\p H$ which is homotopically non-trivial in $H$.

\begin{lem}\label{gib}
Let $W=H\cup (S\times I)$ be an $HS$-manifold. Let $A$ denote the
set of annuli $H\cap (S\times I)=\p S\times I$. Suppose that
 $S$ has no disk components, and that for
 every compressing disk $D$ of $H$, the set
 $D \cap A$ has at least two components.
 Then $W$ has incompressible boundary.
 \end{lem}

\pf Suppose otherwise that $\p W$ is compressible in $W$. Let
$(B,\p B)\subset (W, \p W)$ be a compressing disk. Isotope $B$ so
that it intersects the set of  annuli $A$ in a collection of
properly embedded arcs and simple closed curves.  Since no
component of $S$ is a disk, and since each component of $A$ is
non-trivial in $H$, we can remove, by isotopy of $B$, all simple
closed curve components of $B\cap A$ (by a standard inner-most
argument, using also the fact that $H$ and $S\times I$ are
irreducible $3$-manifolds).

Note that the intersection $A\cap B$ cannot be empty since
otherwise $B$ would be contained in $S\times I-A$ but each
component of  $\p (S\times I)-A$ is incompressible in $S\times I$.

 We may also assume that each arc component of $A\cap B$ is
 essential in $A$. For otherwise we can surger the disk $B$
 along an outer-most such arc in $A$ to get a compressing disk of
 $W$ whose intersection with $A$ has fewer components.

Now $B\cap A$ is a set of arcs, each of which is essential in $A$. Let
$\a$ be a component of $B \cap A$ which is outer-most in $B$.
 Let $\b$ be the component of $\p
B\setminus\p \a$ whose interior is disjoint from $A$, and let
$B_1$ be the sub-disk of $B$ co-bounded by $\a$ and $\b$. Then
$B_1\cap A=\p B_1\cap A= \a$, and thus if $B_1$ is contained in
$H$, then it must be an essential compressing disk in $H$. But by
our assumption no such compressing disk exists. On the other hand,
there is no properly embedded disk in $S\times I$ which intersects
$\p S\times I$ in a single essential arc. \qed

\section{Proof of Theorem \ref{qfs}}\label{pf}

In Section \ref{wraps}, we found, for each $i=1,2$, a finite cover
$\breve  Y_i=\breve S_i\times I$ of $Y_i=S_i\times I$, such that
the map $g_i:J_i^-\ra Y_i^-$ lifts to an embedding $\breve g_i:
 J_i^-\ra \breve Y_i^-$, and the
$d$ components of  $\breve g_i(\p_p J_i^-)$ are evenly spaced in
$\p_p\breve Y_i^-$, far apart from each other in $\p_p\breve
Y_i^-$. Recall  from Section \ref{kiyi} that $K_i$ is an embedded
submanifold of $J_i$ with an $R$-collared neighborhood in $J_i$, and
that $(K_i^-, \p_p K_i^-)$ is properly embedded in the pair
$(J_i^-, \p_p J_i^-)$, with a relative $R$-collared neighborhood. It
follows that the pair $(\breve g_i(K_i^-), \breve g_i(\p_p
K_i^-))$ has a relative $R$-collared neighborhood in $(\breve Y_i^-,
\p_p\breve Y_i^-)$.

Also recall  from Section \ref{kiyi} that $K_1^-$ and $K_2^-$ are
isometric under the isometry $h:K_1\ra K_2$.
Thus  there is a corresponding isometry from  $\breve g_1(K_1^-)$
to $\breve g_2(K_2^-)$,
 which is $\breve g_2\circ h\circ\breve g_1^{-1}$.

Now let $\breve Y^-$ be the union of $\breve Y_1^-$ and $\breve
Y_2^-$ with $\breve g_1(K_1^-)$ and $\breve g_2(K_2^-)$ identified
by the isometry. Let $U_k^-\subset \breve Y^-$ be the
identification of $\breve g_1(K_{1,k}^-)$ with $\breve
g_2(K_{2,k}^-)$, $k=1,...,q$, and let $U^-$ be the disjoint union
of $U_k^-$'s. Then $\breve Y^-$ is a connected metric space,
 with a path metric induced from the metrics on $Y_1^-$ and $Y_2^-$.
 There is an induced local isometry $f:\breve{Y}^- \ra M$.

 Define the parabolic boundary, $\p_p\breve Y^-$, of $\breve
Y^-$ to be the union of $\p_p\breve Y_1^-$ and $\p_p\breve Y_2^-$,
with $\breve g_1(\p_pK_1)$ and $\breve g_2(\p_pK_2)$ identified by
the isometry $\breve g_2\circ h\circ\breve g_1^{-1}$. The
parabolic boundary of $U^-$ is defined to be the
identification of $\breve g_1(\p_p K_1)$ and $\breve g_2(\p_p K_2)$.
 Let $D_j$, $j= 1, ..., d$,
  be the components of the parabolic boundary $\p_pU^-$ of
$U^-$, and let $s_j$ be the topological center point of $D_j$
 (i.e. the $s_j$'s are the intersection points of $\p \breve
S_1^-$ and $\p \breve S_2^-$ in $\p_p \breve Y^-$).
 Since $\breve{g}_i(K_i^-)$ has an
 $R$-collared neighborhood in $\breve Y_i^-$, then $U^-$ has
an $R$-collared neighborhood in $\breve Y^-$.

Recall $f_i:(S_i^-, \p S_i^-)\ra (M^-, \p M^-)$ is a proper map,
  such that $f_i|_{\p S_i^-}:\p S_i^-\ra \p M$ is an
embedding for each $i=1,2$. Let $\b_{i,j}^*=f_i(\b_{i,j})$.  Then
$\D$ is the intersection  number between $\b_{1,1}^*$ and
$\b_{2,1}^*$, and $t_1,...,t_d$ are the $d=n_1n_2\D$ intersection
points between $\{\b_{1,j}^*,j=1,...,n_1\}$ and $\{\b_{2,j}^*,
j=1,...,n_2\}$  (since each $\b_{i,j}^*$ is a Euclidean circle in
the Euclidean torus $\p M^-$). Recall also that $\breve \b_{i,j}$,
$j=1,...,n_i$, are boundary components of $\p \breve S_i^-$, and
each $\breve \b_{i,j}$ is the cyclic covering of $\b_{i,j}$ of
order $m_i=N_id_i+1$. Recall that by our convention, $t_1$ is the
base point for each of $M$, $M^-$, $C$ and $T=\p M=\p C$, and that
$t_1$ is one of intersection points between $\b_{1,1}^*$ and
$\b_{2,1}^*$. We may consider $\b_{1,1}^*$ and $\b_{2,1}^*$ as two
elements in $\pi_1(T,t_1)=\pi_1(C, t_1)$. Now let $\cal A$ be the
subgroup of $\pi_1(T, t_1)$ generated by the two elements
$(\b_{1,1}^*)^{m_1}$ and $(\b_{2,1}^*)^{m_2}$. Then $\cal A$ is a
rank two subgroup of $\pi_1(T,t_1)=\pi_1(C, t_1)$ of
finite index. Let $p_{0}:\breve C\ra C$ be the covering
corresponding to $\cal A$. By our construction, $\p_p\breve Y^-$
can be embedded isometrically in $\breve T=\p \breve C$ such that
$p_0:\breve \b_{i,j}\ra \b_{i,j}^*$ is the map $\breve \b_{i,j}\ra
\b_{i,j}\ra\b_{i,j}^*$ for each $j=1,...,n_i$. Thus the geometric
intersection number in $\breve T$ between $\breve \b_{1,1}$ and
$\breve \b_{2,1}$ is equal to $\D$ and
 there are $d=n_1n_2\D$ intersection points
$\{s_k,k=1,...,d\}$ between $\{\breve\b_{1,j},j=1,...,n_1\}$ and
$\{\breve\b_{2,j},j=1,...,n_2\}$. We may assume that the $s_k$'s are
indexed so that   $p_0(s_k)=t_k$, $k=1,...,d$.
 From the construction of Section \ref{wraps}, the points
$\{s_k,k=1,...,d\}$ divide the circles
$\{\breve\b_{i,j},j=1,...,n_1\}$ into segments,  each of which has wrapping
number $N_i$. Thus  $\{\breve\b_{i,j},
i=1,2,j=1,...,n_i\}$ divides the torus $\breve T$ into a set of
Euclidean parallelograms with
 long sides (because the wrapping numbers $N_1$ and
$N_2$ can be chosen arbitrarily large).

 We now  replace the $R$-collared
neighborhood of $U^-$ in $\breve Y^-$ by a hyperbolic $3$-manifold
$\bar U^-$, whose construction is given below, such that
\newline (i) $\bar U^-$ is a thickening  of
$U^-$;\newline (ii) the new space $ Y^-=\breve Y_1^-\cup \bar
U^-\cup \breve Y_2^-$ is a connected, compact, hyperbolic
$3$-manifold, locally convex everywhere except on its parabolic
boundary, whose metric restricts to
the original metric on $\breve Y^-$;
\newline
(iii) $Y^-$ has a local isometry $f:Y^- \ra M$ which extends
 the local isometries
$f_i \circ \breve{p}_i: \breve{Y}_i^- \ra M$;
\newline
 (iv)  the
parabolic boundary of $Y^-$ is a regular neighborhood of that of
$\breve Y^-$ in the torus $\breve T$ and thus the complement of
$\p_p(Y^-)$ in $\breve T$ is
  a set of ``round-cornered parallelograms''
 in $\breve T$ (cf. Figure \ref{annulus});
\newline
(v) since each such parallelogram given in (iv) has  very long
sides, we can  cap off $\p_p  Y^-$ with a solid cusp $C_0$;
the resulting manifold $Y$ is a convex, hyperbolic 3-manifold with
a cusp, and $Y$ has a local isometry into $M$.

 We  now provide more details.
 First we construct  $\bar U^-$,  component-wise.
 We illustrate the construction of $\bar U_k^-$ for
the component $U_k^-$ of $U^-$.
 Recall the construction of $K_{i,k}$ given in
Section \ref{kiyi}. It is the quotient space of $Z_{i,j_k}\subset
X_i$ under the group $\G_{i,j_k}$. Recall that $\breve q_i:X_i\ra
\breve Y_i$ is the universal covering map. Thus  $Z_{i,j_k}^-$ is
the universal cover of $\breve g_i(K_{i,k}^-)$ under the map
$\breve q_i$. Also there are elements $\g_{i,j_k}\in \G$ such that
$X_{i,j_k}=\g_{i,j_k}(X_i)$, $W_{j_k}=X_{1,j_k}\cap X_{2,j_k}$,
and $Z_{i,j_k}=\g_{i,j_k}^{-1}(W_{j_k})$. The space $W_{j_k}$ is invariant
under the action of the group $\g_{1,j_k}\G_1\g_{1,j_k}^{-1}\cap
\g_{2,j_k}\G_2\g_{2,j_k}^{-1}$.

\begin{figure}[!ht]
{\epsfxsize=3in \centerline{\epsfbox{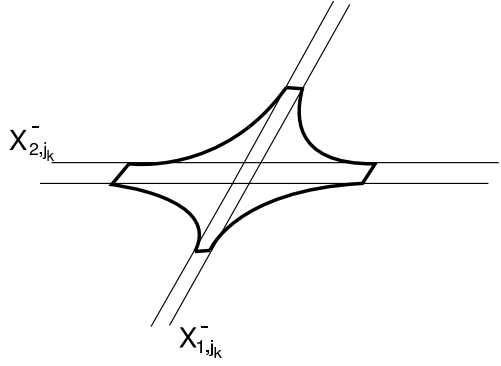}}\hspace{10mm}}
\caption{The plane region enclosed in the thickened curve is the
component of the parabolic boundary of $\bar W_{j_k}^-$ in $\p
B_\infty$.}\label{Kpb}
\end{figure}

 Now let $Hull(X_{1,j_k}\cup  X_{2,j_k})$ be the convex hull of
$X_{1,j_k}\cup  X_{2,j_k}$ in $\mathbb H^3$, and let $N_R(W_{j_k})$ be
the $R$-collared neighborhood of $W_{j_k}$ in $\mathbb H^3$. Then
 $Hull(X_{1,j_k}\cup  X_{2,j_k}) -(X_{1,j_k}\cup  X_{2,j_k}) \subset
   N_R(W_{j_k})$ by
Proposition \ref{convh}. Let $\bar W_{j_k} = N_R(W_{j_k})\cap
Hull(X_{1,j_k} \cup X_{2,j_k})$, and $\bar W_{j_k}^-=\bar
W_{j_k}\setminus {\cal B}$. We call $\bar W_{j_k}^-\cap {\cal B}$
the parabolic boundary of $\bar W_{j_k}$.
 Note that $\bar
W_{j_k}^-$ is invariant under the action of
$\g_{1,j_k}\G_1\g_{1,j_k}^{-1}\cap \g_{2,j_k}\G_2\g_{2,j_k}^{-1}$.
The component of the parabolic boundary of $\bar W_{j_k}^-$ in $\p
B_\infty$ is as shown in Figure \ref{Kpb}. Let $\bar U_{k}^-=\bar
W_{j_k}^-/
 (\g_{1,j_k}\G_1\g_{1,j_k}^{-1}\cap \g_{2,j_k}\G_2\g_{2,j_k}^{-1})$.
We now replace the $R$-collared neighborhood of $U_k^-$ in $\breve
Y^-$ by $\bar U_k^-$; that is, we glue  $\breve Y^-\setminus
N_{(R, \breve Y^-)}( U_k^-)$ with $\bar U_k^-$  along the frontier
of $N_{(R, \breve Y^-)}(U_k^-)$ in $\breve Y^-$ (which is a  part
of the boundary of $\bar U_k$), using   the original gluing map
 $\breve{g}_2 \circ h \circ \breve{g}_1^{-1}$.
 We do this operation for each
component of $U^-$. Because $U^-$ has an $R$-collared neighborhood
in $\breve Y^-$,
 the components $\bar U_k^-$ do not interfere with each other.
 That is,  if we let $Y^-$ denote the
resulting space, then  $\bar U_k^-$, $k=1,...,q$, are mutually
disjoint from each other in $Y^-$. Let $\bar U^-$ be the union of
$\bar U_k^-$, $k=1,...,q$.

\begin{figure}[!ht]
{\epsfxsize=4in \centerline{\epsfbox{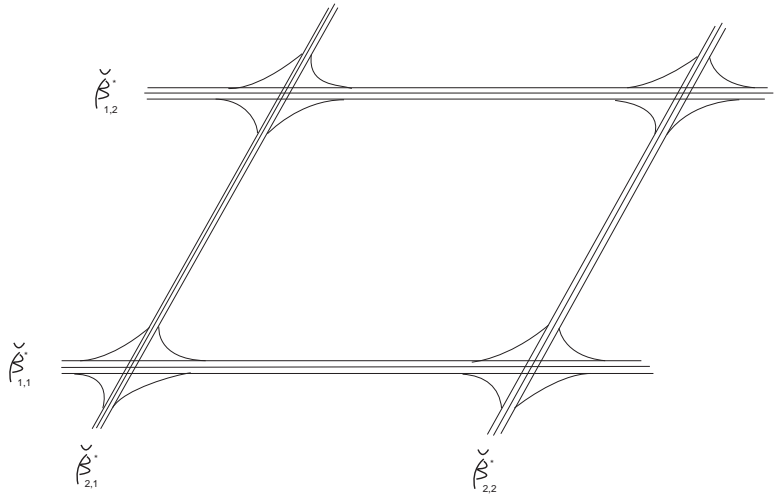}}\hspace{10mm}}
\caption{The parabolic boundary of $Y^-$.}\label{annulus}
\end{figure}

\begin{lem}\label{glue}
$Y^-$ is a connected, compact, hyperbolic $3$-manifold containing
$\breve Y_i^-$, $i=1,2$, as  submanifolds (with their
original hyperbolic structures), and there is a local isometry
 $f:Y^- \ra M^-$
 extending the maps $f_i \circ \breve{p}_i:\breve Y_i^- \ra M^-$.
\end{lem}

\pf By Corollary \ref{subbdle3},  the frontier of $N_{(R, \breve
Y_i^-)}(K_i^-)$ in $\breve Y_i^-$ is a set of (truncated) strips
and annuli (the latter set may be empty) for each $i=1,2$. Note
that the frontier of $N_{(R, \breve Y^-)} (U_k^-)$ in $\breve Y^-$
is the disjoint union of the frontier of $N_{(R, \breve
Y_i^-)}(K_i^-)$ in $\breve Y_i^-$, $i=1,2$. Hence
  $\breve Y^-$ is obtained from gluing two $3$-manifolds along
  subsurfaces in their boundaries and thus is a manifold.
Obviously it is a connected and compact $3$-manifold. We just need
to show that the hyperbolic structures of the gluing pieces
  match up over the
identified region,   forming a global hyperbolic structure on
$Y^-$.

 It is enough to verify  this around each component of $\bar
U^-$. From the construction of $\bar U_k^-$ given above,
 we see that  $X_{1,j_k}^-\cup \bar W_{j_k}^-\cup  X_{2,j_k}^-$ is
 a hyperbolic $3$-submanifold of $\mathbb H^3$.
 Also  $X_{i,j_k}^-$ is a universal cover of
$\breve Y_i^-$, $i=1,2$, and $\bar W_{j_k}^-$ is a universal cover
of $\bar U^-_k$, so there is a natural map
from $X_{1,j_k}^-\cup \bar W_{j_k}^-\cup  X_{2,j_k}^-$ to the
manifold $\breve Y_1^-\cup \bar U_k^-\cup \breve Y_2^-$.  This  provides
the required hyperbolic structure around the component $U_k^-$.

 Finally, the map $f$ can be constructed by piecing together the maps
  $f_i \circ \breve{p}_i$, and then extending to $Y^-$
 in the obvious way.
\qed

\begin{lem}\label{hbdisk}
Each component of $\bar U^-$ is a handlebody.
\end{lem}

\pf Each component $\bar U^-_k$ of $\bar U^-$ is
 homeomorphic to $K_{i,k}^-$, and thus is  compact and irreducible.
 Since the fundamental group of $\bar
U_k^-$ is isomorphic to a subgroup of the free group $\G_i$,
 then $\bar U_k^-$ is a handlebody.
 \qed

The parabolic boundary $\p_p Y^-$ of $Y^-$ in $\breve T$  is the
union of the parabolic boundary of $\breve Y_i^-$, $i=1,2$, and
that of $\bar U^-$ (see Figure \ref{annulus}).

Now we are going to construct  the cusp $C_0$ mentioned in (v)
above. The horosphere $\p B_\infty$ is a universal cover of
$\breve T$. Let $p_*:\p B_\infty\ra \breve T$ be the covering map.
Along each component of $p_*^{-1}(\{\breve \b_{i,j}, i=1,2,
j=1,...,n_i\})$ we place an appropriate translation of $X_i$ by an
element of $\G$, and at each point of $p^{-1}_*(\{s_1,...,s_d\})$
we place an appropriate translation of a component of $\{\bar
W_{j_1},...,\bar W_{j_q}\}$ by an element of $\G$. Let $Q$ denote
the union of these manifolds.

Let $B_\infty^0$ be the horoball based at $\infty$ which is
smaller than $B_\infty$ by distance one, i.e. the horizontal plane
$\p B_\infty^0$ is above $\p B_\infty$ by distance one. Let $V_0$
be the region between the two horizontal planes $\p B_\infty^0$
and $\p B_\infty$, and let $Q_0=Q\cap V_0$. Let $\tilde C_0$ be
the convex hull of $Q_0$ in $\mathbb H^3$. Then obviously $\tilde
C_0$ is contained in $B_\infty$.

\begin{lem}\label{C0}
If $n$ (and thus $N_i>n$) is large enough, then $\tilde C_0\cap
V_0=Q_0$.
\end{lem}

\begin{figure}[!ht]
{\epsfxsize=4in \centerline{\epsfbox{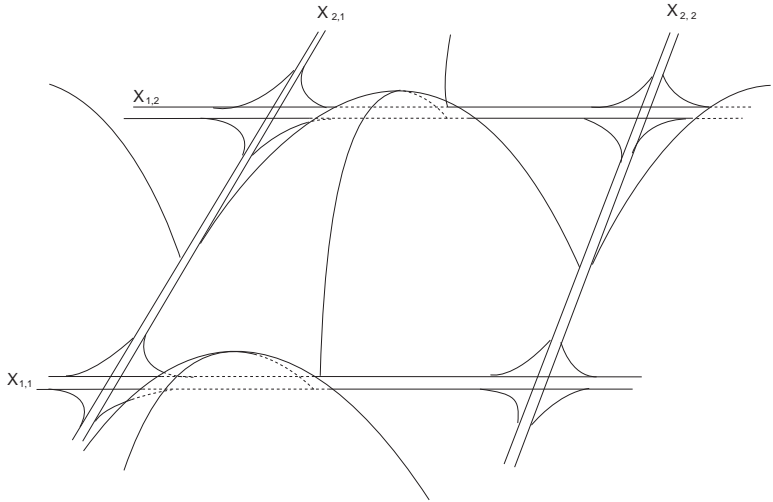}}\hspace{10mm}}
\caption{The convex hull above the plane
$\p B_\infty^0$}\label{annulus1}
\end{figure}

\pf Consider the frontier of $Q_0$ in $V_0$. It is a set of
infinitely many annuli. Let $A_1$ be one of them. Then every point
$x$ in $A_1$ is a point in the boundary of some translation of
$X_i$ or  in the boundary of some translation of  $\{\bar
W_1,...,\bar W_q\}$. The  tangent plane $P_x$ of that manifold at
$x$ (a geodesic plane) is not a vertical plane and thus its
intersection with the horizontal plane $\p B_\infty^0$ is a
Euclidean circle of finite diameter $d_x$. Modulo the action of
$\b_{1,1}^*$ and $\b_{2,1}^*$, the set $\{d_x, x\in A_1\}$ has an
upper bound independent of the integer $n$.  Also modulo the
action of the abelian group ${\cal A}=<
(\b_{1,1}^*)^{m_1},(\b_{2,1}^*)^{m_2}>$, there are only finitely
many different annuli in $Fr_{V_0}(Q_0)$. Hence  the set $\{d_x,
x\in Fr_{V_0}(Q_0)\}$ has an upper bound independent of the
integer $n$. Therefore if $n$ is sufficiently large, each $P_x$,
$x\in Fr_{V_0}(Q_0)$, will only intersect $Q_0$ at $x$.
 Thus  in forming the convex hull of $Q_0$, all the
new points added are above the plane $\p B_\infty^0$. (cf. Figure
\ref{annulus1}). The lemma is proved. \qed

We may assume that $N_i>n$ has been chosen big enough so that the
conclusion of Lemma \ref{C0} holds.

 By our construction, $Q_0$ is invariant under the action of
the abelian group ${\cal
A}=<(\b_{1,1}^*)^{m_1},(\b_{2,1}^*)^{m_2}>$, and so is $\tilde
C_0$. Now let $C_0=\tilde C_0/{\cal A}$. Then $C_0$ is contained
in the cusp $\breve C$ and $C_0\cap \breve T =\p_pY^-$. Let $Y$ be
the manifold which is  the union of $Y^-$ and $C_0$ glued along
the parabolic boundary of $Y^-$. We use the obvious gluing map,
which is locally consistent with the gluing of $Q$ and $\tilde
C_0$  in $\mathbb H^3$. As in the proof of Lemma \ref{glue}, one
can show that $Y$ is a   hyperbolic manifold with a local isometry
$f$ into $M$. Moreover $Y$ is also convex. Indeed, we only need to
check local convexity in a small neighborhood of $\p_p Y^-$ in
$Y$, which holds, since the model space $Q\cup \tilde C_0$ is
locally convex in a small neighborhood of $\tilde C_0$ in $Q\cup
\tilde C_0$.

Thus the local isometry $f$ induces an injection of $\pi_1(Y,
s_1)$ into $\pi_1(M, t_1)$.
We shall show:

\begin{pro} \label{last}
If $\D$  is bigger than one, or
 if both of $n_1$ and $n_2$ are bigger than one, then\\
(1) The boundary of $Y$ is incompressible in $Y$.\\ (2) No
essential loop in $\p Y$ is freely homotopic into $C_0$.
\end{pro}

To prove Proposition \ref{last},
 it is sufficient to show that every Dehn filling of $Y$
along its cusp $C_0$ gives a $3$-manifold with incompressible
boundary.

Let $Y(\a)$ be any  Dehn filling of $Y$ along $C_0$ with slope
$\a$. We claim that $ Y(\a)$ is an $HS$-manifold (see Section
\ref{hsmf}).
 The handlebody part
$H$ of $Y(\a)$ is $\bar U^-\cup C_0(\a)$, where $C_0(\a)$ is the
 filling of the cusp $C_0$ with slope $\a$. Indeed
 by Lemma \ref{hbdisk} each component of $\bar U^-$ is a
handlebody which connects to the solid torus $C_0(\a)$ along its
parabolic boundary $\p_p \bar U^-$ which is a set of disks. Thus
$H=\bar U^-\cup C_0(\a)$ is a connected handlebody. The $S\times
I$ part of $Y(\a)$ is $Y(\a)\setminus H=Y(\a)\setminus (\bar
U^-\cup C_0(\a))$. Indeed $Y(\a)\setminus(\bar U^-\cup C_0(\a))$
is the union of $\breve Y_i^-\setminus N_{(R, \breve
Y_i^-)}(\breve g_i(K_i^-))=\breve Y_i^-\setminus\breve
 g_i((AN_{(R,X_i)}(K_i))^-)$, $i=1,2$. It follows from Corollary \ref{subbdle3}
 that  $\breve
 g_i((AN_{(R,X_i)}(K_i))^-)$ can be considered as  $F_i^-\times I$ for some
compact subsurface $F_i^-$ of $\breve S_i^-$.
  Therefore each component of  $\breve
Y_i^-\setminus N_{(R, \breve Y_i^-)}(\breve g_i(K_i^-))$ can be
given a trivial $I$-bundle structure over a compact surface with
boundary such that the frontier in $\breve Y_i^-$
consists of $I$-fibers (these $I$-fibers may not be consistent
with the old $I$-fibers for $\breve Y_i^-$).
The surface $S$ is compact, but is possibly disconnected.

 Let $A=\p S\times I$, which is the frontier of
 $\bigcup_{i=1,2} \breve  g_i((AN_{(R,X_i)}(K_i))^-)$ in $Y(\a)$ and
 is a set of mutually disjoint, properly  embedded
 annuli in $Y(\a)$.
 By Lemma \ref{gib}, we only need to show that for each
compressing disk $D$ of $H$, $D\cap A$ has at least two
components, and that each component of $S$ is not a disk. We deal
with the latter requirement first.

\begin{lem}If $n$ (and thus $N_i>n$) is sufficiently large, then
$S$ has no disk component. \end{lem}

\pf It is equivalent to show that if $n$ is sufficiently large
then for each $i=1,2$, each component of  $\breve
Y_i^-\setminus\breve g_i((AN_{(R,X_i)}(K_i))^-)$  is not simply
connected.

Suppose otherwise that  $\breve Y_i^-\setminus\breve
g_i((AN_{(R,X_i)}(K_i))^-)$ has a component $E_0$ which is simply
connected (a $3$-ball).  We call the part of the boundary of $E_0$
which lies in $\p_p \breve Y_i^-$ the parabolic boundary of $E_0$
and denote it by $\p_p E_0$. The union of the parabolic boundary and
the frontier of $E_0$ in $\breve Y_i$ is an annulus $A_0$  in the
boundary of $E_0$. The annulus  $A_0$ can be decomposed by a set of
parallel, essential arcs into components which  are alternately
components in $Fr_{\breve Y_i^-}(E_0)$ and $\p_p E$. We call these
components \textit{frontier  faces}  and
 \textit{parabolic faces}
 of $A_0$, respectively. Since the frontier of $\breve
g_i((AN_{(R,X_i)}(K_i))^-)$ in $\breve Y_i^-$ has exactly $d$
components (Lemma \ref{d-strips}), the annulus $A_0$ has at most
$2d$ faces. Note that every parabolic face of the annulus  $A_0$
is  a very long rectangle, depending on $n$, and that every
frontier face of $A_0$ has  a bounded diameter, independent of
$n$.

The  $3$-ball component  $E_0$ has a lift, $\tilde E_0$, to
$X_i^-$, the universal cover of $\breve Y_i^-$. Note that $\tilde
E_0$ is isometric to $E_0$. Let $\tilde A_0$ be an
annulus   in the boundary of $\tilde E_0$ which is a lift of $A_0$.
The annulus $\tilde A_0$ has the corresponding decomposition into
parabolic and frontier faces. Every parabolic face  of $\tilde
A_0$ is a long Euclidean  rectangle contained in $\p_p
X_i^-=X_i\cap {\p \cal B}_i$. Since $\G$ acts transitively on
components of $ {\cal B}$, there is an element $\g$ of $\G$ such
that $\g(\tilde A_0)$ has a parabolic face $D_0$ which  lies in
$\p B_{\infty}$.

\noindent {\bf Claim}. $\g(\tilde A_0)$ has only one parabolic
face which lies in $\p B_\infty$.

Since $\g(\tilde E_0)$ is contained in $\g(X_i^-)$, we only need
to show that $\g(\tilde A_0)$ has   only one parabolic face which
lies in $\p B_\infty\cap \g(X_i^-)$, which is an  infinite
Euclidean strip between two parallel Euclidean lines.
 Note that every frontier face of $\g(\tilde A_0)$
separates $\g(X_i^-)$. It follows that if $\g(\tilde A_0)$
has at least two parabolic faces in
 $ \g(X_i^-)\cap \p B_\infty$, then there must exist a frontier face of
 $\g(\tilde A_0)$ with two opposite sides contained in the strip
 $ \g(X_i^-)\cap \p B_\infty$
 as essential arcs.
 But  this contradicts Lemma \ref{diffcusps}, proving the claim.

Recall that we have assumed that every horoball component in
${\cal B}$, except $B_\infty$, has Euclidean diameter less than one.
It follows that the Euclidean diameter of the set $\g( \tilde
A_0)\setminus D_0$ is some fixed number independent of $n$.   But
the Euclidean diameter of $D_0$ must be very large if $n$ is very
large. Thus the annulus   $\g( \tilde A_0)$ cannot exist if $n$ is
sufficiently large. The  lemma follows. \qed

We may assume that the number $N_i>n$ has been chosen big enough
so that the surface $S$ has no disk components.

 Now for the former  requirement that
 for each compressing disk $D$ of $H$, $D\cap A$ has at least
two components, it is sufficient to show that $\p H\setminus A$ is
incompressible in $H$ (since the genus of $H$ is obviously larger
than one). We show

\begin{lem}\label{d>1}
If either both of $n_1$ and $n_2$ are bigger than one or $\D$ is
bigger than one, then $\p H\setminus A$ is incompressible in $H$.
\end{lem}

\begin{figure}[!ht]
{\epsfxsize=6in \centerline{\epsfbox{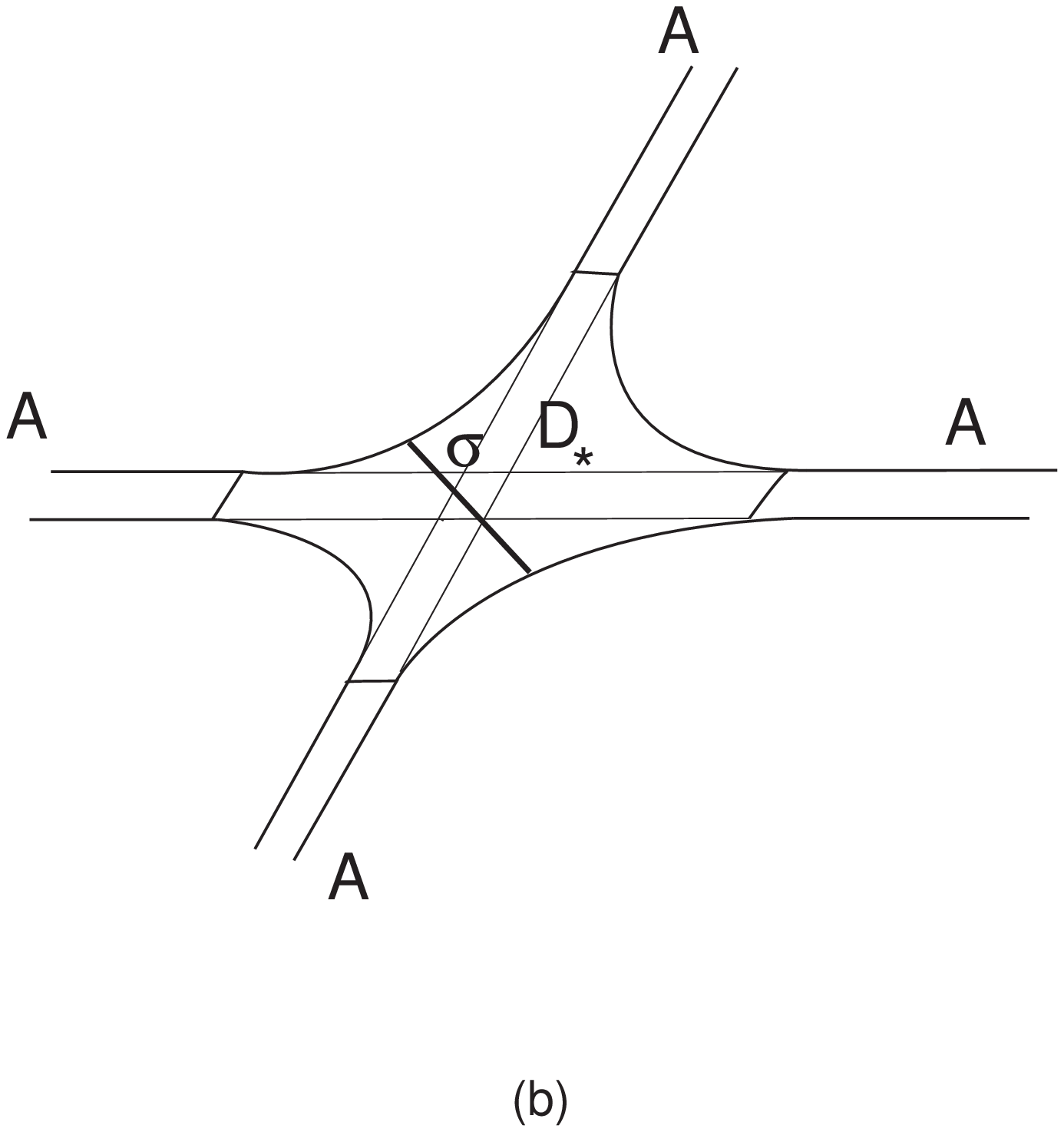}}\hspace{10mm}}
\caption{}\label{Kpb1}
\end{figure}

\pf   We call  $\p \breve Y_i^-\setminus\p_p \breve Y_i^-$ the
{\it horizontal boundary} of $\breve Y_i^-$. It has two components
and is incompressible  in $\breve Y_i^-$. The boundary of $\bar
U^-$ can be divided into three parts: the parabolic boundary $\p_p
\bar U^-$, the frontier of $\bar U^-$ in $\breve Y^-$,  and the
rest which we call the {\it horizontal boundary} of $\bar U^-$
(which we  denote by $\p_h \bar U^-$).
 Figure \ref{Kpb1} illustrates $\p_p \bar U^-$;
 in this figure, the frontier boundary meets $\p_p \bar U^-$ in
 straight segments, and the horizontal boundary meets $\p_p \bar U^-$
 in curved arcs.

\noindent {\bf Claim}. The horizontal boundary of $\bar U^-$ is
incompressible in $\bar U^-$.

 \textit{Proof of Claim:}  We just need to prove the claim for each
 component $\bar U_k^-$ of $\bar U^-$.
 First note that
  the boundary of the $I$-bundle $\breve g_i((AN_{(R,X_i)}(K_{i,k}))^-)$
can be naturally divided into parabolic, frontier and horizontal
boundaries as well.
 Let $A_{i,k}$
 be the frontier boundary of $\breve g_i((AN_{(R,X_i)}(K_{i,k}))^-)$,
 and let $S_{i,k}'$ be the horizontal boundary
 of $\breve g_i((AN_{(R,X_i)}(K_{i,k}))^-)$. Note that $A_{1,k}\cup A_{2,k}$
 is the frontier boundary of $\bar U_k^-$  and that
 $\p_h \bar U_k^- = \p \bar U_k^-\setminus (\p_p \bar U_k^- \cup A_{1,k} \cup A_{2,k})$. Obviously
$S_{i,k}'$ is incompressible in $\breve
g_i((AN_{(R,X_i)}(K_{i,k}))^-)$.  Since  each component of
$S_{i,k}'$ separates
 $\bar U_k^-$ and carries the fundamental group of $\bar U_k^-$,
 each component of  $ \p \bar U_k^-
       \setminus (\p_p \bar U_k^- \cup A_{1,k})$ is parallel
    in $\bar U_k^-$    to a component of $S_{i,k}'$ and thus is
    incompressible in $\bar U_k^-$.
 The components of $A_{2,k}$ are all annuli and strips,
 and the core curve of every annulus component
 of $A_{2,k}$ is essential in $\bar U_k^-$.  Therefore,
  $\p_h \bar U_k^- = (\p \bar U_k^-
       \setminus (\p_p \bar U_k^- \cup A_{1,k})) \setminus A_{2,k}$
 is incompressible in $\bar U_k^-$.
The proof of the claim is finished.

 Returning to the proof of Lemma \ref{d>1},
 suppose that there is a compressing disk $D$
for $H$  which is disjoint from the annuli $A$. We may assume that $D$ is
chosen to minimize   the components of  $D\cap \p_p \bar U^-$.

 If $D\cap \p_p \bar U^-$ is empty, then  $D$ is contained in $\bar
U^-$ (it cannot be in $C_0(\a)$ since  $\p C_0(\a)\setminus
\p_pY^-$ is a set of disks), contradicting the claim.
 Thus $D\cap \p_p \bar U^- \neq \emptyset$. Certainly we may
assume that $D\cap \p_p \bar U^-$ has no circle components. Let
 $\s$ be an arc component of $D\cap \p_p \bar U^-$ which is outermost
in $D$. The arc $\s$ divides $D$ into two disks; let $D_0$ be the
one whose interior  is disjoint from $\p_p \bar U^-$. Let $\b=\p
D_0\cap \p D$. Then $\p D_0=\s\cup\b$. Let $D_*$ be the component of
$\p_p\bar U^-$ which contains the arc $\s$.

 Figure \ref{Kpb1} shows the
parabolic boundary of $Y^-$ near $D_*$. A pair of parallel
straight lines in the figure (including the dotted line segments)
is a part of a pair of circles which bounds a component of the
original parabolic boundary of $\breve Y_i^-$. There are two such
components at $D_*$, one from $\p_p\breve Y_1^-$ and the other
from $\p_p\breve Y_2^-$. We call the components of their intersections
 with $\p D_*$ {\it corners} of $D_*$. Alternately, the four corners
are the intersection components of the annuli $A$ with $D_*$.

We claim that the endpoints of $\s$ cannot separate the four corners
in $\p D_*$, i.e. a case like that shown in Figure \ref{Kpb1} (b) or
(c) is impossible. Indeed,  the endpoints of $\s$ are also the
endpoints of the connected arc $\b$ whose interior is disjoint from
the parabolic boundary of $Y^-$ and the annuli $A$. So if a case
like Figure \ref{Kpb1} (b) or (c) happens, then $\b$ cannot be
contained in $\p C_0(\a)$. For otherwise
 the geometric intersection number $\D$ would be one and
 $n_1$ or $n_2$ would be equal to one.
 The arc $\b$ cannot be contained in the horizontal
boundary of $\bar U^-$ either. For the endpoints of $\s$ lies in
different components of the horizontal boundary of $\bar U^-$.

 Hence $\s$ is contained in $D_*$
as shown in Figure \ref{Kpb1} (a). Let $\b'$ be the sub-arc in $\p
D_*$ which is disjoint from the corners of $D_*$ and  co-bounds a
sub-disk $D_1$ in $D_*$ with $\s$. Then the union of $D_0$ and $D_1$
along $\s$ is a properly embedded disk in  $H$ which we denote by
$D_2$. Suppose that $\b$ is contained in $\p C_0(\a)$. Then $D_0$ is
contained in $C_0(\a)$. Since $\p C_0(\a)\setminus \p_p Y^-$ is a
set of disks,  $\p D_0$ cannot be an essential curve in the torus
$\p C_0(\a)$. Thus $\p D_0$ bounds a disk $D_3$ in $\p C_0(\a)$. The
two disks $D_0$ and $D_3$ form a $2$-sphere in $C_0(\a)$ and thus
bound a $3$-ball in $C_0(\a)$ (since $C_0(\a)$ is irreducible).  Now
it is clear that we can isotope the part of $D$ contained in the
$3$-ball to cross the sub-disk $D_1$ of $D_*$ and thus reduce the
number of intersection components of $D\cap \p_p \bar U^-$. Suppose
then that $\b$ is contained in the horizontal boundary of $\bar
U^-$. Then $D_0$ is contained in $\bar U^-$ and so is the disk
$D_2$. Since the horizontal boundary is incompressible, $\p D_2$ is
not an essential curve in the horizontal boundary, i.e. $\p D_2$
must bound a disk $D_4$ in the horizontal boundary.  The two disks
$D_2$ and $D_4$ form a $2$-sphere in $\bar U^-$ and thus bound a
$3$-ball in $\bar U^- $ (since $\bar U^-$ is irreducible). Again
 we can isotope the part of $D$  contained in
the $3$-ball to cross the sub-disk $D_1$ of $D_*$ and thus reduce
the number of intersection components of $D\cap \p_p \bar U^-$.
\qed

The proof of Proposition \ref{last} is  finished.

 We now are in position to finish the proof of Theorem \ref{qfs}.
 Obviously $Y$ has non-empty boundary.
 Suppose $\D$  is bigger than one, or
 that both of $n_1$ and $n_2$ are bigger than one.
 Then we claim that
 $f|_{\p Y}$ is a quasi-Fuchsian surface.
 Indeed, since $f$ is injective on $\pi_1 Y$,  Part (1) of Proposition \ref{last}
 implies that $f$ is injective in $\pi_1 \p Y$.
 Since $\p Y$ is closed, then $f|_{\p Y}$ is not a virtual fiber.
 Therefore, by Marden/Thurston/Bonahon's
 classification of essential surfaces (see introduction),
 it is enough to show that $f^* \pi_1 \p Y$ contains no non-trivial
 parabolic elements.

 The torus  $\p C_0$ is  incompressible in
$Y$ (otherwise $Y$ would be an open solid torus,
 which is obviously impossible). Hence $f^*(\pi_1(C_0,
s_1))$ is a finite index subgroup of the abelian group $\pi_1(\p
C, t_1)$. Hence if $\a$ is a non-trivial loop of $\p Y$,
 and if $f \a$ is freely homotopic into $C$, then some
non-zero power of $\a$  is freely homotopic into $C_0$,
 contradicting Proposition \ref{last} Part (2).

 Suppose then, that $\D = 1$ and  one of $n_1$ or $n_2$ (say $n_1$)
 is 1. In this case we take the double cover of the manifold $Y$
dual to the non-separating surface $\breve S_1$ in $Y$ (note that
$\breve S_1$ is naturally embedded in $Y$). Let $\hat p:\hat Y\ra
Y$ be the double cover. Then $\hat Y$ is a convex hyperbolic
3-manifold with a single cusp, which maps by a local isometry into $M$.
  Also $\hat p^{-1}(\breve Y_1^-)$
has two components, and so in particular its parabolic boundary
has two components on the boundary of the cusp $\hat C_0=\hat
p^{-1}(C_0)$. Now we just need to show that every Dehn filling  of
$\hat Y$ gives a manifold whose boundary is incompressible. Let
$\hat Y(\a)$ be any  Dehn filling of $\hat Y$ along the cusp $\hat
C_0$ with slope $\a$. We give  $ \hat Y(\a)$ the obvious
$HS$-manifold structure. Obviously the surface cross interval part
of the $HS$-manifold has no simply connected components (since
taking the double cover does not change this property). The rest
of proof is exactly as that of Proposition \ref{last} since the
parabolic boundary of $\hat p^{-1}(\breve Y_i^-)$ has at least two
components now for each of $i=1,2$.
 This completes the proof of Theorem \ref{qfs}.

\vspace{10mm}

\noindent Mathematics Department \\ SUNY at Buffalo \\Buffalo, NY
14290
\\
jdmaster@buffalo.edu

\noindent
Mathematics Department \\ SUNY at Buffalo \\Buffalo, NY
14290
\\
xinzhang@buffalo.edu

\end{document}